\documentclass[11pt]{article}
\usepackage[affil-it]{authblk}

\usepackage[margin=1in]{geometry} 
\usepackage{color}
\definecolor{clemson-orange}{RGB}{234,106,32}
\definecolor{chicago-maroon}{RGB}{128,0,0}
\definecolor{northwestern-purple}{RGB}{82,0,99}
\definecolor{lime-green}{RGB}{50,205,50}
\definecolor{pink}{RGB}{255,0,128}

\usepackage{graphicx}

\usepackage{algorithm}
\usepackage{algpseudocode}
\usepackage{xfrac}
\usepackage{caption}

\RequirePackage[OT1]{fontenc} 
\RequirePackage{graphicx}
\RequirePackage{amsthm}
\RequirePackage[cmex10]{amsmath}
\RequirePackage{natbib}
\RequirePackage{hypernat}
\usepackage{fullpage}
\usepackage{mathrsfs}

\allowdisplaybreaks
\usepackage{setspace}
\definecolor{pink}{RGB}{255,0,128}

\definecolor{cornell-red}{RGB}{179,27,27}

\newcommand{\bfp}{\mathbf{p}}
\newcommand{\bfx}{\mathbf{x}}
\newcommand{\bfmu}{\boldsymbol\mu}
\newcommand{\bfnu}{\boldsymbol\nu}

\usepackage[title]{appendix}%
\usepackage{amsmath}
\usepackage{amssymb}
\usepackage{amsthm}
\usepackage{enumitem}
\setlist[enumerate]{noitemsep, topsep=0pt}

\bibpunct[, ]{(}{)}{,}{a}{}{,}%
\newcommand{\bb}{\mathbb}

\newcommand{\N}{{\bb N}}

\newcommand{\R}{\bb R}
\newcommand{\E}{\bb E}

\newcommand{\cT}{\mathcal{T}}
\newcommand{\calT}{\mathcal{T}}
\newcommand{\calP}{\mathcal{P}}
\newcommand{\calL}{\mathcal{L}}

\newcommand{\reals}{\mathbb{R}}
\newcommand{\mubar}{\bar{\mu}}
\newcommand{\bfphi}{\boldsymbol\phi}
\newcommand{\bfv}{\mathbf{v}}

\DeclareMathOperator{\argmin}{argmin}

\usepackage[small,compact]{titlesec}
\titleformat{\subsection}[runin]
        {\normalfont\bfseries}
        {\thesubsection}
        {5pt}
        {}
        [.]
\titleformat{\subsubsection}[runin]
        {\normalfont\bfseries}
        {\thesubsubsection}
        {0.5em}
        {}
        [.]
\titleformat{\section}
        {\large\bfseries}
        {\thesection}
        {5pt}
        {}

\usepackage{titling}
\usepackage{etoolbox}
\newcommand{\zerodisplayskips}{%
  \setlength{\abovedisplayskip}{2pt}%
  \setlength{\belowdisplayskip}{2pt}%
  \setlength{\abovedisplayshortskip}{2pt}%
  \setlength{\belowdisplayshortskip}{2pt}}
\appto{\normalsize}{\zerodisplayskips}
\appto{\small}{\zerodisplayskips}
\appto{\footnotesize}{\zerodisplayskips}

\theoremstyle{definition}
\newtheorem{theorem}{Theorem}
\newtheorem{lemma}{Lemma}
\newtheorem{corollary}{Corollary}
\newtheorem{proposition}{Proposition}
\newtheorem{definition}{Definition}
\newtheorem{remark}{Remark}
\newtheorem{assumption}{Assumption}
\newtheorem{example}{Example}

\newtheorem{problem}{Problem}

\usepackage[colorlinks,citecolor=northwestern-purple,urlcolor=chicago-maroon,linkcolor=chicago-maroon]{hyperref}

\usepackage[nameinlink]{cleveref}
\crefname{assumption}{Assumption}{Assumptions}
\crefname{lemma}{Lemma}{Lemmas}
\crefname{theorem}{Theorem}{Theorems}
\crefname{corollary}{Corollary}{Corollaries}
\crefname{proposition}{Proposition}{Propositions}
\crefname{claim}{Claim}{Claims}
\crefname{procedure}{Procedure}{Procedures}
\crefname{algorithm}{Algorithm}{Algorithms}
\crefname{figure}{Figure}{Figures}
\crefname{remark}{Remark}{Remarks}
\crefname{section}{Section}{Sections}
\crefname{procedure}{Procedure}{Procedures}
\crefname{definition}{Definition}{Definitions}
\crefname{example}{Example}{Examples}
\crefname{table}{Table}{Tables}
\crefname{equation}{}{}
\crefname{enumi}{}{}
\crefname{problem}{Problem}{Problems}
\crefname{appendix}{Appendix}{Appendices}

\usepackage{subcaption}

\usepackage{mathtools}
\usepackage{bbm}
 
\title{Wasserstein gradient flow for optimal probability measure decomposition}
\author{%
Jiangze Han,\thanks{University of British Columbia, email: jiangze.han@sauder.ubc.ca} \ Christopher Thomas Ryan,\thanks{University of British Columbia, email: chris.ryan@sauder.ubc.ca} \ Xin T. Tong\thanks{National University of Singapore, email: xin.t.tong@nus.edu.sg}
}
\date{May 16, 2024}

\begin{document}
\maketitle

\abstract{We examine the infinite-dimensional optimization problem of finding a decomposition of a probability measure into $K$ probability sub-measures to minimize specific loss functions inspired by applications in clustering and user grouping. We analytically explore the structures of the support of optimal sub-measures and introduce algorithms based on Wasserstein gradient flow, demonstrating their convergence. Numerical results illustrate the implementability of our algorithms and provide further insights.




\section{Introduction}
With the rapid advancement of AI, automated algorithms are increasingly being used to solve routine problems. Particularly intriguing are the applications of AI in social organizations, which have the potential to benefit both private and public sectors. These applications include the organization of markets, allocation of resources, and mechanism design, among others \citep{dai2021learning,chen2021market,niazadeh2023online,zhalechian2022online,agrawal2023learning}. This paper studies a new problem of how to decompose a population of customers or clients into groups to optimize a generic quantitive criterion. 

Consider the following probability measure decomposition problem. Later, we will show how this problem can arise in applications. Individuals in a population are represented by their feature vectors $\bfx\in\mathbb{R}^d$. Feature vectors are distributed according to a probability distribution $\pi$. Let $\calP_2(\R^d)$ be the space of probability measures defined on $\R^d$ with finite second moment; that is, $\calP_2(\R^d)=\{\mu:\int_{\R^d} \|\bfx\|^2d\mu(\bfx)<\infty\}$. When a measure $\mu\in\calP_2(\R^d)$ is absolutely continuous with respect to the Lebesgue measure, we use the same symbol $\mu$ to represent the measure's associated probability density function. We define a decomposition of $\pi$ as follows.

\begin{definition}[Probability measure decomposition]
    Given a probability measure $\pi\in\mathcal{P}_2(\R^d)$, we say the vector  $\bfmu\doteq(\mu_1, \mu_2,\ldots,\mu_K)\in\mathcal{P}_2(\R^d)^{\otimes K}$ of probability measures with weight vector $\bfp=(p_1,\ldots, p_K)\in\mathbb{R}^{\otimes K}$ is a \emph{decomposition} of $\pi$, if $(\bfmu,\bfp)\in\mathcal{P}_\pi$, where
        \begin{equation}\label{eq:feasible-region}
            \calP_\pi\doteq\left\{(\bfmu,\bfp):\sum_{k\in[K]} p_k=1,p_k\geq 0, \sum_{k\in[K]} p_k \mu_k=\pi\right\}
        \end{equation}
        and $[K] \doteq \{1,2,\dots, K\}$. 
        The equality $\sum_{k\in[K]}p_k\mu_k=\pi$ holds in duality with the space $C^\infty_c(\R^d)$ of compactly supported smooth (i.e., has infinitely many derivatives) functions; that is, for all $f\in C^\infty_c(\R^d)$,
        \begin{equation*}
           \sum_{k\in[K]}\int_{\R^d} f(x)d\mu_k(x)=\int_{\R^d} f(x)d\pi(x).
        \end{equation*}
            
\end{definition}

Intuitively, we decompose the population (distributed according to $\pi$) of feature vectors into $K$ sub-populations (distributed according to $\mu_1, \dots, \mu_K$). 
Within each sub-population $k\in[K]$, individual features are distributed according to probability measure $\mu_k$. The population weight of the whole population is normalized to be $1$. Each sub-population $k\in[K]$ has weight $p_k$. 

Among all decompositions of the feature distribution $\pi$, we seek one that minimizes (i) a weighted sum of distribution loss function $L:\mathcal{P}_2(\R^d)\rightarrow\mathbb{R}$ associated with feature distribution $\mu_k$ of each sub-population, and (ii) a weight loss function $R:\mathbb{R}\rightarrow\mathbb{R}$
associated with the population weight $p_k$ of each sub-population. The purpose of this loss is to penalize a sub-population with a small weight, which can be impractical for different reasons detailed in the examples below.  Formally, we consider the following optimal decomposition problem. 

\begin{problem}[Optimal decomposition problem]\label{prob:opt-decomposition}
Let $L:\mathcal{P}_2(\R^d)\rightarrow\mathbb{R}$ be a distribution loss function and $R:\mathbb{R}\rightarrow\mathbb{R}$ a weight loss function. Given a target of $K$ sub-populations in a population with distribution $\pi$, solve
     \begin{equation}\label{eqn:genR}
\min_{(\bfmu,\mathbf{p})\in\calP_\pi}F(\bfmu,\bfp),\quad F(\bfmu,\bfp)\doteq\sum_{k\in[K]} (p_k L(\mu_k)+R(p_k)),
\end{equation}
where the feasible region $\mathcal P_\pi$ is defined in \cref{eq:feasible-region}.
\end{problem}

We consider the following family of distribution loss functions $L$.
\begin{definition}[Coupled loss function]\label{def:loss}
    We say $L:\mathcal{P}_2(\R^d)\rightarrow\mathbb{R}$ is a \emph{coupled loss function} if 
        \begin{equation*}
            L(\mu)=\int_{\R^d} \int_{\R^d} \ell(\mathbf x,\mathbf y)d\mu(\mathbf x)d\mu(\mathbf y)
        \end{equation*}
        for some continuously differentiable function $\ell:\R^d\times\R^d\rightarrow\mathbb{R}$ satisfying
    \begin{equation*}
        |\ell(\mathbf z,\mathbf x)-\ell(\mathbf z,\mathbf y)|\leq\|\mathbf x-\mathbf y\|
    \end{equation*}
for all $\mathbf x, \mathbf y,\mathbf z \in \R^d$. We call $\ell$ the \emph{kernel} of $L$.
\end{definition}

\begin{definition}[Weight loss function]
For some $\theta,\beta>0$, define $R:(0,1)\rightarrow\R$ as
        \begin{equation*}
            R(p) \doteq \frac{\theta}{p^\beta}.
        \end{equation*}
\end{definition}

We present two applications of this general setup.

\begin{example}[League design with Elo rating system \citep{elo1978rating}]\label{exp:elo}
In many competition-based online games, players are grouped into different ``leagues" based on their skill levels, and only players from the same league can compete with each other. League design aims to create competitive gaming environments where players are not overwhelmed by strong opponents or bored by weaker ones. One way to quantify the skill level and competitiveness of games is the Elo-type system.%
\footnote{For descriptions about Elo rating system, please see \url{https://en.wikipedia.org/wiki/Elo_rating_system}.}
For simplicity, we focus on one-on-one competitions, similar to chess. In the Elo-type system, each player is given a skill level $x \in (0,\infty)$ (sometimes called Elo score). The probability of winning for a player with skill level $x$ against a player with skill level $y$ is taken to be $\sfrac{x}{(x+y)}$.%
\footnote{In other variants of the Elo system, people use $\sfrac{\log(x)}{\alpha}$ to represent skill level for some game specific parameter $\alpha$, which is equivalent to our setting by a change-of-variable argument.}
%
A game is deemed more competitive as each player's win rate gets closer to $50\%$. A common practice is to minimize the difference of each player's winning probability with $50\%$. For example, \cite{simonov2023suspense} show in their study using data from the game streaming platform Twitch that the expected game length and viewership can be increased by making the round-win probabilities of games closer to a balanced distribution of $50\%$-$50\%$.

Suppose the skill-level distribution of all players has density $\pi$, and the goal is to decompose players into $K$ leagues. 
Suppose in each league $k\in[K]$, players arrive to join a game according to a Poisson process with arrival rate $p_k$ (i.e., the expected waiting time for players in sub-population $k$ is $\sfrac{1}{p_k}$). Our decomposition aims to maximize the competitiveness of games and minimize the waiting time of each league. This can be achieved by solving \cref{eqn:genR} with weight loss function $R(p)=\sfrac{1}{p}$ and distribution loss function
    \begin{equation}\label{eq:elo-loss}
        \begin{split}
                L(\mu)&=\int_{\R^d}\int_{\R^d} \left(\frac{x}{x+y}-\frac{1}{2}\right)^2d\mu(x)d\mu(y)\\
            &=\frac12\int_{\R^d}\int_{\R^d} \frac{x^2+y^2}{(x+y)^2}d\mu(x)d\mu(y)-\frac14
        \end{split}
    \end{equation}
Note that $L$ is a coupled loss function with kernel $\ell(x,y)=\frac{1}{2}\frac{x^2+y^2}{(x+y)^2}-\frac{1}{4}$. We call this distribution loss function $L$ the $\textit{Elo loss}$.
\end{example}

\begin{example}[Generalized clustering]\label{exp:clustering}
In clustering, the goal is to generate sub-populations according to specific criteria. As in our base setup, suppose a population's feature vectors $\mathbf x\in\R^d$ are distributed according to $\pi\in\calP_2(\R^d)$. This means that the feature vector $X$ of a randomly sampled individual in the population is a random variable with law $\pi$. Suppose a designer aims to decompose this population into sub-populations to maximize a sense of similarity in certain feature dimensions while concurrently maximizing a sense of diversity in other dimensions within each sub-population. Accordingly, we can define a loss function $L$ as follows. Let $W$ be a diagonal matrix with nonzero diagonal entries. Define a distribution loss $L$ by
    \begin{equation}\label{eq:var-loss}
        \begin{split}
                 L(\mu)&=\int_{\R^d}\langle  \mathbf x-\E_{\mu}[X], W(\mathbf x-\E_{\mu}[X])\rangle d\mu(\mathbf x)\\
                 &=\int_{\R^d}\int_{\R^d} \langle \mathbf x-\mathbf y, W(\mathbf x-\mathbf y )\rangle d\mu(\mathbf x)d\mu(\mathbf y),
        \end{split}
    \end{equation}
where $\langle \mathbf x,\mathbf y\rangle$ denotes the standard inner product in $\R^d$. Note that $L$ is a coupled loss function with the kernel $\ell(\mathbf x,\mathbf y)=\langle \mathbf x-\mathbf y, W(\mathbf x-\mathbf y)\rangle$. We call this distribution loss function $L$ the $\textit{variance loss}$. 

The matrix $W$ specifies the weights assigned to each dimension. When $W_{i,i}>0$, this minimizes the dissimilarity of features in dimension $i$. Conversely, when $W_{i,i}<0$, this maximizes the diversity of features in dimension $i$. Notably, if $W$ is the identity matrix, then $L(\mu)$ corresponds to the trace of the covariance matrix. In this special case, we decompose the distribution $\pi$ into $K$ sub-distributions $\mu_k$ to minimize the variance of each sub-distribution $\mu_k$.

To create sub-populations with sufficiently large sizes, we can also impose penalties for selecting smaller sub-populations. This can be achieved, for example, by setting the weight loss function $R$ to $R(p)=\sfrac 1p$.
\end{example}

 In certain scenarios, the population size $p_k$ of each sub-population $k\in[K]$ is predefined. For example, to make the decomposition more balanced, one can require $p_k=\sfrac1K$ for all $k\in[K]$.  In this scenario, we can define the set of feasible decompositions for a given $\mathbf p$: 
    \begin{equation}\label{eq:feasible-region-given-p}
        \calP_{\pi,\bfp}
        \doteq\{\bfmu\in\mathcal{P}_2(\R^d)^{\otimes K}: \sum_{k\in[K]} p_k \mu_k=\pi\}. 
    \end{equation}
In this case, we consider the following simplified problem.
\begin{problem}[Optimal decomposition problem with specified weights]\label{prob:decomposition-2}
Let $L:\mathcal{P}_2(\R^d)\rightarrow\mathbb{R}$ be a given distribution loss function. Given the number $K$ of sub-populations, sub-population weights $\bfp=(p_1,\ldots,p_K)$, and feature distribution $\pi$, solve
  \begin{equation}
\label{eqn:fixpR}
\min_{\bfmu\in\calP_{\pi,\bfp}}F_{\bfp}(\bfmu),\quad F_{\bfp}(\bfmu)\doteq\sum_{k\in[K]} p_k L(\mu_k).
\end{equation}
\end{problem}
\cref{prob:decomposition-2} can be thought of a sub-problem of \cref{prob:opt-decomposition}. Recall that we denote the objective functions in \cref{prob:opt-decomposition} and \cref{prob:decomposition-2} by $F(\bfmu,\bfp)$ and $F_\bfp(\bfmu)$, respectively. \cref{prob:opt-decomposition} can be decomposed as follows
    \begin{equation*}
        \min_{(\bfmu,\mathbf{p})\in\calP_{\pi}} F(\bfmu,\bfp)=\min_{\mathbf{p}\in \mathcal{S}} \min_{\bfmu\in \calP_{\pi,\bfp}}F_\bfp(\bfmu),
    \end{equation*}
where $\mathcal{S}=\{\bfp:\sum_{k\in[K]} p_k=1,p_k\geq 0\}$. The inner problem in the above decomposition is exactly \cref{prob:decomposition-2}.  

\subsection{Literature review}\label{ss:literature}
In recent years, distribution-oriented optimization has been an active direction in operations research. Its applications span over robust optimization, matching of social networks, and stochastic games; see for examples \cite{bertsimas2019adaptive,hu2022graph,light2022mean}.   Our paper falls in this general research direction. Our discussion relies on the methodologies of optimal transport theory and Wasserstein gradient flow. Optimal transport theory introduces a metric---known as the Wasserstein metric---on the set of probability measures, effectively transforming it into a metric space called Wasserstein space. This metric enables the development of calculus and geometric concepts, such as gradient and geodesics, on the space of probability measures. The standard references for optimal transport are \cite{villani2009optimal,santambrogio2015optimal}. With calculus and geometry defined by optimal transport theory, we can establish the concept of gradient flow on the Wasserstein space. Wasserstein gradient flow has emerged as a popular methodology for tackling optimization problems formulated on the space of probability measures. The standard references are \cite{ambrosio2005gradient,santambrogio2017euclidean}. Various optimization algorithms leveraging Wasserstein gradient flow have been developed to address diverse optimization problems \citep{salim2020wasserstein,chewi2020gradient,chizat2018global,zhang2018policy}. In contrast to prior research, one of our methodological contributions is introducing a Wasserstein gradient flow tailored to a constrained optimization problem. We demonstrate that, in the limit, the proposed Wasserstein gradient flow converges to a feasible solution satisfying the optimality condition.

We also give two applications of our paper in \cref{exp:clustering,exp:elo}. In \cref{exp:elo} and \cref{ss:case-study-elo}, we show that our model and algorithm can be applied to address the league design problem. It is a common practice for online games to group players into different groups (called ``league" in our paper) \citep{agarwal2009matchmaking,manweiler2011switchboard,francillette2013players}. One important consideration for grouping is to form groups where players share similar skill levels while ensuring that sessions have an adequate number of players \citep{francillette2013players}, which is the major setting we consider. 

Another closely related problem is clustering, where $N$ data points, $\bfx_1,\ldots, \bfx_N$, need to be split into $K$ groups. The most well known version is $K$-means clustering problem, which seeks the optimal label $l: [N]\to [K]$ that minimizes
\[
L_K(l)\doteq  \sum_{k\in [N]} \|\bfx_{k}-\mathbf m_{l(k)}\|^2,\quad \mathbf m_{i}\doteq \frac{\sum_{k\in [K]} 1_{l(k)=i}\bfx_{k}}{\sum_{k\in [K]} 1_{l(k)=i}}. 
\]
Suppose we view each cluster as a decomposition component. In that case, K-means clustering is similar to the problem in \cref{exp:clustering} with $W=I_d$ and $R=0$, where $I_d$ is a $d$ by $d$ identity matrix. The main difference is that K-means clustering focuses on handling problems where finitely many data samples are present, while our decomposition problem focuses on the case where the data distribution is available. Given its popularity, algorithmic studies of K-means began in the 1960s \citep{macqueen1967some} and remain active until today. Interested readers can refer to a recent survey \cite{ikotun2023k} for its development. Due to its combinatorial nature, algorithms for K-means clustering often rely on greedy heuristic arguments. To the best of our knowledge, the associated optimality conditions of these algorithms are seldom available, except for continuous relaxation of K-means \citep{blomer2020complexity}. This is also a difference between our problem and the classical K-means clustering problem, as we will formulate and prove the associated optimality conditions.   

Our paper investigates the problem of optimally decomposing a probability measure into a finite number of probability measures based on a specific objective function. In statistics, a related problem is the identifiability of a finite mixture model. A finite mixture model can be seen as a probability density function that is a convex combination of $K$ probability density functions called ``component densities". A mixture model is said to be identifiable if this convex combination can be uniquely determined. For literature in this direction, see \cite{teicher1963identifiability,yakowitz1968identifiability,kim2015empirical,mclachlan2019finite}. Our research problem differs in two key aspects. First, our problem focuses on finding the optimal decomposition of a probability measure that minimizes a specific loss functional among all possible decompositions, whereas ``identifiability" concerns the uniqueness of the decomposition. Second, the ``identifiability" literature mainly focuses on parametric probability densities such as Gaussian and Gamma distributions, whereas our problem deals with nonparametric probability measures.

The main algorithm we proposed is a constraint controlled gradient flow (CCGF) type of algorithm. The notion of a CCGF was first proposed in \cite{liu2021sampling} to generate samples with moment constraints. The decomposition problems we study here imposes a very different kind of constraint, which has not been studied before in the literature. Moreover, we provide a construction of CCGF (more details in \cref{rem:CCGFequality}) and explain it in the Euclidean space setting.

\subsection{Our contributions}

In short, this work's main contributions are threefold. First, we formulate the optimal decomposition problems and discuss the geometric properties of the optimal solution. Second, we develop the constrained Wasserstein optimality condition for the optimal decomposition problem. Third, we design gradient flow-based algorithms that can approximately achieve the optimality conditions and numerically test their efficacy on various problems. To the best of our knowledge, these aspects have not been studied in the literature. 

The paper is organized as follows. In \cref{s:structure}, we prove several structural properties of the optimal solutions to \cref{prob:opt-decomposition,prob:decomposition-2}. In \cref{s:preliminary}, we introduce the necessary technical ingredients to understand our analysis, including gradient flow, calculus in Wasserstein space, and geodesic geometry of Wasserstein space. In \cref{s:optimality-condition}, we show that existing optimality conditions do not readily apply to our setting, and we propose our own optimality conditions. In \cref{s:algorithms}, we develop Wasserstein gradient flow procedures to solve \cref{prob:opt-decomposition,prob:decomposition-2} and show their convergence. In \cref{s:implementation}, we implement  Wasserstein gradient flow in the settings described by \cref{exp:elo,exp:clustering}. \cref{s:conclusion} concludes and points to future research directions. Proofs of all technical results can be found in the appendix.

\section{Structural results and interpretation}\label{s:structure}
Before transitioning to our gradient flow design, we take some time to reflect on the type of structures of the solutions to \cref{prob:opt-decomposition,prob:decomposition-2} we hope to explore. In particular, we are interested in the geometric properties we might see arise in optimality. We will see these geometric properties clearly in numerical examples in \cref{s:implementation}.

We first note that $\bfmu$ defined by $\mu_k=\pi$ for all $k\in[K]$ is a trivial feasible solution to \cref{prob:opt-decomposition,prob:decomposition-2}. We do not want this ``boring" solution in practice. Hence, a natural question is: to what extent, and in what sense, are the optimal sub-population measures $\mu^*_k$ different? For simplicity, we focus on the optimal solution $\bfmu^*$ to \cref{prob:decomposition-2} in this section, although the results also hold for the optimal solutions to \cref{prob:opt-decomposition}.

We answer this question by analyzing a generalized notion of the ``support'' (i.e., $(\delta,c)$-interior support explained later) of an optimal solution. First, we check if the ``supports" of the $\mu_k$ intersect or are disjoint. If the ``supports" are indeed disjoint, we can already rule out the trivial solution $\mu_k=\pi$ for all $k\in[K]$. Moreover, in this case, when we decompose the probability density $\pi$ in \cref{prob:decomposition-2}. We are effectively partitioning the underlying space $\R^d$. 

We are also interested in the relative positions of these ``supports" in $\R^d$. For example, in the one-dimensional case, is it possible that the ``support" of $\mu_1$ is $[0,1]\cup[2,3]$ and the ``support" of $\mu_2$ is $(1,2)$? Motivated by this example, we say a collection of sets $\{S_k\}_{k\in[K]}$, where $S_k\subseteq\R^d$ for $k\in[K]$, is \textit{convex in pairs} if $\text{conv}(S_i)\cap S_j=\emptyset$ for any pair $i\neq j\in[K]$. In the previous example, the ``support" $[0,1]\cup[2,3]$ of $\mu_1$ and the ``support" $(1,2)$ of $\mu_2$ are not convex in pairs. This concept speaks to the practical implementability of the results that arise from our approach. For instance, in \cref{exp:elo}, a league design that is not convex in pairs would be unnatural because it implies grouping high-end and low-end players into one league while placing middle-level players in another. This type of design would be hard for game designers and players to justify. 



We formally define our notion of $(\delta,c)$-interior support as follows.

\begin{definition}[$(\delta,c)$-interior support]\label{def:support}
For real numbers $\delta,c>0$, $S(\delta,c)\subseteq \R^d$ is the $(\delta,c)$-\emph{interior support} of a probability density $\mu\in\mathcal{P}_{2,ac }(\R^d)$ if for any $\mathbf x\in S(\delta,c)$, 
\begin{equation*}
    \mu(\mathbf y)>c,\quad \forall \mathbf y\in\R^d\text{ s.t. } \|\mathbf x-\mathbf y\|\leq \delta.
\end{equation*}
\end{definition}

Recall that the \emph{support} of a probability density function $\mu:\R^d\rightarrow\R_+$ is the set $S\doteq\{\mathbf x\in\R^d:\mu(\mathbf x)>0\}$. Note that $S(\delta,c)\subseteq S$ for any $\delta,c>0$. Our \cref{def:support} of $S(\delta,c)$ generalizes the support $S$ of a density function in the following two ways. First, in $S(\delta,c)$ we only consider points with a density larger than a threshold $c$. We use this threshold to determine if the density $\mu$ is too small. A point $\mathbf x$ is called $c$-\textit{negligible} if $\mu(\bfx)\leq c$. Points in the support $S$ are $0$-negligible while points in $(\delta,c)$-interior supports are $c$-negligible. 

Second, a point $\bfx\in\R^d$ is in $S(\delta,c)$ if all points $\mathbf y$ in the $\delta$-neighbourhood are not negligible (i.e., $\mu(\mathbf y)>c$). In comparison, a point $\bfx\in\R^d$ is in the support $S$ as long as this point itself is not negligible. Hence, the support can be seen as a $(\delta,c)$-interior support with $\delta=0$ and $c=0$.

Let $S_k(\delta,c)$ be the $(\delta,c)$-interior support of the optimal solution $\mu_k$ to \cref{prob:decomposition-2} with the distribution loss functions presented in \cref{exp:elo,exp:clustering}. In the rest of this section, we check under what conditions the $\{S_k(\delta,c)\}_{k\in[K]}$ are disjoint and convex in pairs. We first show that $\{S_k(\delta,c)\}_{k\in[K]}$ are disjoint under mild conditions.
\begin{proposition}\label{prop:disjoint-support}
      Let $\{\mu_k\}_{k\in [K]}$ be the optimal solution to \cref{prob:decomposition-2}. Suppose the distribution loss function $L$ is a coupled loss function (\cref{def:loss}) with kernel $\ell:\R^d\times\R^d\to\R$ and let $\delta > 0$ and $c > 0$ be given. Suppose $\nabla^2_{x,y} \ell(\mathbf x_0,\mathbf x_0)$ is not positive semi-definite for any $\mathbf x_0\in\R^d$. Then, $S_i(\delta,c)\cap S_j(\delta,c)=\emptyset$ for any $i\neq j\in[K]$.
\end{proposition}




\begin{corollary}\label{cor:disjoint-support}
      Suppose $L$ is either Elo loss (equation \cref{eq:elo-loss}) or variance loss (equation \cref{eq:var-loss}). Then, the $(\delta,c)$-interior supports of the optimal densities to \cref{prob:decomposition-2} are disjoint for any $\delta,c>0$; that is, $S_i(\delta, c)\cap S_j(\delta,c)=\emptyset$ for any $\delta,c>0$ and any $i,j\in[K]$ such that $i\neq j$.
\end{corollary}

Hence, when we decompose $\pi$ into $\mu_1,\ldots,\mu_K$, we in fact almost partition (note that $\cup_{k=1}^KS_k(\delta,c)$ may not be $\R^d$) the feature space $\R^d$ into $K$ disjoint sets.

Next, we show in \cref{exp:non-convex} below that under Elo loss, $\{S_k(\delta,c)\}_{k\in[K]}$ may not be convex in pairs; that is, for any two points $\mathbf x,\mathbf y\in S_i(\delta,c)$, there are possibly some points in the line segement $\{\lambda \mathbf x+(1-\lambda)\mathbf y:\lambda\in(0,1)\}$ that is not in $S_i(\delta,c)$ but in $S_j(\delta,c)$ for some $j\neq i$. Note that this also implies $S_i(\delta,c)$ is not convex. 

\begin{example}[Non-convexity under Elo loss]\label{exp:non-convex}
We model the application described in \cref{exp:elo} as an instance of \cref{prob:decomposition-2}. Let $K=2$ and $p_1=p_2=\sfrac{1}{2}$. That is, we group players into two leagues with equal sub-populations. Consider the Elo loss function $L$ (equation \cref{eq:elo-loss}). Suppose the skill level distribution $\pi$ is given by  
    \begin{equation*}
        \pi = \frac 14\delta_{a-1}+\frac 12\delta_a+\frac 14\delta_{a+1},\quad a>1,
    \end{equation*}
    where $\delta_x$ is the dirac measure at point $x$. That is, only $3$ possible skill levels are given by $a-1,a$, and $a+1$. The optimal solution is $\mu^*_1=\delta_{a}$ (with support $S_1\doteq\{a\}$) and $\mu^*_2=\tfrac{1}{2}\delta_{a-1}+\tfrac{1}{2}\delta_{a+1}$ (with support $S_2\doteq\{a-1,a+1\}$). Note that $\tfrac{1}{2}\mu^*_1+\tfrac{1}{2}\mu^*_2=\pi$. Hence, $\tfrac{1}{2}(a-1)+\tfrac{1}{2}(a+1)\in S_1$, while $a-1,a+1\in S_2$. The proof is presented in the appendix. 
\end{example}

Note that the discrete distribution $\pi$ in the example does not admit a probability density function, and hence, not a $(\delta, c)$-interior support. However, the same outcome can be shown if we smooth the discrete distribution to get a density function. Our numerical experiments (\cref{exp:mixedlognormal}) illustrate this using a mixed lognormal distribution with three peaks. In \cref{exp:mixedlognormal} below we show that under a certain distribution $\pi$, it is indeed better to group high-end players and low-end players into one league while middle-level players in another league. 

We show in the following lemma that, under a mild condition, $\{S_k(\delta,c)\}_{k\in[K]}$ are convex in pairs if the distribution loss in \cref{prob:decomposition-2} is the variance loss \cref{eq:var-loss}.

\begin{lemma}\label{lem:variance-convex}
Suppose $L(\mu)$ is the distribution loss (equation \cref{eq:var-loss}) in \cref{exp:clustering} with a positive semidefinite $W$. Given any $\delta,c>0$, the collection of sets $\{S_k(\delta,c)\}_{k\in[K]}$ are convex in pairs.
\end{lemma}

Note that \cref{lem:variance-convex} does not imply that $S_i(\delta,c)$ is convex. Some points in the line segement $\{\lambda \mathbf x+(1-\lambda)\mathbf y:\lambda\in(0,1)\}$ for some $\mathbf x,\mathbf y\in S_i(\delta,c)$ may not be in $S_i(\delta,c)$. This may happen if the support of $\pi$ is not convex by itself. However, \cref{lem:variance-convex} shows that at least the line segment joining any two points $\mathbf x,\mathbf y\in S_i(\delta,c)$ is not in any other $S_j(\delta,c)$. 

\section{Preliminaries on optimal transport and gradient flow}\label{s:preliminary}

Let $\calP_{2,ac}(\R^d)$ be all probability measures that are absolutely continuous with respect to the Lebesgue measure $m$ and have a finite second moment, i.e. 
\[
\calP_{2,ac}(\R^d)=\left\{\mu: \mu\ll m,\quad \int_{\R^d} \|\mathbf x\|^2 d\mu(\mathbf x)<\infty\right\}.
\]
Let $(\reals^d,\Sigma)$ be the Boreal measurable space. Given a measurable mapping $T: (\reals^d,\Sigma)\to \reals^d$ and a measure $\mu:\Sigma\to[0,\infty]$, the pushforward measure $T\sharp\mu:\Sigma\to[0,\infty]$ of $\mu$ under $T$ is defined to be the measure given by $T\sharp\mu(B)=\mu(T^{-1}(B))$ for $B\in\Sigma$. The main property we need for the push-forward measure is the change-of-variables formula; that is, $\E_{T\sharp \mu}[f]=\E_{\mu}[f\circ T]$ for any measurable function $f:\R^d\rightarrow\R$ and where $\circ$ denotes function composition.

This paper assumes that $\R^d$ has the usual topology. Let $\calL^2(\R^d;\R^d)$ be the $L_2$ space of functions $f:(\R^d,\mathcal{B})\to(\R^d,\mathcal{B},m)$,  where $\mathcal{B}$ is the Borel $\sigma$-algebra and $m$ is the Lebesgue measure; that is,  $\calL^2\doteq \calL^2(\R^d, \R^d)\doteq\{f:\R^d\to\R^d:\int_{\R^d}\|f(\mathbf x)\|^2dm(\mathbf x)<\infty\}$, where $\|\cdot\|$ is the $2$-norm in $\R^d$.

Let $\calL^2_{\mu}(\R^d;\R^d)$ be the $L_2$ space of functions $f:(\R^d,\mathcal{B})\rightarrow(\R^d,\mathcal{B}, \mu)$; that is, 
    \[
    \calL^2_\mu(\R^d;\R^d) \doteq\{f:\R^d\rightarrow\R^d:\int_{\R^d}\|f(\mathbf x)\|^2d\mu(\mathbf x)<\infty\}.
    \]
The inner product in $\calL^2_\mu\doteq \calL^2_\mu(\R^d;\R^d)$ is defined as $\langle f, g\rangle_\mu\doteq\int_{\R^d}\langle f(\mathbf x),g(\mathbf x)\rangle d\mu(\mathbf x)$ for any $f,g\in\calL^2_\mu$. The $\calL^2_\mu$ norm is defined as $\|f\|_\mu\doteq\sqrt{\int_{\R^d}\|f(\mathbf x)\|^2d\mu(\mathbf x)}$ for any function $f\in\calL^2_\mu(\R^d;\R^d)$. 

We also need the following product space: for a vector $\bfmu=(\mu_1,\ldots,\mu_K)$ of Borel measures,
    \begin{equation*}
    \calL^2_{\bfmu}(\R^d;\R^d) \doteq     \prod_{k\in[K]}\calL^2_{\mu_k}(\R^d;\R^d)=\{(f_1,\ldots,f_K):f_k:(\R^d,\mathcal{B})\rightarrow(\R^d,\mathcal{B}, \mu_k)\in\calL^2_{\mu_k},\forall k\in[K]\}.
    \end{equation*}
The inner product in $\calL^2_{\bfmu}(\R^d;\R^d)$ is $\langle f, g\rangle_{\bfmu}\doteq\sum_{k\in[K]}\int_{\R^d}\langle f(\mathbf x),g(\mathbf x)\rangle d\mu_k(\mathbf x)$ for any $f,g\in\calL^2_{\bfmu}$. The corresponding norm is given by $\|(f_1,\ldots,f_K)\|_{\bfmu}=\sqrt{\sum_{k\in[K]}\|f_k\|^2_{\mu_k}}$. 

\subsection{Preliminaries on gradient flow}

Standard references on gradient flow include \cite{ambrosio2005gradient,santambrogio2015optimal,santambrogio2017euclidean}. Our treatment will be contained, only setting the necessary development to make sense of our algorithms. 

Given a time-changing velocity field $\phi:[0,1]\times\R^d\rightarrow\R^d$, suppose a particle is moving according to this velocity field $\phi(t,\mathbf x)$; that is, at time $t$ and position $\mathbf x$, the velocity of this particle is given by $\phi(t,\mathbf x)$. Then, the particle's trajectory is the solution to the ordinary differential equation (ODE): $\dot{\mathbf x}(t)=\phi(t,\mathbf x(t))$. In particular, if $\phi(t,\mathbf x)$ is the gradient of some function $F:\R^d\rightarrow\R$, then the solution function $\mathbf x:[0,1]\rightarrow\R^d$ is called the \emph{gradient flow} of $F$. 

Now, suppose a population of particles, whose positions in $\R^d$ are initially distributed according to probability measure $\mu(0)$, are moving together according to the velocity field $\phi(t,\mathbf x)$; that is, the trajectory $\mathbf x(t)$ of a particle initially positioned at point $\mathbf x_0$ is given by the solution to 
 \begin{equation}\label{eq:trajectory-ODE}
         \dot{\mathbf x}(t)= \phi(t,\mathbf x(t)),\quad \bfx(0)=\mathbf x_0.
 \end{equation}
Under this velocity field, $\phi(t,\mathbf x)$, particles at different initial positions $\mathbf x_0$ have different trajectories. Hence, we can define a mapping $T_\phi^t: \mathbf x(0)\mapsto \mathbf x(t)$, which maps the initial position $\mathbf x(0)$ of a particle to its position $\mathbf x(t)$ at time $t$. Here, $\mathbf x(t)$ is the solution to ODE \cref{eq:trajectory-ODE}.

With this mapping $T^t_\phi$, we can also study the change in the distribution of the particles' positions. Since these particles are initially distributed according to $\mu(0)$, and their positions at time $t$ are given by the map $T^t_\phi$, their positions at time $t$ are distributed according to the pushforward measure $\mu(t)=T^t_{\phi}\sharp\mu(0)$. It turns out, as a curve in the space of probability measures, $\{\mu(t) : t \in [0,1]\}$ can be characterized by the following partial differential equation (PDE) called the \emph{continuity equation}, which is expressed as follows:
    \begin{equation}\label{eq:continuity-equation}
        \frac{d}{dt}\mu(t,\mathbf x)=-\nabla \cdot (\phi(t,\mathbf x)\mu(t,\mathbf x)),
    \end{equation}
where $\nabla\cdot:\R^d\rightarrow\R$ denotes the divergence operator and where we may think of $\mu(t,\mathbf x)$ as a density function evaluated at the point $\mathbf x$. This equation should be interpreted with some caution. Each term in the continuity equation \cref{eq:continuity-equation} has no meaning in isolation since each term is not even defined if $\mu(t)$ is a measure (and not a density function). The continuity equation only holds in the following weak sense; that is, for any test function $f\in C_c^\infty(\R^d)$, 
    \begin{equation}\label{eq:weak-sol}
         \frac{d}{dt}\int_{\R^d}f(\mathbf x)d\mu(t,\mathbf x)=\int_{\R^d}\langle\nabla f(\mathbf x),\phi(t,\mathbf x)\rangle d\mu(t,\mathbf x),
    \end{equation}
     where $C^\infty_c(\reals^d)$ are smooth functions with compact support. We note that this equation \cref{eq:weak-sol} is well-defined even if $\mu(t)$ is a measure, where here $\mu(t,\mathbf x)$ is short-hand for $\mu(t)(\mathbf x)$. We formalize this argument in the following lemma.

\begin{lemma}[Theorem 1.3.17 in \citep{chewilog}]\label{lem:implementation}
    Let $\phi(t,\mathbf x):[0,1]\times\R^d\rightarrow\R^d$ be a time-dependent velocity field. Suppose random variables $t\mapsto X_t$ evolve according to $\dot{X}_t=\phi(t,X_t)$. Then, the law $\mu_t$ of $X_t$ solves the continuity equation \cref{eq:continuity-equation} in the weak sense; that is, equation \cref{eq:continuity-equation} holds for $(\mu_t,\phi(t,\cdot))$ in the sense of equation \cref{eq:weak-sol}. 
\end{lemma}

\subsection{Calculus in Wasserstein space}\label{ss:calculus-W}

To develop gradient flow on our space $\calP_2(\R^d)$ of probability measures, we must define an appropriate notion of ``gradient". The standard definition is inadequate since it is only defined in vector spaces, but $\calP_2(\R^d)$ is not a vector space. We define our notion of gradient over the space of probability measures using optimal transport theory, as we explain now.

Although a large part of optimal transport theory can be developed in a more general framework, we focus on the $\calP_2(\R^d)$ for simplicity. The \emph{Wasserstein-2 distance} on $\calP_2(\R^d)$ is defined as follows. 
\begin{definition}[Wasserstein-2 distance]
For $\mu,\nu\in\calP_2(\R^d)$,
    \begin{equation}\label{eq:Wasserstein}
        W^2_2(\mu,\nu)\doteq\min\left\{\int_{\R^d\times\R^d}\|\mathbf x-\mathbf y\|^2d\gamma:\gamma\in\Pi(\mu,\nu)\right\},
    \end{equation}
    where $\|\cdot\|$ is the standard Euclidean norm and $\Pi(\mu, \nu)$ is the set of \emph{transport plans} defined as
    \begin{equation}
        \Pi(\mu, \nu)\doteq\{\gamma\in\calP(\R^d\times\R^d):(\pi_{\mathbf x})\sharp\gamma=\mu,(\pi_{\mathbf y})\sharp\gamma=\nu\},
    \end{equation}
    where $\pi_{\mathbf x}$ and $\pi_{\mathbf y}$ are the two coordinate projections of $\R^d\times\R^d$ onto $\R^d$. The minimizer $\gamma^*$ to \cref{eq:Wasserstein} is called the \emph{optiaml transport plan}. Moreover, if there exists a measurable map $T:\R^d\to\R^d$ such that $\gamma^*= (\text{id}, T)\sharp\mu$, then this map $T$ is called the \emph{optimal transport map}. 
\end{definition}

The set $\calP_2(\R^d)$ is called \emph{Wasserstein space} when endowed with the Wasserstein-2 metric $W_2$. It is well-known that $(\calP_2(\R^d), W_2)$ is a metric space but not a vector space since the sum of two probability measures need not be a probability measure (see 
 Proposition 5.1 in \cite{santambrogio2015optimal}). 

In Euclidean space, the gradient of a differentiable function $f:\R^d\to\R$ is defined to be the unique vector field $\nabla f:\R^d\to\R^d$ such that, for any $\mathbf x_0,\mathbf v\in\R^d$, 
 \begin{equation}\label{eq:directional-gradient}
        \lim_{t\rightarrow 0}\frac{f(x_0+t\mathbf v)-f(\mathbf x_0)}{t}=\langle\nabla f(\mathbf x_0),\mathbf v\rangle.
    \end{equation}
The left-hand side is the directional derivative of function $f$ at the point $\mathbf x_0$ along direction $\mathbf v$. The right-hand side is the vector product of gradient $\nabla f(\mathbf x_0)$ at point $\mathbf x_0$ and the direction $\mathbf v$. That is, in Euclidean space, the instantaneous change of function value (i.e., left-hand side of equation \cref{eq:directional-gradient}) along any direction $\mathbf v$ can be easily computed by the inner product between gradient and the direction (i.e., right-hand side of equation \cref{eq:directional-gradient}). 

To define the Wasserstein gradient, we first define a similar notion of directional derivative in Wasserstein space. In Euclidean space, the term $\mathbf x(t)\doteq \mathbf x_0+t\mathbf v$ in the directional derivative \cref{eq:directional-gradient} defines a trajectory of movement $\mathbf x(t):[0,1]\rightarrow\R^d$ that starts at $\mathbf x_0$ and moves in direction $\mathbf v$. In Wasserstein space, the vector addition $\mathbf x_0+t\mathbf v$ is not valid since Wasserstein space is not a vector space. Instead of using vector addition, we use absolutely continuous curves $\mu(t):[0,1]\rightarrow (\calP_2(\R^d), W_2)$ to describe the trajectory of movement. This ``absolutely continuous curve" is defined using the Wasserstein metric. We can define a notion of ``direction'' of the curve $\mu(t)$ similarly to how $\mathbf v$ defines a direction for $\mathbf x(t)$.

The direction of the trajectory $\mu(t)$ can be described as follows. As shown in the following \cref{lem:ac-curve}, in Wasserstein space, any absolutely continuous curve can be characterized by the continuity equation \cref{eq:continuity-equation}. Moreover, there exists a vector field that can be understood as the ``direction" of the trajectory $\mu(t)$.

\begin{lemma}[Theorem 8.3.1 \citep{ambrosio2005gradient}]\label{lem:ac-curve}
    Let $\mu(t):(0,1)\rightarrow(\calP_2(\R^d),W_2)$ be an absolutely continuous curve. Then there exists a vector field $\phi:(0,1)\times\R^d\rightarrow\R^d$ such that the continuity equation \cref{eq:continuity-equation} holds in the weak sense, and for $m$-a.e. $t\in (0,1)$, 
       \[
      \phi(t, \cdot)\in  \overline{\{\nabla \Psi: \Psi\in C^\infty_c(\reals^d)\}}^{\calL^2_{\mu_t}},
        \]
        where $\bar A^{\calL_{\mu_t}^2}$ represents the topological closure of set $A$ in $\calL_{\mu_t}^2\doteq \{f:\R^d\rightarrow\R^d: \int_{\R^d}\|f(\mathbf x)\|^2d\mu_t(\mathbf x)<\infty \}$.
\end{lemma}

In \cref{lem:ac-curve}, at any $t\in(0,1)$, the vector field $v(t):\R^d\to\R^d$ is called the ``tangent vector" to the curve $\mu(t)$ at time $t$ because it describes the direction of trajectory $\mu(t)$ through the continuity equation \cref{eq:continuity-equation}. \cref{lem:ac-curve} indicates that starting at any point $\mu_0\in(\calP_2(\R^d),W_2)$, all directions that we can move along in Wasserstein space should be contained in the following set
 \[
 \calT \calP_2(\mu_0)=\overline{\{\nabla \Psi: \Psi\in C^\infty_c(\reals^d)\}}^{\calL^2_{\mu_0}}.
 \]
 The set $\calT \calP_2(\mu_0)$ is called the tangent space, which is known to be a vector space \citep[Lemma 8.4.2]{ambrosio2005gradient}. With this tangent space, the Wasserstein-2 gradient can be defined through equation \cref{eq:directional-gradient} as follows.
 
\begin{definition}[Wasserstein-2 gradient]
Given a functional $F:\calP_2(\R^d)\rightarrow\R$, the Wasserstein gradient of $F$ at $\mu_0\in\calP_2(\R^d)$ is defined to be the unique element $\nabla F(\mu_0)\in\calT\calP_2(\mu_0)$ such that, for every curve $\mu(\cdot):[0,\infty)\rightarrow\calP_2(\R^d)$ with $\mu(0)=\mu_0$ and tangent vector $v_0\in \calT\calP_2(\mu_0)$ at time $0$, it holds that   
    \begin{equation}\label{eq:Wasserstein-directional}
          \frac{d}{dt}F(\mu(t))|_{t=0}=\langle\nabla F(\mu_0),v_0\rangle_{\mu_0},
    \end{equation}
where $ \langle \nabla F(\mu_0), v_0\rangle_{\mu_0}\doteq \int_{\R^d} \langle\nabla F(\mu_0), v_0\rangle d\mu_0$.
\end{definition}

Note that equation \cref{eq:Wasserstein-directional} recovers the nice property described in \cref{eq:directional-gradient} for Euclidean space. Next, we compute the Wasserstein-2 gradient of the following relevant functions.

\begin{lemma}
\label{lem:U}
The Wasserstein gradient $\nabla L:\R^d\rightarrow\R^d\in \calL^2_{\mu}$ of the coupled loss function $L$ (defined in \cref{def:loss}) at $\mu\in\calP_2(\R^d)$ is
\[
\nabla L(\mu)(\mathbf x)=\int_{\R^d} (\nabla_1 \ell(\mathbf x,\mathbf z)+\nabla_2 \ell(\mathbf z,\mathbf x))d\mu(\mathbf z).
\]
\end{lemma}

In our work, we are optimally deciding a vector $\bfmu=(\mu_1,\mu_2,\ldots, \mu_K)\in \calP_2(\R^d)^{\otimes K}$ of probability measures. It is straightforward to generalize the notion of Wasserstein calculus onto $\calP_2(\R^d)^{\otimes K}$. In particular, the tangent space at $\bfmu$ is
\[
 \calT \calP^{\otimes K}_2(\bfmu)=\bigotimes_{k\in[K]}\overline{\{\nabla \Psi: \Psi\in C^\infty_c(\reals^d)\}}^{\calL^2_{\mu_k}}.
 \]
 
 Given a functional $F:\calP_2(\R^d)^{\otimes K}\rightarrow\R$,
we say it is differentiable at $\bfmu_0$ with gradient  
$\nabla F(\bfmu_0)=(\nabla_{\mu_1} F(\bfmu_0),\ldots, \nabla_{\mu_K}F(\bfmu_0))$  if for every curve $\bfmu(\cdot):[0,\infty)\rightarrow\calP^{\otimes K}_2$ with $\bfmu(0)=\bfmu_0$ and tangent vector $\bfv_0\in \calT\calP^{\otimes K}_2(\bfmu_0)$ at time $0$, it holds that   
    \begin{equation}\label{eq:Wasserstein-gradient-K}
          \frac{d}{dt}F(\bfmu(t))|_{t=0}=\langle\nabla F(\bfmu_0),\bfv_0\rangle_{\bfmu_0},
    \end{equation}
where $ \langle \nabla F(\bfmu_0), \bfv_0\rangle_{\bfmu_0}\doteq \sum_{k\in[K]} \langle\nabla_{\mu_k} F(\bfmu_0), (\bfv_0)_k\rangle_{(\bfmu_0)_k}$.

By the definition of objective function $F_{\bfp}(\bfmu)$ of \cref{prob:decomposition-2}  and $F(\bfmu,\bfp)$ of \cref{prob:opt-decomposition}, the Wasserstein gradient $\nabla_{\bfmu}F_{\bfp}(\bfmu_0)$ and $\nabla_{\bfmu} F(\bfmu_0,\bfp_0)$ at $\bfmu_0\in\calP^{\otimes K}_2$ and $\bfp_0$ are given by 
    \[
   \nabla_{\mu_k} F(\bfmu_0,\bfp_0)= \nabla_{\mu_k}F_{\bfp}(\bfmu_0)=p_k\nabla L(\mu_k).
    \]

In later development, we will have occasion to use Kullback-Leibler (KL) divergence to measure constraint violations in our problem. Therefore, we refresh the readers' understanding of this concept and derive its Wasserstein gradient. 

Recall that the Kullback-Leibler (KL) divergence between two probability densities $\mu$ and $\pi$ is defined as 
    \begin{equation*}
        \text{KL}(\mu\|\pi)\doteq\int_{\R^d}\mu(x)\log\frac{\mu(x)}{\pi(x)}dx.
    \end{equation*}
KL divergence is a type of statistical distance of how one probability density $\mu$ differs from a second, reference probability density $\pi$. It is $0$ if and only if $\mu=\pi$ almost surely.

\begin{lemma}
\label{lem:KLder}
For any $\bfmu\in\calP^{\otimes K}_{2,ac}(\R^d)$ and $\bfp\in\R^d_+$ with $\sum_{k\in[K]}p_k=1$, define $\mubar\doteq\sum_{k\in[K]}p_k\mu_k$. For any $\pi\in\calP_{2,ac}(\R^d)$, the Wasserstein gradient $\nabla_{\bfmu}\text{KL}(\mubar\|\pi)=(\nabla_{\mu_k}\text{KL}(\mubar\|\pi))_{k\in[K]}$ is given by 
\[
\nabla_{\mu_k} \text{KL}(\bar{\mu}\|\pi)(x)=p_k (s_{\mubar}(x)-s_{\pi}(x)),
\]
for all $x\in\R^d$, where $s_{\mubar}(x)\doteq \nabla\log\mubar(x)$ and $s_{\pi}(x)\doteq \nabla\log\pi(x)$.
\end{lemma}

\subsection{Geodesic geometry of Wasserstein space}
Our work also needs a concept analogous to the vector-space notion of a ``convex set" but defined in Wasserstein space. Recall that in Euclidean space, a set $M$ is convex if for any $\mathbf x,\mathbf y\in M$, the whole line segment $(1-t)\mathbf x+t\mathbf y$ for $t\in[0,1]$ connecting $\mathbf x$ and $\mathbf y$ lies in the set $M$. Since Wasserstein space is not a vector space, we cannot define the ``line segment" by the vector addition $(1-t)\mathbf x+t\mathbf y$. Nevertheless, we can define an analogous notion of ``line segment" called a ``geodesic."

In a nutshell, a geodesic $\mu(t):[0,1]\rightarrow\calP_2(\R^d)$ in Wasserstein space is the shortest path connecting $\mu(0)\in\calP_2(\R^d)$ and $\mu(1)\in\calP_2(\R^d)$. The term ``shortest path" refers to the fact that, if we measure the ``total length" of this curve $\mu(t)$ with Wasserstein metric, it is the shortest one among all curves connecting $\mu(0)$ and $\mu(1)$. This generalizes the fact that a line segment is the shortest path connecting two points in Euclidean space. In Wasserstein space, geodesics are characterized as follows. 

\begin{lemma}[Theorem 5.27 \citep{santambrogio2015optimal}]
    Suppose $T:\R^d\rightarrow\R^d$ is the optimal transport map (to problem \cref{eq:Wasserstein}) from $\mu\in\calP_{2,ac}(\R^d)$ to $\nu\in\calP_{2}(\R^d)$. Then the curve $\mu(t):[0,1]\rightarrow\calP_{2}(\R^d)$ defined by $\mu(t)\doteq ((1-t)\text{id}+tT)\sharp\mu$ is a geodesic connecting $\mu$ and $\nu$, where $\text{id}:\R^d\rightarrow\R^d$ is the identity map, i.e., $\text{id}(\mathbf x)=\mathbf x$.
\end{lemma}

A set $M\subseteq\calP_2(\R^d)$ is \emph{geodesically convex} in Wasserstein space if, for any $\mu,\nu\in M$, the geodesic $\{\mu(t) : t \in [0,1]\}$ connecting $\mu$ and $\nu$ is contained in the set $M$.

\section{Optimality condition}\label{s:optimality-condition}

In this section, we derive optimality conditions for \cref{prob:opt-decomposition,prob:decomposition-2}. First, we note below that the classical Karush-Kuhn-Tucker (KKT) conditions do not apply in our setting. Accordingly, we propose two new optimality conditions (i.e., \cref{prop:opt-condition-1,prop:opt-condition-2}) for \cref{prob:opt-decomposition,prob:decomposition-2}. The proposed optimality conditions generalize the geometric form of the KKT condition in Euclidean space.

We begin with \cref{prob:decomposition-2}. Recall that the feasible set of \cref{prob:decomposition-2} is 
\[
\calP_{\pi,\bfp}=\left\{(\mu_1,\ldots, \mu_K)\in \calP_2(\R^d)^{\otimes K}:\sum_{k\in [K]} p_k \mu_k=\pi\right\}.
\]
Note that $\calP_{\pi,\bfp}$ is not geodesically convex in Wasserstein space $(\calP_2(\R^d), W_2)$, even though the constraint looks very linear (but not linear in Wasserstein space!). To see this, consider the following simple example.

\begin{example}[Non-convexity of the feasible set $\calP_{\pi,\bfp}$]
Let $K=2$, $\bfp=(\frac12,\frac12),$ and $\pi=\text{Uniform}[-1,1]$. We have two simple decompositions $\bfmu=(\mu_1,\mu_2)\in \calP_{\pi,\bfp}$ and $\bfnu=(\nu_1,\nu_2)\in \calP_{\pi,\bfp}$ given by 
\[
\mu_1=\text{Uniform}[-1,0],\quad \mu_2=\text{Uniform}[0,1],\quad\text{and}
\]
\[
\nu_1=\text{Uniform}[0,1],\quad \nu_2=\text{Uniform}[-1,0].
\]
The optimal transport map $T_1:\R^d\rightarrow\R^d$ from $\mu_1$ to $\nu_1$ is given by $T_1(x)=x+1$, and the optimal transport map $T_2:\R^d\rightarrow\R^d$ from $\mu_2$  to $\nu_2$ is $T_2(x)=x-1$. Hence, the geodesic $(\mu_1,\mu_2)(\cdot):[0,1]\rightarrow\calP^{\otimes 2}_2$ connecting $\bfmu$ and $\bfnu$ are given by 
\[
\mu_1(t)=\text{Uniform}[-1+t,t],\quad \mu_2(t)=\text{Uniform}[-t,1-t]. 
\]
It is easy to check that $(\mu_1(t),\mu_2(t))\in \calP_{\pi,\bfp}$ only when $t=0$ or $1$. Hence, the feasible set $\calP_{\pi,\bfp}$ is not geodesically convex in Wasserstein space.  
\end{example}

Despite this unusual geometry, we can still consider a KKT-type optimality condition on $\calP_{\pi,\bfp}$. Our optimality condition is motivated by the following observation about the usual KKT condition in Euclidean space. 

Suppose we seek to minimize a smooth function $f:\R^d\rightarrow\R$ subject to constraint $g(\mathbf x)=0$, where $g:\R^d\rightarrow\R$ is also a smooth function. The standard KKT condition is given by 
\begin{equation}
\label{tmp:KKT} 
\nabla f(\mathbf x)+\lambda \nabla g(\mathbf x)=0,
\end{equation} for some $\lambda\in \R$. We call equation \cref{tmp:KKT} the algebraic form of the KKT condition. The algebraic form can be extended to optimization problems in infinite-dimensional vector spaces \citep[Theorem 1, Section 9.3]{luenberger1997optimization}. However, this version does not apply to our setting since Wasserstein space is not a vector space. Hence, we need to propose our own KKT condition.

The standard KKT condition in Euclidean space also has a geometric interpretation due to the manifold structure of the feasible region, which we now describe. The constraint $g(\mathbf x)=0$ defines a smooth manifold $M\doteq g^{-1}(0)\cap\{\bfx\in\R^d:\nabla g(\bfx)\neq 0\}$. For each point $\mathbf x$ on $M$, the manifold $M$ has a \textit{tangent space} $\calT M(\mathbf x)$ which is a vector space consisting of all tangent vectors at $\bfx$ to differentiable curves on $M$ passing through $\bfx$. Moreover, the tangent space can be charaterized by $\calT M(\mathbf x)=\{\mathbf v\in\R^d: \mathbf v^\top \nabla g(\mathbf x)=0\}$.  For details of above discussion about manifold structure, see \cite[Theorem 5.4, 5.5]{edwards1994advanced}. From a geometric perspective, equation \cref{tmp:KKT} represents the fact that, at optimal solution $\bfx$, $\nabla f(\bfx)$ is in the straight line $\{\mathbf v\in\R^d:\mathbf v=\lambda\nabla g(\bfx)\text{ for some $\lambda\in\R$}\}$ defined by $\nabla g(\bfx)$. To avoid the need for a parameter $\lambda$ in the optimality condition, equivalently, we can say the projection of $\nabla f(\bfx)$ onto the space $\calT M(\bfx)$, which consists of vector $\bfv$ that is orthogonal to $\nabla g(\bfx)$, is $0$. We introduce the projection norm onto $\calT M(\bfx)$ to quantify this ``orthogonal" relation. The projection norm of $M$ is defined as follows: for each $\phi \in \R^d$, 
\begin{equation}\label{eq:tangent-norm}
    \|\phi\|_{\calT M(\mathbf x)}=\sup_{\bfv\in \mathcal{T} M(\mathbf x)}\frac{\langle \mathbf v, \phi\rangle}{\|\mathbf v\|}.
\end{equation}
It is easy to see that algebraic form \cref{tmp:KKT} is equivalent to 
    \begin{equation}\label{tmp:KKT-2}
        \|\nabla f(\mathbf x)\|_{\calT M(\mathbf x)}=0.
    \end{equation}
We term this form \cref{tmp:KKT-2} the geometric form as it stems from the geometric perspective outlined above.


The optimality condition \cref{tmp:KKT-2} can be easily extended to Wasserstein space. To extend the geometric form \cref{tmp:KKT-2} of KKT condition to Wasserstein space, we need to extend the notion of ``tangent vector" first. As explained in \cref{ss:calculus-W} (particularly in \cref{lem:ac-curve}), for an absolutely continuous curve $\mu(\cdot):(0,1)\to(\calP_2(\R^d),W_2)$, there exists a time-dependent velocity field $v(\cdot,\cdot):(0,1)\times\R^d\to\R^d$ such that the continuity equation \cref{eq:continuity-equation} holds in weak sense. The ``tangent vector" to this curve $\mu(\cdot)$ at time $t$ is the instantaneous velocity field $v(t,\cdot):\R^d\to\R^d$ (velocity field at time $t$). Hence, it is natural to define ``tangent space" by the continuity equation \cref{eq:continuity-equation} as follows. Define the tangent space at $\bfmu=(\mu_1,\ldots,\mu_K)\in \calP_2(\R^d)^{\otimes K}$ as
    \begin{equation}\label{eq:tangent-space-W}
        \begin{split}
    \mathcal{T}\calP_{\pi,\bfp}(\bfmu_0)\doteq&\{(\phi_1,\ldots,\phi_K)\in\calL^2_{\bfmu}: 
    \exists \bfmu(t):[0,1]\rightarrow\calP_{\pi,\bfp},\text{ s.t. }\forall k\in[K],
    \\&\tfrac{d}{dt}\mu_k(t)|_{t=0}=-\nabla\cdot(\mu_k(0) \phi_k);\bfmu(0)=\bfmu_0 \}. 
        \end{split}
    \end{equation}
 The equation $\frac{d}{dt}\mu_k(t)\large|_{t=0}=-\nabla\cdot(\mu_k(0) \phi_k)$ holds in the weak sense; that is, for all $f\in C^\infty_c(\R^d)$, 
    \[
    \frac{d}{dt}\int_{\R^d} f(\mathbf x)d\mu_k(t,\mathbf x)|_{t=0}=\int_{\R^d}\langle\nabla f(\mathbf x),\phi_k(\mathbf x)\rangle d\mu_k(0,\mathbf x).
    \]
The \emph{projection norm} of tangent space $\calT\calP_{\pi,\bfp}(\bfmu)$ is 
\[
\|\bfphi\|_{\calT\calP_{\pi,\bfp}(\bfmu)}\doteq\sup_{u\in \mathcal{T} \calP_{\pi,\bfp}(\bfmu)}\frac{\langle u, \bfphi\rangle_{\bfmu}}{\|u\|_{\bfmu}},
\]
for all $\bfphi\in\calL^2_{\bfmu}$. We prove that the optimal solution to \cref{prob:decomposition-2} satisfies the following condition.

\begin{proposition}[Optimality condition for \cref{prob:decomposition-2}]\label{prop:opt-condition-1}
Suppose $\bfmu^*$ is the optimal solution to \cref{prob:decomposition-2}. Then, 
\[
\|\nabla F_\bfp(\bfmu^*)\|_{\mathcal{T}\calP_{\pi,\bfp}(\bfmu^*)}=0.
\]
\end{proposition}

Similarly, for \cref{prob:opt-decomposition}, we define the following tanget space, 
   \begin{equation*}
        \begin{split}
        \mathcal{T}\calP_{\pi}(\bfmu_0,\bfp_0)\doteq&\{(\bfphi,\bfv)\in\calL^2_{\bfmu}\times\calL^2:
\exists(\bfmu(t),\bfp(t))\in\calP_\pi,\text{ s.t. } \\
&\bfmu(t)\in \calP_{\pi,\bfp(t)}, \bfmu(0)=\bfmu_0, \text{ and } \bfp(0)=\bfp_0;\\
&\forall k\in[K], \frac{d}{dt}\mu_k(t)\large|_{t=0}=-\nabla\cdot(\mu_k(0) \phi_k) \text{ and } \tfrac{d}{dt}p_k(t)|_{t=0}=v_k\},
        \end{split}
    \end{equation*}
for all $(\bfmu_0,\bfp_0)\in\calP_2(\R^d)^{\otimes K}\times\R^{K}$. Note that this tangent space is not empty, since for any $(\bfmu_0,\bfp_0)$, it contains the zero element $(\bfp,\bfv)=(\mathbf{0},\mathbf{0})$.

The projection norm of $\calT\calP_\pi(\bfmu,\bfp)$ is, 
\[
\|(\bfphi,\bfv)\|_{\calT\calP_{\pi}(\bfmu,\bfp)}\doteq\sup_{(u_1,u_2)\in \mathcal{T} \calP_{\pi}(\bfmu,\bfp)}\frac{\langle u_1, \bfphi\rangle_{\bfmu}+\langle u_2, \bfv\rangle}{\|u_1\|_{\bfmu}+\|u_2\|}
\]
for all $(\bfphi,\bfv)\in\calL^2_{\bfmu}\times\calL^2_{\bfmu}$. We prove that the optimal solution to \cref{prob:opt-decomposition} satisfies the following condition.

\begin{proposition}\label{prop:opt-condition-2}
Suppose $(\bfmu^*,\bfp^*)$ is the optimal solution to \cref{prob:opt-decomposition}. Then, 
\[
\|\nabla F\|_{\mathcal{T}\calP_{\pi}(\bfmu^*,\bfp^*)}=0.
\]
\end{proposition}

As a remark, the main difficulty of proving \cref{prop:opt-condition-1,prop:opt-condition-2} by contradiction is constructing feasible solutions that are strictly better. This is nontrivial since the tangent space is not invariant, so moving along an element in the tangent space may leave the feasible set, and one needs to ``project" the infeasible solution back. All these procedures have to be done within the Wasserstein metric.

\section{Constraint controlled gradient flow (CCGF) for optimal probability measure decomposition}\label{s:algorithms}
We aim to develop a concept of gradient flow for \cref{prob:opt-decomposition,prob:decomposition-2} in  Wasserstein space. We first illustrate the notion of gradient flow in Euclidean space in \cref{ss:flow-Euclidean}. The purpose of elucidating this Euclidean gradient flow is to illustrate the main idea behind our flow design in Wasserstein space. Then, with the machinery introduced in \cref{s:preliminary}, we extend Euclidean gradient flow to Wasserstein space in \cref{ss:sub-problem,ss:main-problem}. We note that, in \cref{s:implementation}, we implement Wasserstein gradient flow via discretization. This yields an iterative algorithm.

\subsection{CCGF in Euclidean space}\label{ss:flow-Euclidean}

In this subsection, we design a constraint-controlled gradient flow (CCGF) in Euclidean space to solve the following constrained minimization problem. If we implement the CCGF, we can get an iterative algorithm. The concept of CCGF has also inspired the design of particle-based algorithms to solve constrained sampling problems in \cite{liu2021sampling} and \cite{zhang2022sampling}.
While the CCGF is an intuitive idea, we could not find any other references to it, even in standard Euclidean space settings. Therefore, we will first consider the following finite-dimensional optimization problem:

\begin{problem}[Finite-dimensional optimization problem]\label{opt:finite}
Let $f:\R^n\rightarrow\R$ and $g:\R^n\rightarrow[0,\infty]$ be two smooth functions and solve
     \begin{equation*}   
        \begin{split}
            \min_{\mathbf x\in\mathbb{R}^n}&\quad f(\mathbf x)\\
            \text{s.t.}&\quad g(\mathbf x)= 0,
        \end{split}
    \end{equation*}
    where $C_0\doteq\{\mathbf x\in\mathbb{R}^n: g(\mathbf x)= 0\}$ be the feasible region.
\end{problem}
     We will impose necessary assumptions when they are needed in later development. We do not assume that $f,g$ are convex or concave.
    
 To solve Problem \ref{opt:finite}, a first thought is to use the projected gradient descent algorithm. Let $d:\R^n\times\R^n\rightarrow [0,\infty)$ be the Euclidean metric. The project gradient descent algorithm is given by 
    \begin{equation*}
        \begin{split}
            \mathbf x'_{k+1} &= \argmin_{y\in\R^n} f(\mathbf x_k)+\langle\nabla f(\mathbf x_k),y-\mathbf x_k\rangle + \frac{1}{2}d(y,\mathbf x_k)^2\\
            x_{k+1} &= \Pi_{C_0}(\mathbf x'_{k+1}),
        \end{split}
    \end{equation*}
where $\Pi_{C_0}(\mathbf x'_{k+1})$ is the projection of $\mathbf x'_{k+1}$ to the feasible set $C_0$. However, such an algorithm can be hard to implement for two reasons: (1) projections can be difficult to find, and (2) the projected point $\mathbf x_{k+1}$ can be far away from $\mathbf x'_{k+1}$, which makes it difficult to derive the gradient flow in the limit. 

We propose the following ``variational interpolation" approach to design our gradient flow to avoid these two challenges. The idea of the variational interpolation approach is as follows. We aim to design a vector field $\phi:\R^d\rightarrow\R^d$ to specify the ``velocity" of movement at each position in $\R^d$. Accordingly, we call it a ``velocity field." The purpose of our design is to guarantee that if a particle is moving in $\R^d$ according to this velocity field $\phi$, it will eventurally converge to a feasible solution $\mathbf x^*$ (i.e., $g(\mathbf x^*)=0$) which satisfies optimality condition \cref{tmp:KKT-2}. Recall that the trajectory of the particle is given by the solution to the ODE: $\dot{\mathbf x}(t)=\phi(\mathbf x(t))$. Hence, mathematically speaking, our challenge is to design $\phi$ such that the solution to the ODE $\dot{\mathbf x}(t)=\phi(\mathbf x(t))$ converges to a feasible solution to \cref{opt:finite} (i.e., $\lim_{t\to\infty}g(\mathbf x(t))=0$) satisfying optimality condition \cref{tmp:KKT-2} in the end. To design such a velocity field $\phi$, we start with a sequence of discrete points generated by an iterative scheme explained as follows.

Fixing any time step parameter $\tau>0$, for some constraint-control parameter $\alpha\in(0,\sfrac{1}{\tau})$, we look for a sequence of points $(\mathbf x^\tau_k)_{k\in\mathbb{N}}$ defined through the following iterated scheme,
\begin{equation}\label{eq:temp-iterative}
     \mathbf x^\tau_{k+1} \doteq \argmin_{\mathbf x\in C_{(1-\alpha \tau)g(\mathbf x^\tau_k)}} f(\mathbf x^\tau_k)+\langle\nabla f(\mathbf x^\tau_k),\mathbf x-\mathbf x^\tau_k\rangle + \frac{1}{2\tau}d(\mathbf x,\mathbf x^\tau_k)^2,
\end{equation}
where $C_{(1-\alpha \tau)g(x^\tau_k)}=\{\mathbf x\in\mathbb{R}^n:g(\mathbf x)= (1-\alpha \tau)g(\mathbf x^\tau_k)\}$.  We can interpret this sequence of points as the positions of the trajectory $\mathbf x(t)$ at time points $t=0,\tau,2\tau,\ldots,k\tau,\ldots$; that is, $\mathbf x(k\tau) = \mathbf x_k^\tau$. Hence, the parameter $\tau$ represents the time length between two consecutive points $\mathbf x^\tau_k$ and $\mathbf x^\tau_{k+1}$.

If $g(x^\tau_k)=0$, i.e., $\mathbf x^\tau_k\in C_0$, the $(k+1)$-st step produced by equation \cref{eq:temp-iterative} is the same as the original projected gradient descent algorithm. If $g(\mathbf x^\tau_k)>0$, i.e., $\mathbf x^\tau_k\notin C_0$, we try to find a solution $\mathbf x^\tau_{k+1}$ which is closer to the actual feasible region $C_0$ in the sense that $g(\mathbf x^\tau_{k+1})= (1-\alpha\tau)g(\mathbf x^\tau_k)$, where $ (1-\alpha\tau)\in(0,1)$.

Moreover, we can further simplify the equation \cref{eq:temp-iterative} by replacing the nonlinear constraint $g(\mathbf x)= (1-\alpha\tau)g(\mathbf x^\tau_k)$ by its linear approximation. This linearization simplifies the problem and makes it easily solvable using Lagrangian-type arguments as shown in \cref{lem:discrete-iteration}. It is clear that, as $\tau\rightarrow 0^+$, we have $\mathbf x^\tau_{k+1}\rightarrow \mathbf x^\tau_k$. Hence, if time step $\tau$ is small enough, the $(k+1)$-st step $\mathbf x^\tau_{k+1}$ should be close to the $k$-th step $\mathbf x^\tau_k$. Therefore, we can replace the function $g$ in the left-hand side of the constraint defining $C_{(1-\alpha \tau) g(\mathbf x_k^\tau)}$ with its linear approximation $\Tilde{g}(\mathbf x)=g(\mathbf x^\tau_k)+\langle \mathbf x-\mathbf x^\tau_k,\nabla g(\mathbf x^\tau_k)\rangle$. Thus, now we consider the sequence $(\mathbf x^\tau_k)_{k\in\mathbb{N}}$ generated by the following iterative scheme
    \begin{equation}\label{eq:discrete-scheme}
        \begin{split}
          \mathbf x^\tau_k\doteq\argmin_{\mathbf x\in\R^n}  &\quad f(\mathbf x^\tau_k)+\langle \nabla f(\mathbf x^\tau_k),\mathbf x-\mathbf x^\tau_k\rangle+\frac{1}{2\tau}d(\mathbf x,\mathbf x^\tau_k)^2\\
            \text{s.t.}&\quad \Tilde{g}(\mathbf x)= (1-\alpha \tau)g(\mathbf x^\tau_k),
        \end{split}
    \end{equation}
where the constraint is equivalent to $\alpha \tau g(\mathbf x^\tau_k)+\langle\nabla g(\mathbf x^\tau_k),\mathbf x-\mathbf x^\tau_k\rangle= 0$. We solve \cref{eq:discrete-scheme} with the Lagrangian method as follows. 

\begin{lemma}\label{lem:discrete-iteration}
    For any $\tau>0$, at any iteration $t$, the optimal solution $\mathbf x^\tau_{k+1}$ to equation \cref{eq:discrete-scheme} is given by 
    \begin{equation*}
         \mathbf x^\tau_{k+1}=\mathbf x^\tau_k-\tau (\nabla f(\mathbf x^\tau_k)+\lambda^* \nabla g(\mathbf x^\tau_k)),\quad 
\lambda^*=\frac{-\langle \nabla g(\mathbf x^\tau_k), \nabla f(\mathbf x^\tau_k)\rangle +\alpha g(\mathbf x^\tau_k)}{\|\nabla g(\mathbf x^\tau_k)\|^2}.
    \end{equation*}
\end{lemma}
Suppose a particle travels at a constant velocity from point $\mathbf x^\tau_k$ to point $\mathbf x^\tau_{k+1}$. Since the duration between two consecutive points is $\tau$, the velocity $\phi^\tau_k$ of moving from position $\mathbf x^\tau_{k}$ to position $\mathbf x^\tau_{k+1}$ is given by
\[
 \phi^\tau_k\doteq\frac{\mathbf x^\tau_{k+1}-\mathbf x^\tau_k}{\tau}=-(\nabla f(\mathbf x^\tau_k)+\lambda^* \nabla g(\mathbf x^\tau_k)).
\]
This motivates us to consider the following design of the velocity field $\phi$. 

\begin{definition}[Euclidean CCGF]\label{def:Euclidean-CCGF}
    The \emph{Euclidean constraint controlled gradient flow (CCGF)} is defined to be the solution to the ODE: \[\dot{\mathbf x}(t)=\phi(\mathbf x(t)),\] where the velocity field $\phi:\R^n\rightarrow\R^n$ is
        \begin{equation}\label{eqn:CCGFEuclidean}  
          \phi(\mathbf x)\doteq-(\nabla f(\mathbf x)+\lambda(x)\nabla g(\mathbf x))\text{ and }\lambda(\mathbf x)=\frac{-\langle\nabla g(\mathbf x),\nabla f(\mathbf x)\rangle+\alpha g(\mathbf x)}{\|\nabla g(\mathbf x)\|^2}.
    \end{equation}
\end{definition}



We show the convergence of the Euclidean CCGF in \cref{def:Euclidean-CCGF} under the following assumptions.

\begin{assumption}\label{as:CCGFEuclidean-g}
$\|\nabla g(\mathbf x)\|>0$ if $g(\mathbf x)\neq 0$.
\end{assumption}

\cref{as:CCGFEuclidean-g} guarantees that $\lambda(\mathbf x)$ in \cref{eqn:CCGFEuclidean} is well-defined. We also need the following assumptions to prove convergence.

\begin{assumption}\label{as:CCGFEuclidean}
    \begin{enumerate}[label=(\roman*)]
        \item $f$ is bounded from below by a real number $f_{\text{min}}$;
        \item there exists an $L>0$ such that $\|\nabla f(\mathbf x)\|\leq L$ for all $\mathbf x\in\R^n$;
        \item (Polyak-\L ojasiewicz condition) there exists a $\kappa>0$ such that $\|\nabla g(\mathbf x)\|^2\geq\kappa g(\mathbf x)$ for all $\mathbf x\in\R^n$.
    \end{enumerate}
\end{assumption}

As we discussed in \cref{s:optimality-condition}, the KKT condition has a geometric form described in equation \cref{tmp:KKT-2}, using the tangent norm defined in equation \cref{eq:tangent-norm}. We show that, in the limit (i.e., as $t\to\infty$), the KKT condition \cref{tmp:KKT-2} is satisfied by the Euclidean CCGF in \cref{def:Euclidean-CCGF}.

\begin{theorem}\label{thm:convergence-Euclidean}
    Suppose \cref{as:CCGFEuclidean-g} holds. Let $\mathbf x(t)$ be the Euclidean CCGF defined in \cref{def:Euclidean-CCGF}. Then, 
    \begin{enumerate}[label=(\roman*)]
        \item $g(\mathbf x(t))= e^{-\alpha t}g(\mathbf x(0))$ for all $t>0$;
        \item If \cref{as:CCGFEuclidean} also holds, then, for any $T > 0$, 
        \[
        \min_{t\leq T}\|\nabla f(\mathbf x(t))\|_{\cT M(\mathbf x)}\leq\frac{C}{\sqrt{T}},
        \]
    where $C\doteq f(\mathbf x(0))-f_{\text{min}}+\frac{2g(\mathbf x(0))}{\kappa}+\frac{L}{\alpha\sqrt{\kappa}}\sqrt{g(\mathbf x(0))}$.
    \end{enumerate}
\end{theorem}

\cref{thm:convergence-Euclidean} shows that the CCGF $\mathbf x(t)$ in \cref{def:Euclidean-CCGF} converges in the limit to a feasible solution (i.e., $\lim_{t \to \infty}g(\mathbf x(t))=0$) exponentially fast with speed controlled by the parameter $\alpha$. Furthermore, the tangent norm of the gradient decreases to $0$ at the speed of $O(\frac{1}{\sqrt{T}})$, i.e., the KKT condition \cref{tmp:KKT-2} is satisfied in the limit as $t\to\infty$.

To implement the Euclidean CCGF approach in practice, we can select a short time length $\tau$ and implement the iterative scheme in \cref{lem:discrete-iteration}. After enough iterations (i.e., $k$ is large enough, approximating $t\to\infty$), the solution $\mathbf x^\tau_k$ will be close to feasible and optimal.  

\begin{remark}
\label{rem:CCGFequality}
Another version of the CCGF can be derived with the constraint in \cref{eq:discrete-scheme} is replaced with $\tilde{g}(x)\leq (1-\alpha \tau)g(\bfx_k^\tau)$. The corresponding CCGF remains largely the same, except $\lambda(\bfx)$ is replaced by its positive part, i.e $\max\{\lambda(\bfx),0\}$. Gradient flow approaches that use a structure analogous to $\max\{\lambda(\bfx),0\}$ can be found in \cite{liu2021sampling,zhang2022sampling}, but in a Wasserstein context, not in Euclidean space. Similar bounds as in \cref{thm:convergence-Euclidean} can be obtained, except one can only obtain $g(\mathbf x(t))\leq e^{-\alpha t}g(\mathbf x(0))$.
While this provides faster convergence to the feasible set, it can also leads to instability in the computation of $\lambda(x)$, since $\nabla g(\bfx)$ appears in the denominator. Using the current version of CCGF resolves this issue.
\end{remark}

\subsection{Wasserstein CCGF for \cref{prob:decomposition-2}}\label{ss:sub-problem}
In this subsection, we extend the notion of CCGF in \cref{def:Euclidean-CCGF} to Wasserstein space to solve \cref{prob:decomposition-2}.

For ease of development, from now on, we assume $\pi$ is absolutely continuous with respect to Lebesgue measure; that is, we think of $\pi$ as a probability density function. 

Now, $F_{\bfp}$ in \cref{prob:decomposition-2} plays the role of the objective function $f$ in \cref{opt:finite}. We use Kullback-Leibler (KL) divergence to measure deviations from the feasibility constraint $\sum_{k\in[K]}p_k\mu_k=\pi$. Recall that the KL divergence between two densities is $0$ if and only if the two densities are equal almost everywhere. Define $\mubar\doteq\sum_{k\in[K]} p_k\mu_k$. The constraint $g(x)=0$ in \cref{opt:finite} becomes $\text{KL}(\mubar\|\pi)=0$ in \cref{prob:decomposition-2}.

Since \cref{prob:decomposition-2} is defined on the product Wasserstein space $(\calP_{2,ac}(\R^d),W_2)^{\otimes K}$, our approach is to design a velocity field $\bfphi=(\phi_1,\ldots,\phi_K)$ on $(\calP_{2,ac}(\R^d),W_2)^{\otimes K}$ to solve \cref{prob:decomposition-2}. Given any point $\bfmu=(\mu_1,\ldots,\mu_K)\in(\calP_{2,ac}(\R^d),W_2)^{\otimes K}$, the corresponding velocity field $\bfphi(\bfmu)$ should specify the ``velocity of movement" at position $\bfmu$. Recall that, as we discussed in \cref{ss:calculus-W}, in a single Wasserstein space $(\calP_{2,ac}(\R^d),W_2)$, the ``velocity" of any trajectory at any position $\mu\in (\calP_{2,ac}(\R^d),W_2)$ is a vector field $\phi(\mu)(\cdot):\R^d\rightarrow\R^d$ on Euclidean space. Given such a velocity field, the trajectory of movement in Wasserstein space is the solution to the continuity equation \cref{eq:continuity-equation}. Hence, the velocity field $\bfphi$ in the product Wasserstein space $(\calP_{2,ac}(\R^d),W_2)^{\otimes K}$ should map any $\bfmu\in(\calP_{2,ac}(\R^d),W_2)^{\otimes K}$ to $K$ vector fields $\bfphi(\bfmu)=(\phi_1(\bfmu),\ldots,\phi_K(\bfmu))$, where each $\phi_k(\bfmu)(\cdot):\R^d\to\R^d$ in $\calL^2_{\mu_k}$ describes the velocity of $\mu_k\in (\calP_{2,ac}(\R^d),W_2)$.

\cref{prob:decomposition-2} can be solved using our velocity field $\bfphi$ by the following procedure. Starting from arbitrary initial densities $\bfmu(0)=(\mu_1(0),\ldots,\mu_K(0))$, the solution $\bfmu(t):\R\to(\calP_{2,ac}(\R^d),W_2)^{\otimes K}$ of the following $K$ continuity equations
    \begin{equation}\label{eq:continuity}
        \frac{d}{dt}\mu_k(t) = -\nabla\cdot(\phi_k(\bfmu(t))\mu_k(t)),\quad\forall k\in[K]
    \end{equation}
    converges to a ``stationary point" of \cref{prob:decomposition-2} in the sense that 
    \begin{itemize}
        \item $\lim_{t\rightarrow\infty}\text{KL}(\mubar(t)\|\pi)=0$ with $\mubar(t)\doteq\sum_{k\in[K]} p_k\mu_k(t)$;
        \item the optimality condition in \cref{prop:opt-condition-1} is approximately satisfied.
    \end{itemize}

Following the Euclidean CCGF in \cref{def:Euclidean-CCGF}, we design the Wasserstein CCGF as follows.

\begin{definition}[Wasserstein CCGF]\label{def:Wassertein-CCGF}
 The \emph{Wasserstein constraint controlled gradient flow (CCGF)} is the solution $\bfmu(t):\R\rightarrow(\calP_{2,ac}(\R^d),W_2)^{\otimes K}$ to the $K$ continuity equations \cref{eq:continuity}, where the velocity field $\phi_k(\bfmu)(\cdot):\R^d\rightarrow\R^d$ is
    \begin{equation*}
             \phi_k(\bfmu)=-(\nabla_{\mu_k}F_\bfp(\bfmu)+\lambda(\bfmu)\nabla_{\mu_k}\text{KL}(\mubar\|\pi)),\quad\forall k\in[K]
    \end{equation*}
    with
    \begin{equation*}
             \lambda(\bfmu)=\frac{-\langle\nabla_{\bfmu} F_\bfp(\bfmu),\nabla_{\bfmu}\text{KL}(\mubar\|\pi)\rangle_{\bfmu}+\alpha\text{KL}(\mubar\|\pi)}{\|\nabla_{\bfmu}\text{KL}(\mubar\|\pi)\|^2_{\bfmu}}.
    \end{equation*}
\end{definition}

By \cref{lem:U,lem:KLder},  we have 
\begin{equation}\label{eq:WCCGF}
        \begin{split}
             \phi_k(\bfmu)&=-(p_k\nabla L(\mu_k)+\lambda(\bfmu)p_k(s_{\mubar}-s_\pi)))\in\calL^2_{\mu_k},\quad\forall k\in[K]\text{ and}\\
             \lambda(\bfmu)&=\frac{-\sum_{k\in[K]}p_k\langle\nabla L(\mu_k),p_k(s_{\mubar}-s_\pi)\rangle_{\mu_k}+\alpha\text{KL}(\mubar\|\pi)}{\sum_{k\in[K]}\|p_k(s_{\mubar}-s_\pi)\|^2_{\mu_k}}\in\R,
        \end{split}
    \end{equation}
where for any $\mathbf x\in\R^d$, $\nabla L(\mu_k)(\mathbf x),s_{\mubar}(\mathbf x),s_\pi(\mathbf x)$ are vectors in $\R^d$ given by 
    \begin{equation*}
        \begin{split}
            \nabla L(\mu_k)(\mathbf x)&=\int_{\R^d} (\nabla_1 \ell(\mathbf x,\mathbf z)+\nabla_2 \ell(\mathbf z,\mathbf x))d\mu_k(\mathbf z),\\
            s_{\mubar}(\mathbf x)&=\nabla\log\mubar(\mathbf x)\text{ and }
            s_{\pi}(\mathbf x)=\nabla\log\pi(\mathbf x).
        \end{split}
    \end{equation*}

We make the following assumptions similar to \cref{as:CCGFEuclidean,as:CCGFEuclidean-g} in Euclidean space. 

\begin{assumption}\label{as:WCCGF}
    \begin{enumerate}[label=(\roman*)]
        \item the kernel $\ell$ is bounded, i.e., for some $\ell_{\text{max}}>0$, $|\ell(\mathbf x,\mathbf y)|<\ell_{\text{max}}$ for all $\mathbf x,\mathbf y\in\R^d$;
        \item the gradient of kernel $\ell$ is bounded, i.e., for some $L_{\text{max}}>0$, $\|\nabla\ell(\mathbf x,\mathbf y)\|<L_{\text{max}}$ for all $\mathbf x,\mathbf y\in\R^d$.
        \item $\pi$ follows $\kappa$-log Soblev inequality, i.e., for some $\kappa>0$, $\|s_\nu-s_\pi\|_\nu^2\geq\kappa\text{KL}(\nu\|\pi)$ for all $\nu\in\calP_{2,ac}(\R^d)$.
    \end{enumerate}
\end{assumption}
It is known that $\pi$ with bounded support or is strongly-log concave satisfies the $\kappa$-log Sobolev inequality. We prove the following convergence result for the Wasserstein CCGF defined in \cref{def:Wassertein-CCGF}.


\begin{theorem}\label{thm:convergence-1}
  Suppose \cref{as:WCCGF} holds. Let $\bfmu(t)$ be the Wasserstein CCGF in \cref{def:Wassertein-CCGF}. Then,
        \begin{enumerate}[label=(\roman*)]
            \item KL-divergence between $\mubar(t)$ and $\pi$ decreases exponentially, i.e., $\text{KL}(\mubar(t)\|\pi)= e^{-\alpha t}\text{KL}(\mubar(0)\|\pi)$;
            \item optimality condition in \cref{prop:opt-condition-1} is approximately satisfied in the sense that for $T>0$,
                \[
                \min_{t\leq T}\|\nabla_{\bfmu} F_\bfp(\bfmu(t))\|_{\mathcal{T}\calP_{\pi,\bfp}(\bfmu(t))}\leq \frac{C}{\sqrt{T}},
                \]
                where 
                    \[
                    C\doteq \tfrac{1}{\sqrt{K}}(F_\bfp(\bfmu(0))+\ell_{\text{max}}+\frac{4\alpha L_{\text{max}}\sqrt{K}}{p_{\text{min}}\sqrt{\kappa}}\sqrt{\text{KL}(\mubar(0)\|\pi)} + \alpha\frac{1}{p_{\text{min}}\kappa}\text{KL}(\mubar(0)\|\pi)),\]
                    with $p_{\text{min}}=\min_{k\in[K]}p_k$.
        \end{enumerate}
\end{theorem}

\cref{thm:convergence-1} shows that the Wasserstein CCGF $\bfmu(t)$ in \cref{def:Wassertein-CCGF} converges in the limit to a feasible solution (i.e., $\lim_{t \to \infty}\text{KL}(\mubar(t)\|\pi)=0$) exponentially fast with speed controlled by the parameter $\alpha$. Furthermore, the tangent norm of the Wasserstein gradient decreases to $0$ at the speed of $O(\frac{1}{\sqrt{T}})$, i.e., the KKT condition in \cref{prop:opt-condition-1} is satisfied in the limit as $t\to\infty$. We discuss an algorithmic implementation of the Wasserstein CCGF in \cref{s:implementation}.  

\subsection{Wasserstein CCGF for \cref{prob:opt-decomposition}}\label{ss:main-problem}
In this subsection, we extend the Wasserstein CCGF algorithm to solve \cref{prob:opt-decomposition}. We must also design a velocity field for the extra decision variable $p_k$. 

The decision variable $(\bfmu,\bfp)=((\mu_1,\ldots,\mu_K),(p_1,\ldots,p_K))$ of \cref{prob:opt-decomposition} is defined in the space $(\calP_{2,ac}(\R^d),W_2)^{\otimes K}\times\R_+^K$. Similar to \cref{ss:sub-problem}, we design a velocity field $\bfphi$ in the space $ (\calP_{2,ac}(\R^d),W_2)^{\otimes K} $ for $\bfmu$ and a velocity field $\bfv$ in the space $\R_+^K$ for $\bfp$. Specifically, given any $(\bfmu,\bfp)\in(\calP_{2,ac}(\R^d),W_2)^{\otimes K}\times\R_+^K$, $\bfphi(\bfmu,\bfp)=(\phi_1(\bfmu,\bfp),\ldots,\phi_K(\bfmu,\bfp))$ is a vector of velocity fields, where each $\phi_k(\bfmu,\bfp)(\cdot):\R^d\to\R^d\in\calL^2_{\mu_k}$ describes the velocity of decision variable $\mu_k$ in Wasserstein space; and $\bfv(\bfmu,\bfp)(\cdot):\R_+^K\rightarrow\R_+^K\in\calL^2$ is a velocity field of the decision variable $\bfp\in\R_+^K$.

Our chosen velocity fields $(\bfphi,\bfv)$ should solve \cref{prob:opt-decomposition} in the following sense. Starting with initial densities $\bfmu(0)=(\mu_1(0),\ldots,\mu_K(0))$ and initial weights $\bfp(0)=(p_1(0),\ldots,p_K(0))$, the solution of the following system of partial differential equations, 
    \begin{equation}\label{eq:pde}
        \begin{split}
               \frac{d}{dt}\mu_k(t) &= -\nabla\cdot(\phi_k(\bfmu(t),\bfp(t))\mu_k(t))\\
               \frac{d}{dt}p_k(t) &= v_k(\bfmu(t),\bfp(t))
        \end{split}
    \end{equation}
    for $k \in [K]$ converges to a ``stationary point" of \cref{prob:opt-decomposition} in the sense that 
    \begin{itemize}
        \item $\sum_{k\in[K]}p_k(t)=1$, and $p_k(t)\geq 0$, for all $k\in[K]$ and all $t\geq 0$;
        \item $\lim_{t\rightarrow\infty}\text{KL}(\mubar(t)\|\pi)=0$ with $\mubar(t)=\sum_{k\in[K]} p_k(t)\mu_k(t)$;
        \item the optimality condition in \cref{prop:opt-condition-2} is approximately satisfied.
    \end{itemize}

To see this, we first need gradients of the objective function $F(\bfmu,\bfp)=\sum_{k\in[K]}p_kL(\mu_k)+\sfrac{\theta}{p_k^\beta}$ and the constraint $\text{KL}(\mubar\|\pi)$. The Wasserstein gradient of both $F(\bfmu,\bfp)$ and $\text{KL}(\mubar\|\pi)$ with respect to $\bfmu$ remains the same as given by \cref{lem:U,lem:KLder}. The following lemma gives the gradients with respect to $\bfp$. 

\begin{lemma}
    The gradients $\nabla_{\bfp}F(\bfmu,\bfp)\in\R^d$ and $\nabla_{\bfp}\text{KL}(\mubar\|\pi)\in\R^d$ are 
    \begin{equation*}
        \begin{split}
            \nabla_{\bfp}F(\bfmu,\bfp)_k&=L(\mu_k)-\frac{\theta\beta}{p_k^{\beta+1}},\\
            \nabla_{\bfp}\text{KL}(\bfmu\|\pi)_k&=\int_{\R^d}\mu_k(x)\log\frac{\mubar(x)}{\pi(x)}dx+1.
        \end{split}
    \end{equation*}
\end{lemma}

We define the following objects to extend the Wasserstein CCGF for \cref{prob:decomposition-2}. Recall in \cref{s:preliminary}, given any $\bfmu\in\calP_{2,ac}(\R^d)^{\otimes K}$, we define inner product $\langle\cdot,\cdot,\rangle_{\bfmu}$ in the space $\calL^2_{\bfmu}(\R^d;\R^d)$ by
$\langle f,g\rangle_{\bfmu}=\sum_{k\in[K]}\int_{\R^d}\langle f_k(x),g_k(x)\rangle d\mu_k(x)$
for any $f,g\in\calL^2_{\bfmu}(\R^d;\R^d)$. Particularly, since $(\nabla_{\bfmu}F(\bfmu,\bfp)$, $\nabla_{\bfmu}\text{KL}(\mubar\|\pi))$ are both in $\calL^2_{\bfmu}(\R^d;\R^d)$, we have
 \[
    \langle\nabla_{\bfmu}F(\bfmu,\bfp),\nabla_{\bfmu}\text{KL}(\mubar\|\pi)\rangle_{\bfmu} = \sum_{k\in[K]}\int_{\R^d}\langle\nabla_{\mu_k} F(\bfmu,\bfp),\nabla_{\mu_k}\text{KL}(\mubar\|\pi)\rangle d\mu_k.
    \]

To ensure that $\sum_{k\in[K]}p_k(t)=1$ for all $t\geq 0$, we project the velocity field $\bfv(t)$ for $\bfp(t)$ so that $\sum_{k\in[K]}v_k=0$. Let $\mathbf{1}$ denote the column vector $(1,\ldots,1)^\top\in\R^K$. We define the projection operator $P:\R^K\rightarrow\R^K$ by
    \[
    Pv = (\text{id} - \frac{1}{K}\mathbf{1}\mathbf{1}^\top)v = v - \frac{\sum_{k\in[K]}v_k}{K}\mathbf{1}.
    \]
for $v\in\R^K$. This projection ensures that $\sum_{k\in[K]}Pv_k=0$. With this projection operator, we can define an inner product and norm: given any $\bfmu\in\calP_2(\R^d)^{\otimes K}$ and $\mathbf u=(u_1,u_2),\mathbf w=(w_1,w_2)$ with $u_1,w_1\in\calL^2_{\bfmu}(\R^d;\R^d)$ and $u_2,w_2\in\R^K$,
    \begin{equation*}
        \begin{split}
                   \langle \mathbf u,\mathbf w\rangle_{\bfmu,P} &= \langle u_1,w_1\rangle_{\bfmu}+\langle Pu_2,Pw_2\rangle,\\
                   \|\mathbf u\|^2_{\bfmu,P} &= \|u_1\|^2_{\bfmu} + \|Pu_2\|^2,
        \end{split}
    \end{equation*}
where $\langle Pu_2,Pw_2\rangle$ and $\|Pu_2\|$ are the Euclidean inner product and norm, respectively. 

In particular, since $ \nabla F(\bfmu,\bfp)=(\nabla_{\bfmu}F(\bfmu,\bfp), \nabla_{\bfp}F(\bfmu,\bfp))$ and $\nabla\text{KL}(\mubar\|\pi)=(\nabla_{\bfmu}\text{KL}(\mubar\|\pi),\nabla_{\bfp}\text{KL}(\mubar\|\pi))$ are both in $\calL^2_{\bfmu}\times\R^K$, we have 
\begin{equation*}
    \begin{split}
    \langle \nabla F(\bfmu,\bfp),\nabla\text{KL}(\mubar\|\pi)\rangle_{\bfmu,P} &=    \langle\nabla_{\bfmu}F(\bfmu,\bfp),\nabla_{\bfmu}\text{KL}(\mubar\|\pi))\rangle_{\bfmu} + \langle P\nabla_{\bfp}F(\bfmu,\bfp),P\nabla_{\bfp}\text{KL}(\mubar\|\pi) \rangle,\\
    \|\nabla\text{KL}(\mubar\|\pi)\|^2_{\bfmu,P} &=  \|\nabla_{\bfmu}\text{KL}(\mubar\|\pi)\|^2_{\bfmu} + \|P\nabla_{\bfp}\text{KL}(\mubar\|\pi)\|^2.
    \end{split}
\end{equation*}

\begin{definition}[Wasserstein CCGF with dynamic weights]\label{def:Wasserstein-CCGF2}
    The \emph{Wasserstein constraint controlled gradient flow with dynamic weights} is defined to be the solution $\bfmu(t):\R\rightarrow(\calP_{2,ac}(\R^d),W_2)^{\otimes K}$ and $\bfp(t):\R\to\R^K$ to the system \cref{eq:pde} of PDEs, where the velocity fields $\phi_k(\bfmu,\bfp)(\cdot):\R^d\rightarrow\R^d$ and $\bfv(\bfmu,\bfp)(\cdot):\R^K\to\R^K$ are, for all $k\in[K]$,
     \begin{equation}\label{eq:WCCGF2}
        \begin{split}
             \phi_k(\bfmu,\bfp)&=-(\nabla_{\mu_k}F(\bfmu,\bfp)+\lambda(\bfmu,\bfp)\nabla_{\mu_k}\text{KL}(\mubar\|\pi)),\\
            v_k(\bfmu,\bfp) &= -P(\nabla_{\bfp}F(\bfmu,\bfp)_k+\lambda(\bfmu,\bfp)\nabla_{\bfp}\text{KL}(\mubar\|\pi)_k), \text{ and }\\
             \lambda(\bfmu,\bfp)&=\frac{-\langle \nabla F(\bfmu,\bfp),\nabla\text{KL}(\mubar\|\pi)\rangle_{\bfmu,P} +\alpha\text{KL}(\mubar\|\pi)}{\|\nabla\text{KL}(\mubar\|\pi)\|^2_{\bfmu,P}}.
        \end{split}
    \end{equation}
\end{definition}

We prove the following convergence result for Wasserstein CCGF with dynamic weights defined in \cref{def:Wasserstein-CCGF2}.

\begin{theorem}\label{thm:convergence-2}
    Suppose \cref{as:WCCGF} holds. Let $(\bfmu(t),\bfp(t))$ be the Wasserstein CCGF with dynamic weights defined in \cref{def:Wasserstein-CCGF2}. Then,
        \begin{enumerate}[label=(\roman*)]
            \item $\sum_{k\in[K]}p_k(t)=1$ and $\bfp(t)\geq p_{\min}$ $m$-almost surely, for some $p_{\min}>0$;
            \item the KL-divergence between $\mubar(t)$ and $\pi$ decreases exponentially, i.e., $\text{KL}(\mubar(t)\|\pi)= e^{-\alpha t}\text{KL}(\mubar(0)\|\pi)$;
            \item optimality condition in \cref{prop:opt-condition-2} is approximately satisfied in the sense that for all $T > 0$,
                \[
                \min_{t\leq T}\|\nabla F(\bfmu(t),\bfp(t))\|_{\mathcal{T}\calP_{\pi}(\bfmu(t),\bfp(t))}\leq \frac{C}{\sqrt{T}},
                \]
                where 
                    \[
                    C\doteq F(0)+\ell_{\text{max}}+\frac{4}{\alpha p_{\min}\sqrt{\kappa}}\sqrt{\text{KL}(0)}\sqrt{L_{\text{max}}^2K^2+K(\ell_{\text{max}}+\tfrac{\theta\beta}{p_{\min}^{\beta+1}})^2} + \frac{1}{p_{\min}\kappa}\text{KL}(0).\]
                   
        \end{enumerate}
\end{theorem}

\cref{thm:convergence-2} shows that the Wasserstein CCGF with dynamic weights $(\bfmu(t),\bfp(t))$ in \cref{def:Wasserstein-CCGF2} (i) guarantees that $\bfp(t)\geq 0$ and $\sum_{k\in[K]}p_k(t)=1$ almost surely for all $t>0$; (ii) converges in the limit to a feasible solution (i.e., $\lim_{t \to \infty}\text{KL}(\mubar(t)\|\pi)=0$) exponentially fast with speed controlled by the parameter $\alpha$. Furthermore, the tangent norm of the Wasserstein gradient decreases to $0$ at the speed of $O(\frac{1}{\sqrt{T}})$, i.e., the KKT condition in \cref{prop:opt-condition-2} is satisfied in the limit as $t\to\infty$. An algorithmic implementation of this flow is discussed in \cref{s:implementation}.  

\section{Implementation of Wasserstein CCGF}\label{s:implementation}

In this section, we implement \cref{def:Wassertein-CCGF,def:Wasserstein-CCGF2} in the settings described by \cref{exp:clustering,exp:elo}. In particular, we conduct a detailed analysis of the league design problem described in \cref{exp:elo}.

\subsection{Implementation of Wasserstein flow in \cref{def:Wassertein-CCGF}}

We begin with \cref{def:Wassertein-CCGF}. Recall that we aim to generate a curve of probability densities $\bfmu(\cdot):[0,T]\rightarrow\calP_{2,ac}(\R^d)^{\otimes K}$ (i.e., $\bfmu(t)=(\mu_1(t),\ldots,\mu_K(t))$) to solve \cref{prob:decomposition-2}. By \cref{thm:convergence-1}, this is achieved by solving the system of $K$ continuity equations \cref{eq:continuity} given by our velocity field design $\bfphi(\bfmu)=(\phi_1(\bfmu),\ldots,\phi_K(\bfmu))$ in \cref{def:Wassertein-CCGF}. Implementing this Wasserstein CCGF involves two challenges. First, each $\mu_k(t)$ is infinite-dimensional and non-parametric, thus hard to represent. Second, we need to decide $\mu_k(t)$ for each $t$ in a continuous interval $[0,T]$. We address these challenges via the ``particles method". Specifically, we use a large population of particles sampled from $\mu_k(t)$ to represent the distribution with density $\mu_k(t)$ and discretize the time horizon $[0,T]$ by considering time steps $t=0,1,2,\ldots$ to get an iterative algorithm.

For any curve $\bfmu(\cdot):[0,T]\rightarrow\calP_{2,ac}(\R^d)^{\otimes K}$, at each $t\in[0,T]$, let $\phi_k(t,\cdot):\R^d\to\R^d$ denote the velocity field $\phi_k(\bfmu(t))$ given in \cref{def:Wassertein-CCGF}. For each $k\in[K]$, we sample $N$ particles from an initial probability density $\mu_k(0)$. If all particles move together according to velocity field $\phi_k$ (i.e., the trajectory of each particle initially located at $\mathbf x_0$ is given by the solution of $\dot{\mathbf x}(t)=\phi_k(t,\mathbf x(t))$ with $\mathbf x(0)=\mathbf x_0$), then by \cref{lem:implementation}, the curve $\mu_k(t)$ solves the continuity equation \cref{eq:continuity} together with the velocity field $\phi_k$.

To move these $N$ particles according to velocity field $\phi_k(t)$, we discretize the interval by considering time step $t=0,1,2,\ldots$ and apply the following iterative scheme: for particle located at $\mathbf x_t$ at time step $t$, we generate its position $\mathbf x_{t+1}$ at time step $t+1$ by 
    \[
    \mathbf x_{t+1}=\mathbf x_t+\eta\phi_k(t, \mathbf x_t),
    \]
for some small step size $\eta>0$; that is, at each discrete time step $t$, we move the position of particle located at $\mathbf x_t$ along the direction given by $\phi_k(t,\mathbf x_t)$ by a small step. If this step size $\eta$ is small enough, this discrete sequence $\{\mathbf x_t\}_{t\in\N}$ should approximate the original continuous curve $\mathbf x(t)$. Therefore, we get the following particle algorithm (Algorithm 1).

\begin{algorithm}
\caption{Wasserstein CCGF particle algorithm}\label{alg:WCCGF}
\begin{algorithmic}[1]
\Require Input the number of iterations $T$, the step size $\eta>0$, the number $K$ of sub-populations, and $k$ initial probability distributions $\mu_k(0)$.
\State For each $k\in[K]$, randomly sample a population of $N$ particles $X_k(0)\doteq (\mathbf x_k^i(0))_{i=1}^N$ with $\mathbf x_k^i(0)\in\R^d$ according to $\mu_k(0)$. Similarly, define $X_k(t)\doteq (\mathbf x_k^i(t))_{i=1}^N$ for all $t\geq 0$.
\For {each time step $t=1,\ldots, T$}
\For {each sub-population $k=1,\ldots,K$}
    \For {each particle $i=1,\ldots, N$}
        \State Compute the velocity $\phi_k(t,\mathbf x_k^i(t))\in\R^d$ according to equation \cref{eq:WCCGF} under \cref{def:Wassertein-CCGF}.
        \State Update the particle with: 
            $\mathbf x_k^i(t+1)=\mathbf x_k^i(t)+\eta\phi_k(t,\mathbf x_k^i(t)).$
    \EndFor
\EndFor
\EndFor
\State Output $X_k(T)$ for all $k\in[K]$. 
\end{algorithmic}
\end{algorithm}

One final question remains: how can we utilize the output? We need the probability density function $\mu_k$ for each $k\in[K]$. We can approximate this density function $\mu_k$ with the population $X_k$ of particles. If one is interested in computing the integral $\int_{\R^d}f(\mathbf x)d\mu_k(\mathbf x)$ for some function $f:\R^d\to\R$ with respect to $\mu_k$, this integral can be approximated by 
\[
    \int_{\R^d}f(\mathbf x)d\mu_k(\mathbf x)\approx \frac{1}{N}\sum_{i=1}^Nf(\mathbf x_k^i).
\]
The density function value $\mu_k(\mathbf x)$ at any point $\mathbf x\in\R^d$ can also be approximated by the method of ``kernel density estimation". Let $\psi:\R^d\times\R^d\to\R$ be 
\[
    \psi(\mathbf x,\mathbf y)\doteq\frac{1}{\sqrt{2\pi}}e^{-\frac{\|\mathbf x-\mathbf y\|^2}{2}}.
\]
Then, the density $\mu_k(\mathbf x)$ can be approximated by 
\[
    \mu_k(\mathbf x)\approx \frac{1}{N}\sum_{i=1}^N\psi(\mathbf x, \mathbf x_k^i),\quad\forall \mathbf x\in\R^d.
\]

In the following example, we demonstrate how \cref{alg:WCCGF} works.
\begin{example}\label{exp:particles}
    Let $\pi\in\calP_{2,ac}(\R^2)$ be the density function of bivariate normal distribution $\mathcal{N}(\begin{bmatrix}
        0\\
        0
    \end{bmatrix}, \begin{bmatrix}
        4^2, 24\rho\\ 
        24\rho, 6^2
    \end{bmatrix})$ with $\rho=0.6$. In \cref{prob:decomposition-2}, suppose we want to decompose $\pi$ into two densities $\mu_1$ and $\mu_2$ with equal weights $p_1=p_2=\sfrac{1}{2}$. Consider the loss 
        \[
        L(\mu)=\int_{\R^2}\int_{\R^2}\|\mathbf x-\mathbf y\|^2d\mu(\mathbf x)d\mu(\mathbf y).
        \]
    We implement \cref{alg:WCCGF} to solve this instance of \cref{prob:decomposition-2} with $N=200$ particles. \cref{fig:gaussian-1} shows the movement of particles in our velocity field design. As shown in \cref{fig:gaussian-1}(a), for each $k = 1, 2$ we sample $200$ particles to represent $\mu_k$. Each colored dot is one particle and the $200$ blue (resp. orange) particles represent probability density $\mu_1$ (resp. $\mu_2$). Then, we implement \cref{alg:WCCGF} to compute the velocity of each particle (line 5 in \cref{alg:WCCGF}) and move each particle along the velocity by a short step (line 6 in \cref{alg:WCCGF}). After 100 such iterations, the 400 particles are arranged in an ellipse, as depicted in \cref{fig:gaussian-1}(b). This elliptical shape is expected since the underlying distribution $\pi$ is a bivariate normal. The contour plot of a bivariate normal density resembles an ellipse, and these points are intended to represent densities $\mu_1$ and $\mu_2$, with $\sfrac{\mu_1}{2} + \sfrac{\mu_2}{2} = \pi$. However, particles of different colors are still intermingled after 100 iterations. As we demonstrated in \cref{cor:disjoint-support}, particles of different colors should separate from each other. This indeed occurs after a sufficient number of iterations, as illustrated in \cref{fig:gaussian-1}(c). Particles are separated into two disjoint parts, each representing the support of one density $\mu_k$. The area with denser particles has a higher density value. Additionally, particles of both colors are arranged in an ellipse, illustrating the contour plot  of the bivariate normal density $\pi$. The final positions of the particles illustrate how the distribution $\pi$ is optimally decomposed into $\mu_1$ and $\mu_2$. 
\begin{figure}
  \begin{subfigure}[t]{.3\textwidth}
    \centering
    \includegraphics[width=\linewidth]{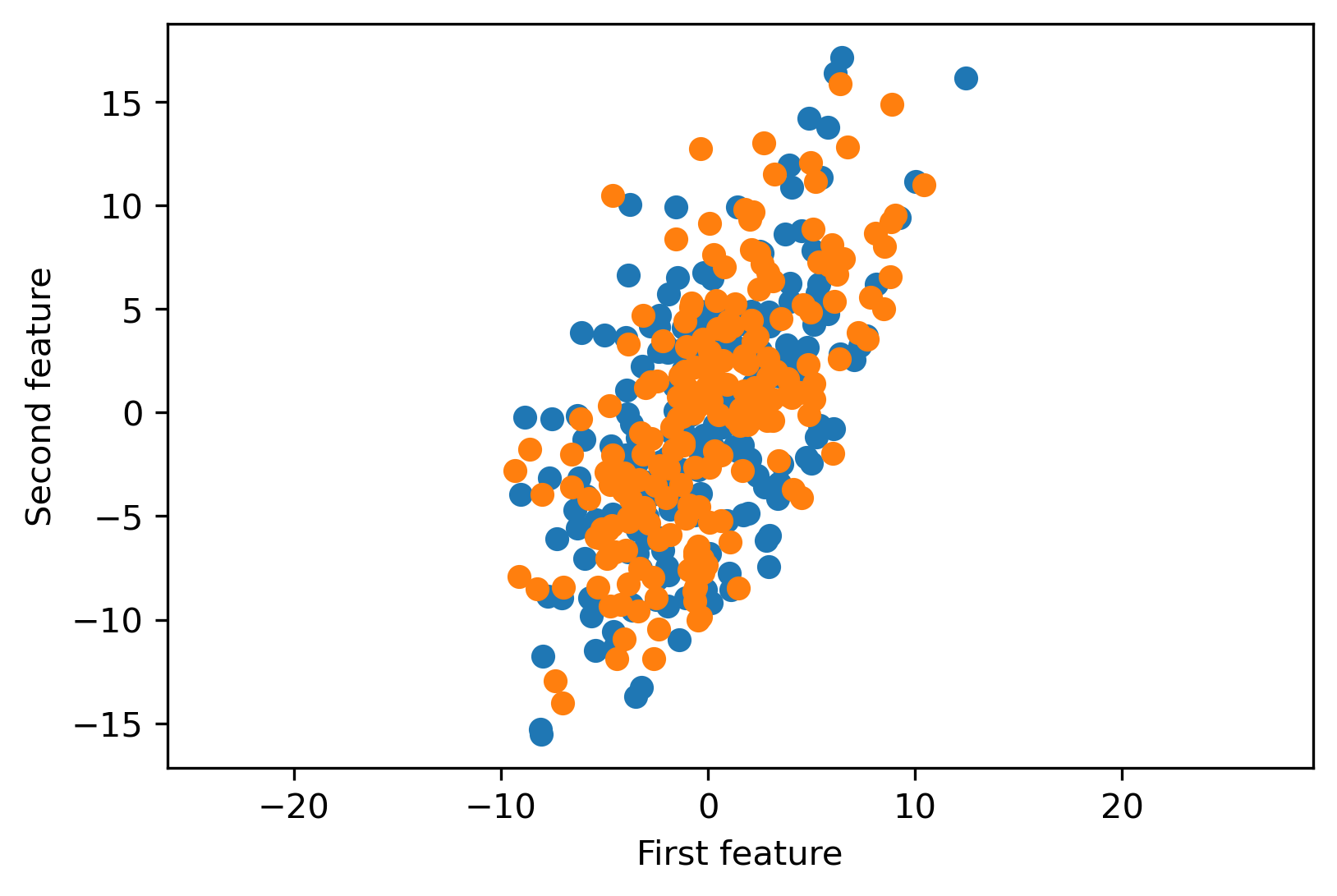}
    \caption{Initial positions of particles.}
  \end{subfigure}
  \hfill
  \begin{subfigure}[t]{.3\textwidth}
    \centering
    \includegraphics[width=\linewidth]{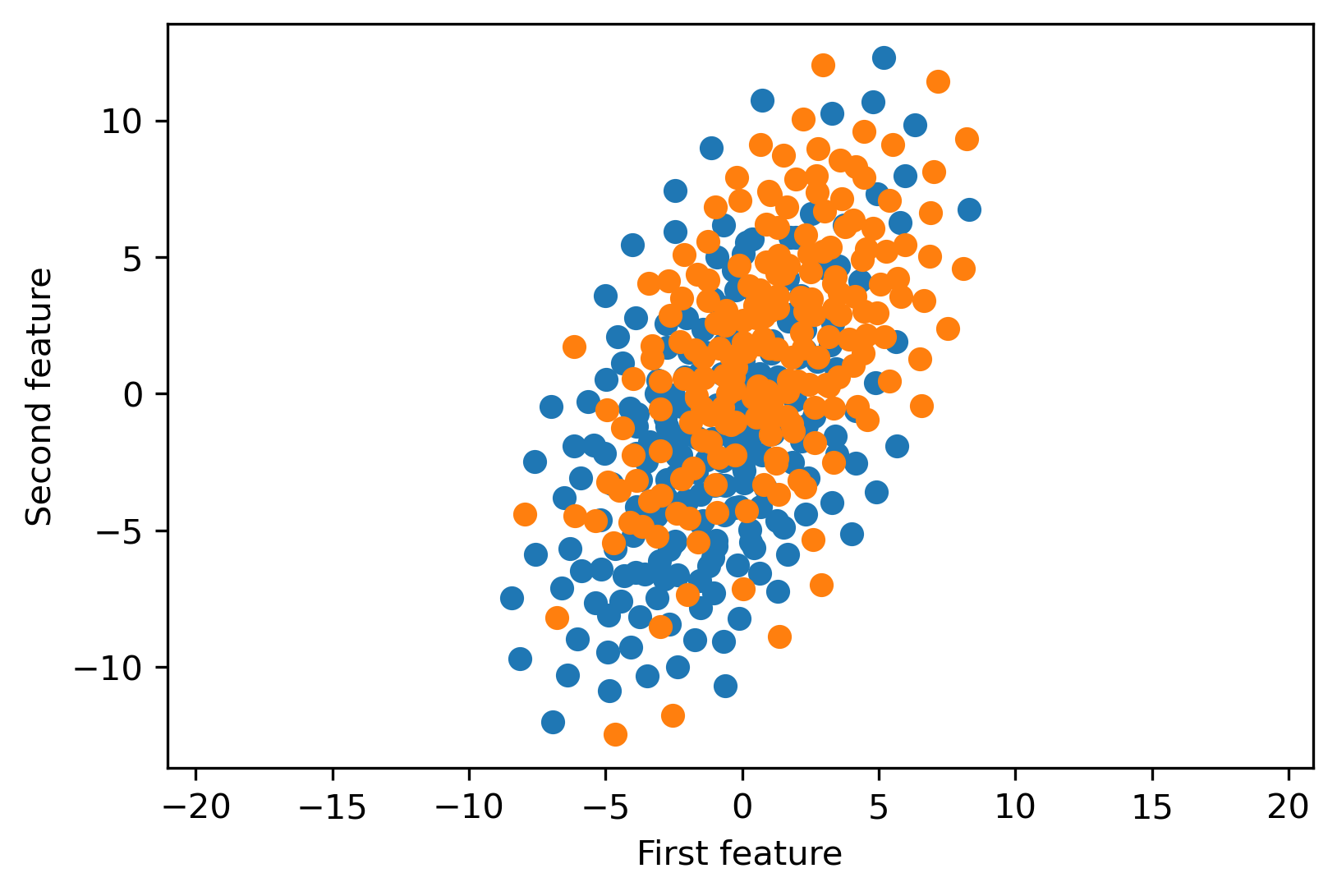}
    \caption{Particles at $t=100$.}
  \end{subfigure}
  \hfill
  \begin{subfigure}[t]{.3\textwidth}
    \centering
    \includegraphics[width=\linewidth]{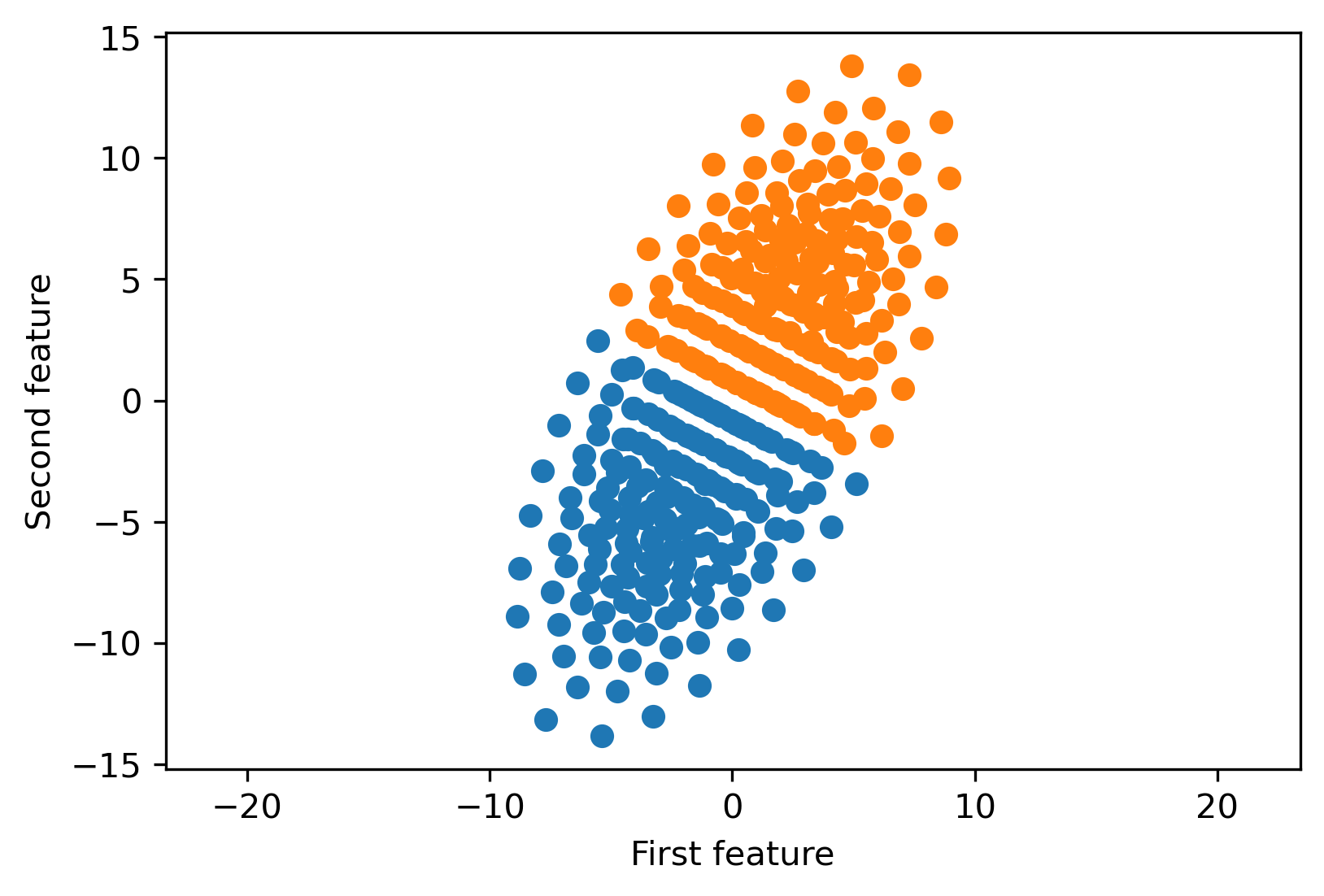}
    \caption{Final positions of particles.}
  \end{subfigure}
  \caption{Movement of particles in Wasserstein CCGF.}\label{fig:gaussian-1}
\end{figure}
\end{example}

\begin{example}[Multiple classes]\label{exp:multi-class}
  We can also decompose the probability density $\pi\in\calP_{2,ac}(\R^2)$ into multiple ($K > 2$) probability densities. Let $\pi\in\calP_{2,ac}(\R^2)$ be the density function of the bivariate normal distribution as given in \cref{exp:particles}. Let $\sigma_{X_1}$, $\sigma_{X_2}$ and $\rho$ be the standard deviation of feature one, feature two, and the correlation between feature one and two, respectively. Let $\mu_{X_1}=\mu_{X_2}=0$. In \cref{prob:decomposition-2}, suppose we want to decompose $\pi$ into $K$ densities $\mu_1,\ldots,\mu_K$ with equal weights $p_k=\sfrac{1}{K}$ for all $k\in[K]$. Consider the loss function $L$ given by \cref{exp:clustering} with $W=\begin{bmatrix}
      1 & 0 \\
      0 & 1
  \end{bmatrix}$. \cref{fig:gaussian-3} shows the decomposition found by \cref{alg:WCCGF} under different setups. Again, each population of particles of one color represents one density $\mu_k$ for $k\in[K]$. The positions of particles of one color resemble the contour plot of the corresponding density function; that is, an area with denser particles has a larger density value. According to \cref{cor:disjoint-support}, the populations of particles of different colors should be disjoint, as illustrated in \cref{fig:gaussian-3}. Additionally, we observe that the optimal decompositions may not always appear intuitive. For example, the decompositions shown in \cref{fig:gaussian-3}(b), (c), and (f) are not immediately obvious. By contrast, \cref{fig:gaussian-3}(a), (d), and (e) are more natural-looking partitions of the space due to their simple geometric nature.
  
   \begin{figure}
  \begin{subfigure}[t]{.3\textwidth}
    \centering
    \includegraphics[width=\linewidth]{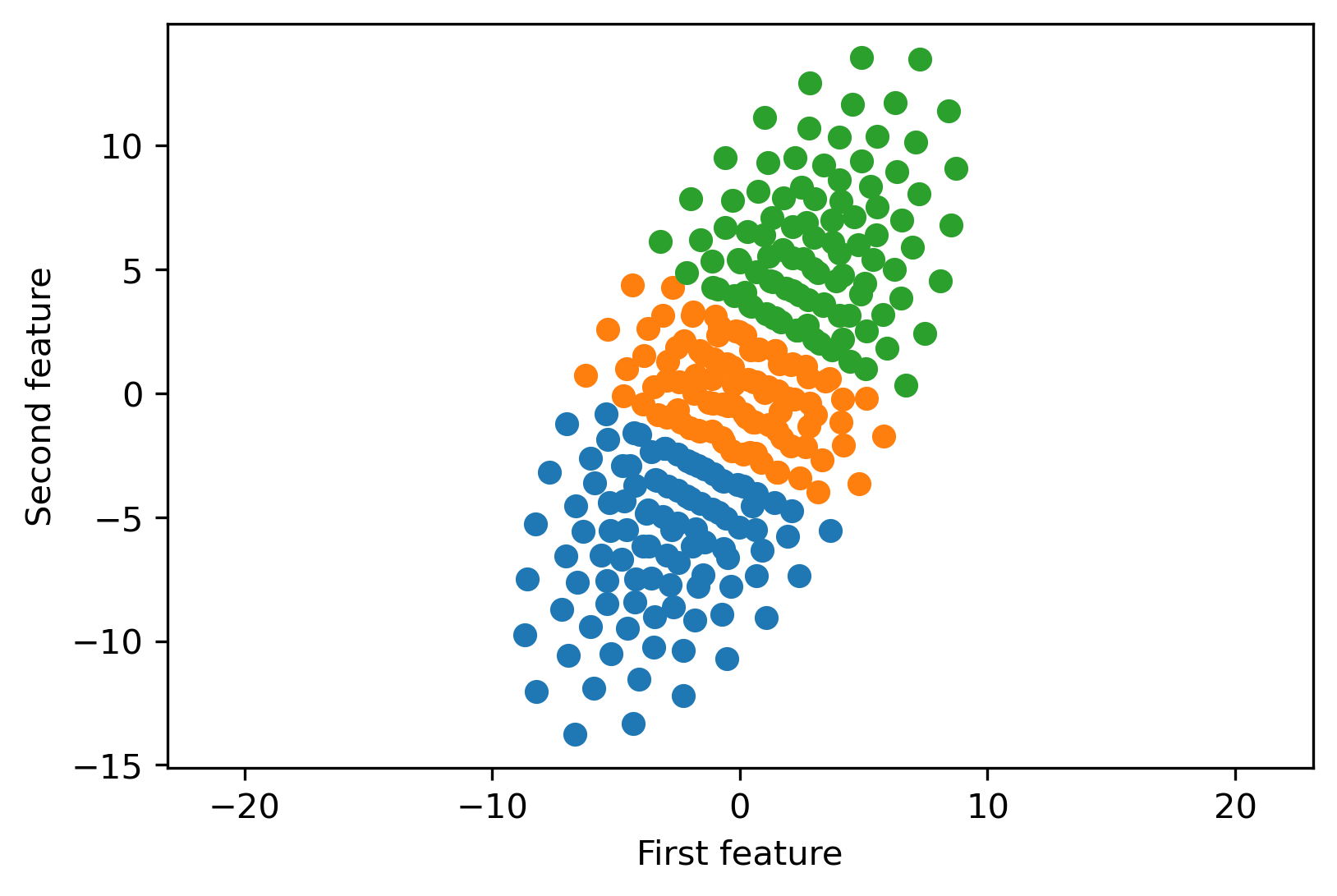}
    \caption{$\sigma_{X_1}=4,\sigma_{X_2}=6,\rho=0.6$ and $K=3$}
  \end{subfigure}
    \hfill
  \begin{subfigure}[t]{.3\textwidth}
    \centering
    \includegraphics[width=\linewidth]{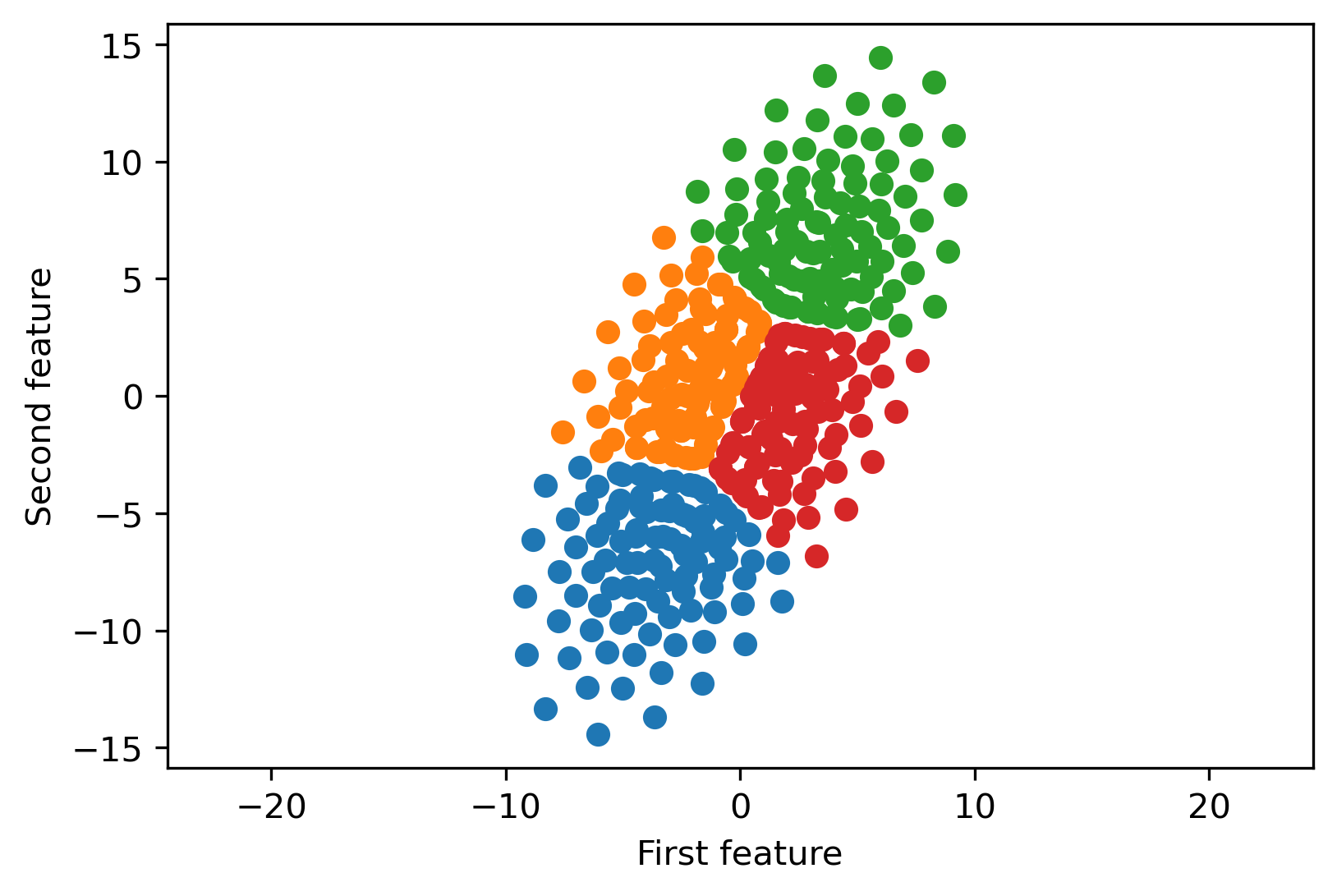}
     \caption{$\sigma_{X_1}=4,\sigma_{X_2}=6,\rho=0.6$ and $K=4$}
  \end{subfigure}
  \hfill
  \begin{subfigure}[t]{.3\textwidth}
    \centering
    \includegraphics[width=\linewidth]{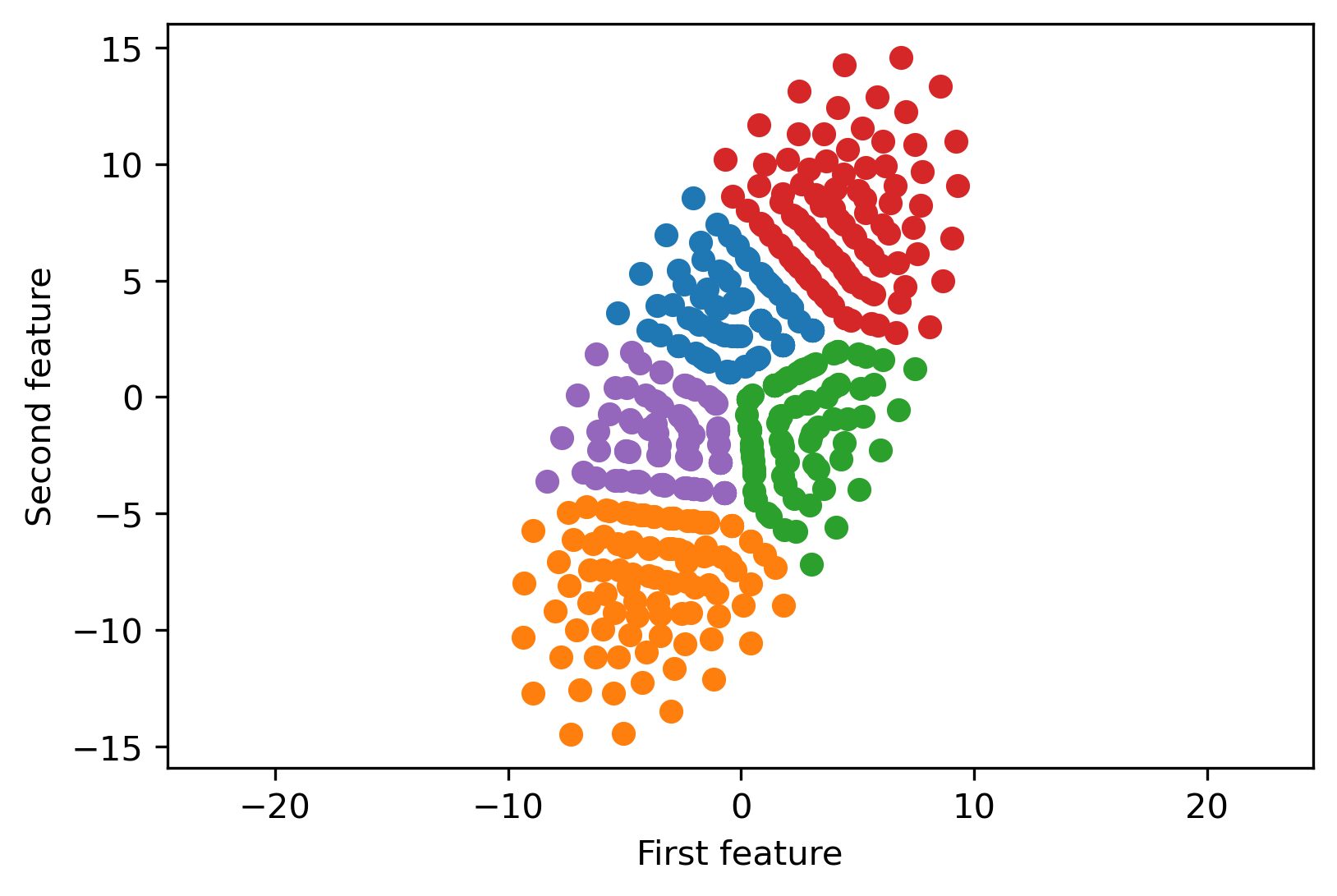}
     \caption{$\sigma_{X_1}=4,\sigma_{X_2}=6,\rho=0.6$ and $K=5$}
  \end{subfigure}
    \hfill
  \begin{subfigure}[t]{.3\textwidth}
    \centering
    \includegraphics[width=\linewidth]{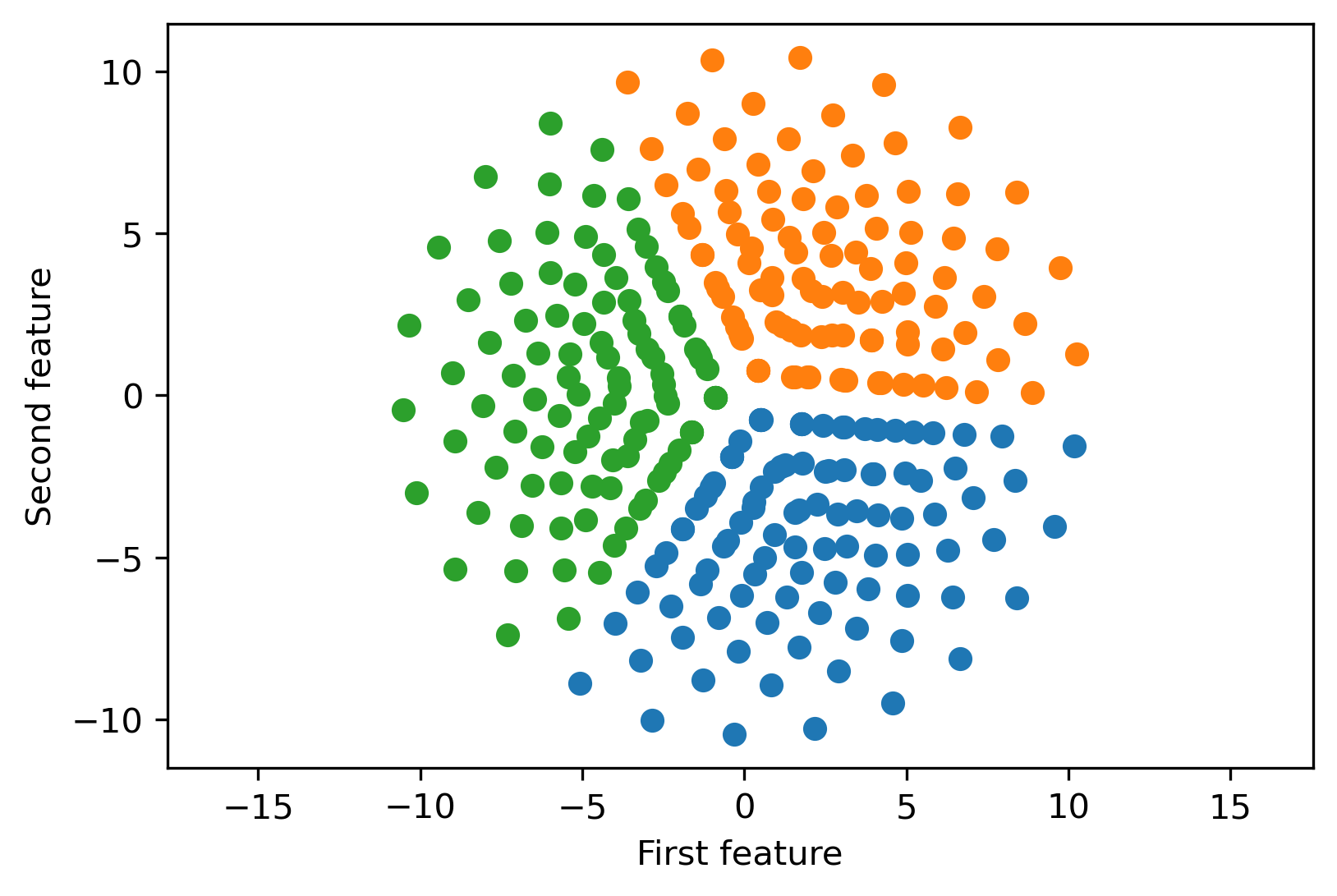}
     \caption{$\sigma_{X_1}=5,\sigma_{X_2}=5,\rho=0$ and $K=3$}
  \end{subfigure}
    \hfill
  \begin{subfigure}[t]{.3\textwidth}
    \centering
    \includegraphics[width=\linewidth]{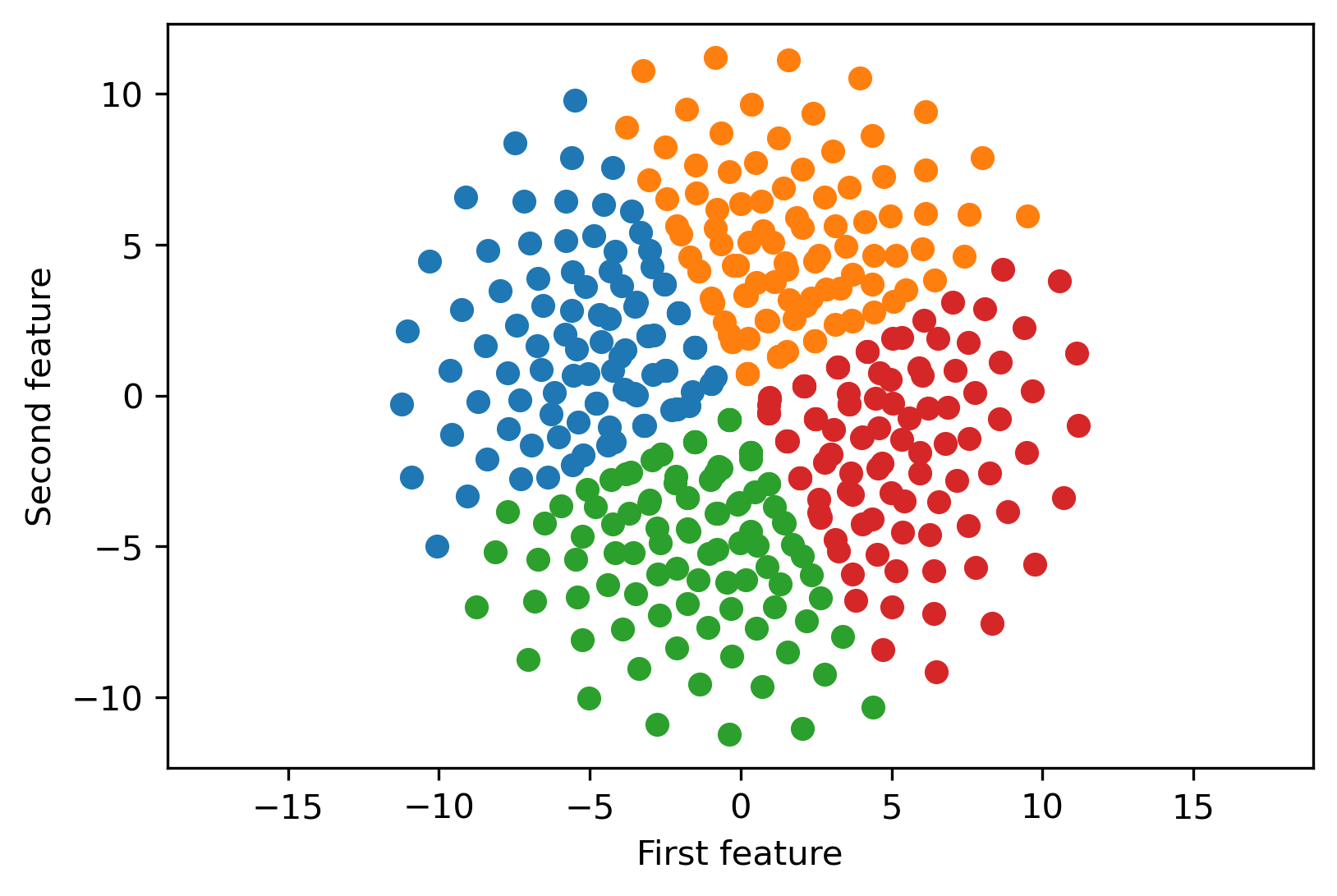}
    \caption{$\sigma_{X_1}=5,\sigma_{X_2}=5,\rho=0$ and $K=4$}
  \end{subfigure}
  \hfill
  \begin{subfigure}[t]{.3\textwidth}
    \centering
    \includegraphics[width=\linewidth]{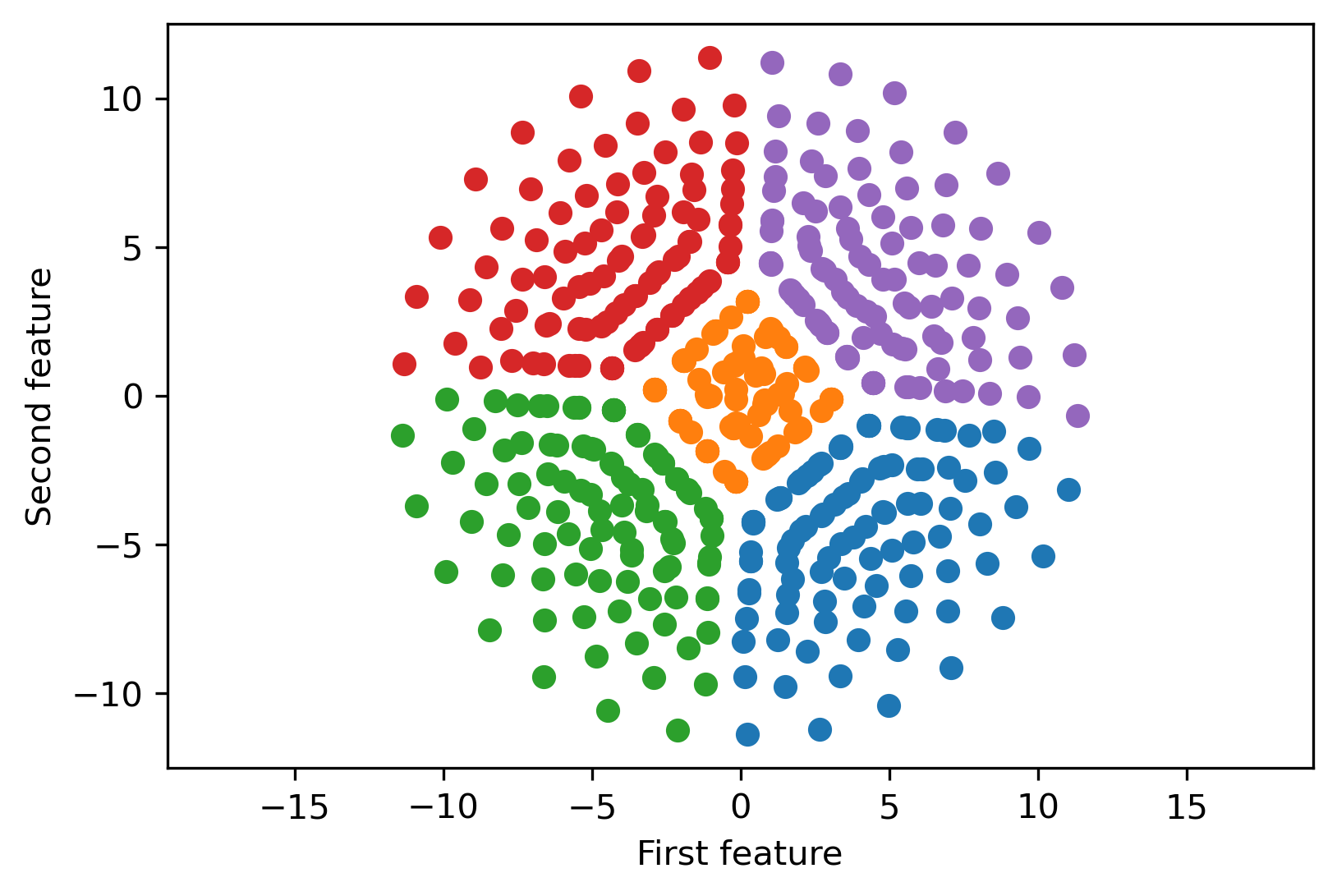}
    \caption{$\sigma_{X_1}=5,\sigma_{X_2}=5,\rho=0$ and $K=5$}
  \end{subfigure}
  \caption{Decomposition with multiple groups.}\label{fig:gaussian-3}
\end{figure}
\end{example}

We use the following natural decompositions as benchmarks to show the efficacy of the decompositions discovered by \cref{alg:WCCGF}.

\begin{example}[Comparision with natural decomposition]\label{exp:com}
    Consider the settings described in \cref{exp:multi-class}. We consider the following simple ways (\cref{fig:com}) using parallel slices to decompose the underlying distribution $\pi$ into $K$ densities $\mu_1,\ldots,\mu_K$. Recall that we require $\sum_{k=1}^K\sfrac{\mu_k}{K}=\pi$. Compared with the alternative decompositions shown in \cref{fig:com}, the decompositions discovered by \cref{alg:WCCGF} shown in \cref{fig:gaussian-3} improves the objective value (i.e., \cref{prob:decomposition-2} with variance loss described in \cref{exp:clustering}) by $7.90\%, 12.33\%$ and $47.67\%$, respectively. 

    \begin{figure}
  \begin{subfigure}[t]{.32\textwidth}
    \centering
    \includegraphics[width=\linewidth]{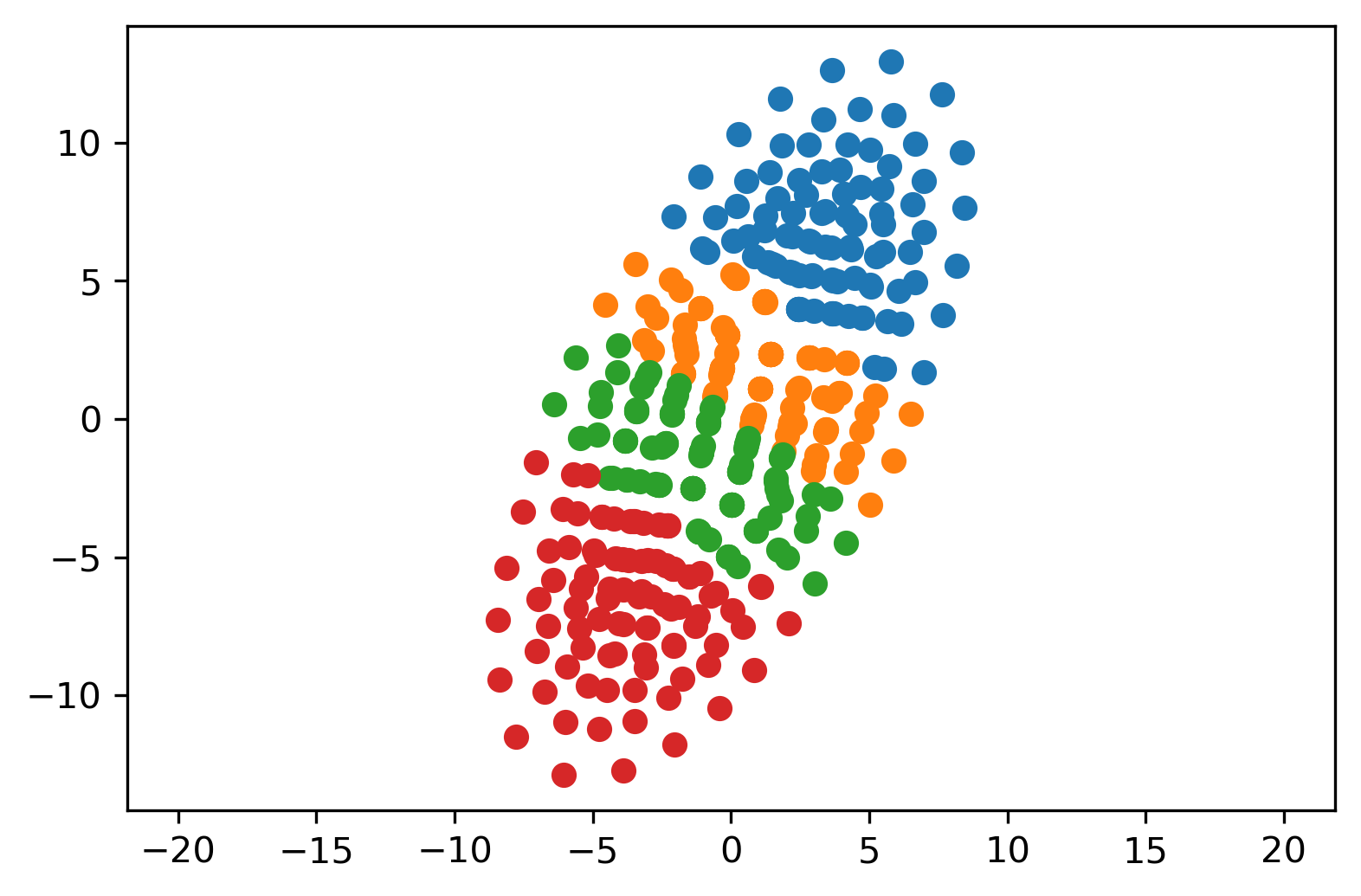}
    \caption{$\sigma_{X_1}=4,\sigma_{X_2}=6,\rho=0.6$ and $K=4$}
  \end{subfigure}
  \hfill
  \begin{subfigure}[t]{.32\textwidth}
    \centering
    \includegraphics[width=\linewidth]{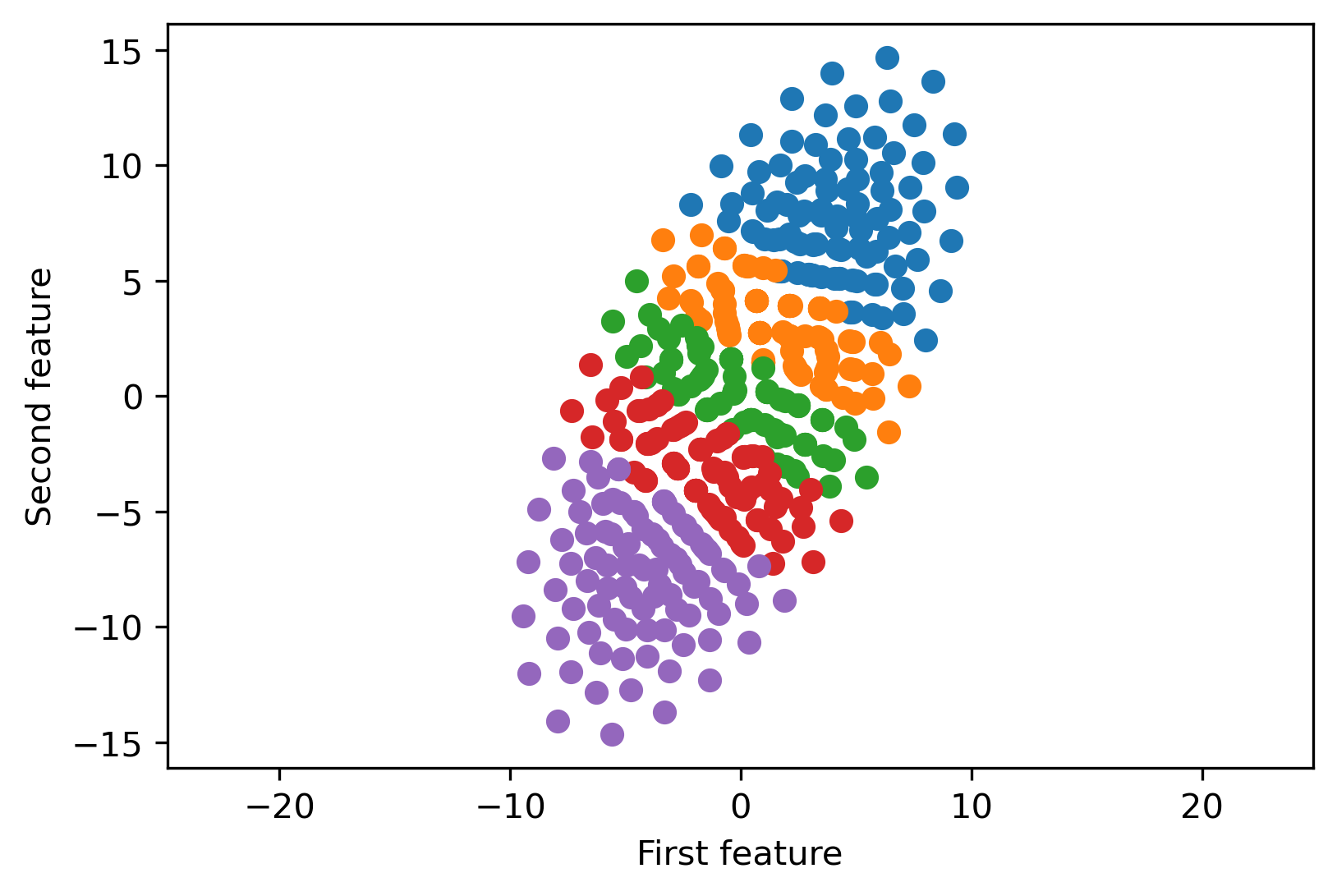}
    \caption{$\sigma_{X_1}=4,\sigma_{X_2}=6,\rho=0.6$ and $K=5$}
  \end{subfigure}
  \hfill
   \begin{subfigure}[t]{.32\textwidth}
    \centering
    \includegraphics[width=\linewidth]{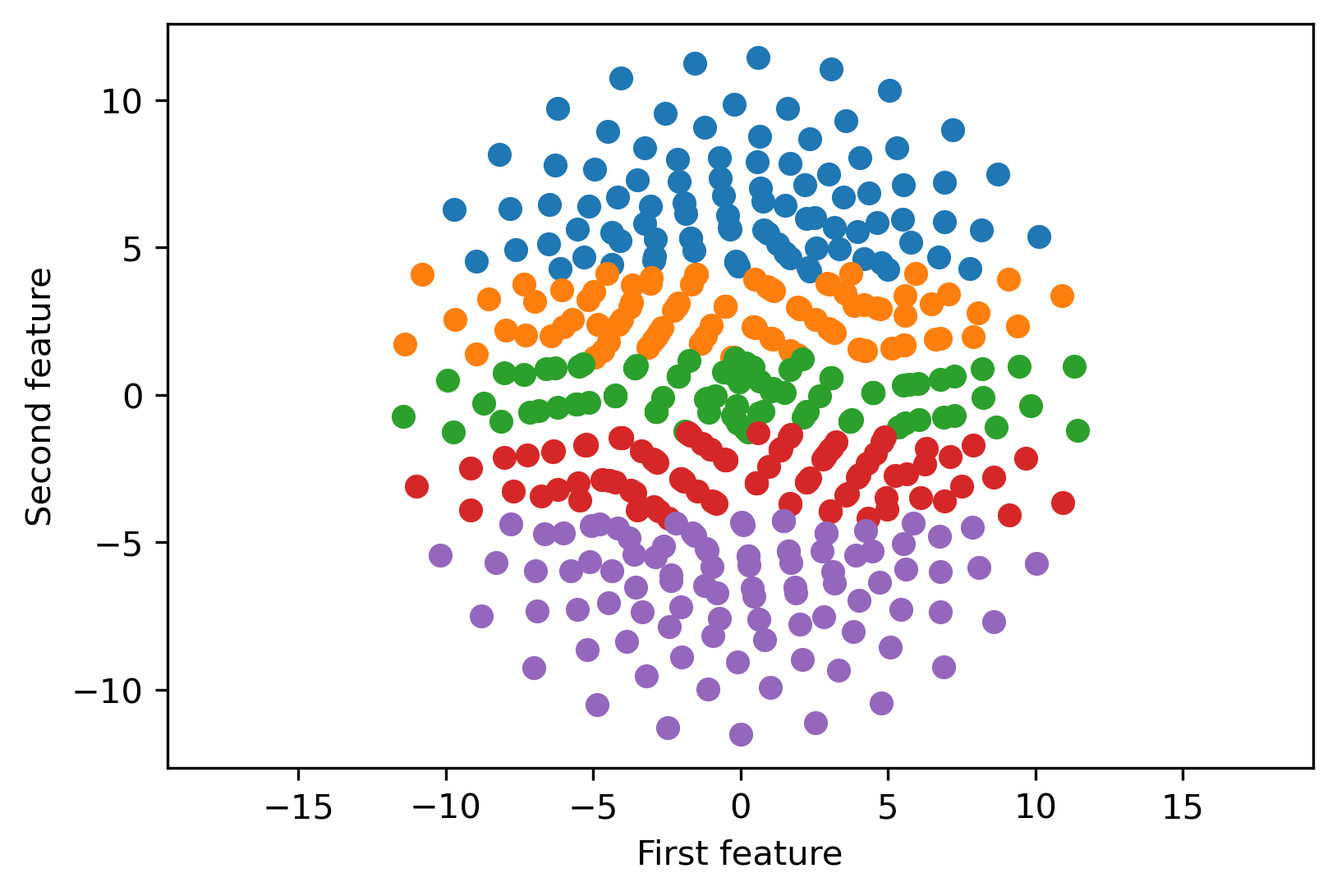}
    \caption{$\sigma_{X_1}=5,\sigma_{X_2}=5,\rho=0$ and $K=5$}
  \end{subfigure}
  \caption{Results of \cref{exp:com}}\label{fig:com}
\end{figure}
\end{example}

In certain scenarios, decision-makers need to divide the users into several groups that are homogeneous in one feature and diverse in other features. For example, one might want to group people with similar interests but have diverse demographic backgrounds.

\begin{example}\label{exp:multi-objective}
We can also maximize the similarity of feature one and maximize the diversity of feature two simultaneously. This can be achieved by considering the loss function $L$ given by \cref{exp:clustering} with $W=\begin{bmatrix}
      1 & 0 \\
      0 & -1
  \end{bmatrix}$. In this case, our algorithm minimizes the variance of feature one minus the variance of feature two. Let $\pi\in\calP_{2,ac}(\R^2)$ be the density function of a mixed bivariate normal distribution, i.e., $\pi=\sfrac{1}{3}\pi_1+\sfrac{1}{3}\pi_2+\sfrac{1}{3}\pi_3$, where $\pi_1,\pi_2,\pi_3$ are densities of normal distributions with parameters $(\mu_{X_1}, \mu_{X_1}, \sigma_{X_1},\sigma_{X_1},\rho)=(-10,5,6,5,-0.7), (0,10,7,5,0.8)$ and $(10, 15, 8, 5, -0.9)$, respectively. \cref{fig:gaussian-2}(a) shows particles sampled from $\pi$. \cref{fig:gaussian-2}(b) shows the decomposition generated by \cref{alg:WCCGF}. The result in \cref{fig:gaussian-2}(b) is not surprising since, with this weight matrix $W$, we aim to minimize the dissimilarity of the first feature and maximize the diversity of the second feature. Notice that the regions in \cref{fig:gaussian-2}(b) are ``narrow'' horizontally but ``tall'' vertically. For comparison, \cref{fig:gaussian-2}(c) shows the decomposition when $W$ is the identity matrix, which aims to define regions that are both ``narrow'' and ``short''. In this case, we minimize the variance of features one and two. 
    
\begin{figure}
  \begin{subfigure}[t]{.32\textwidth}
    \centering
    \includegraphics[width=\linewidth]{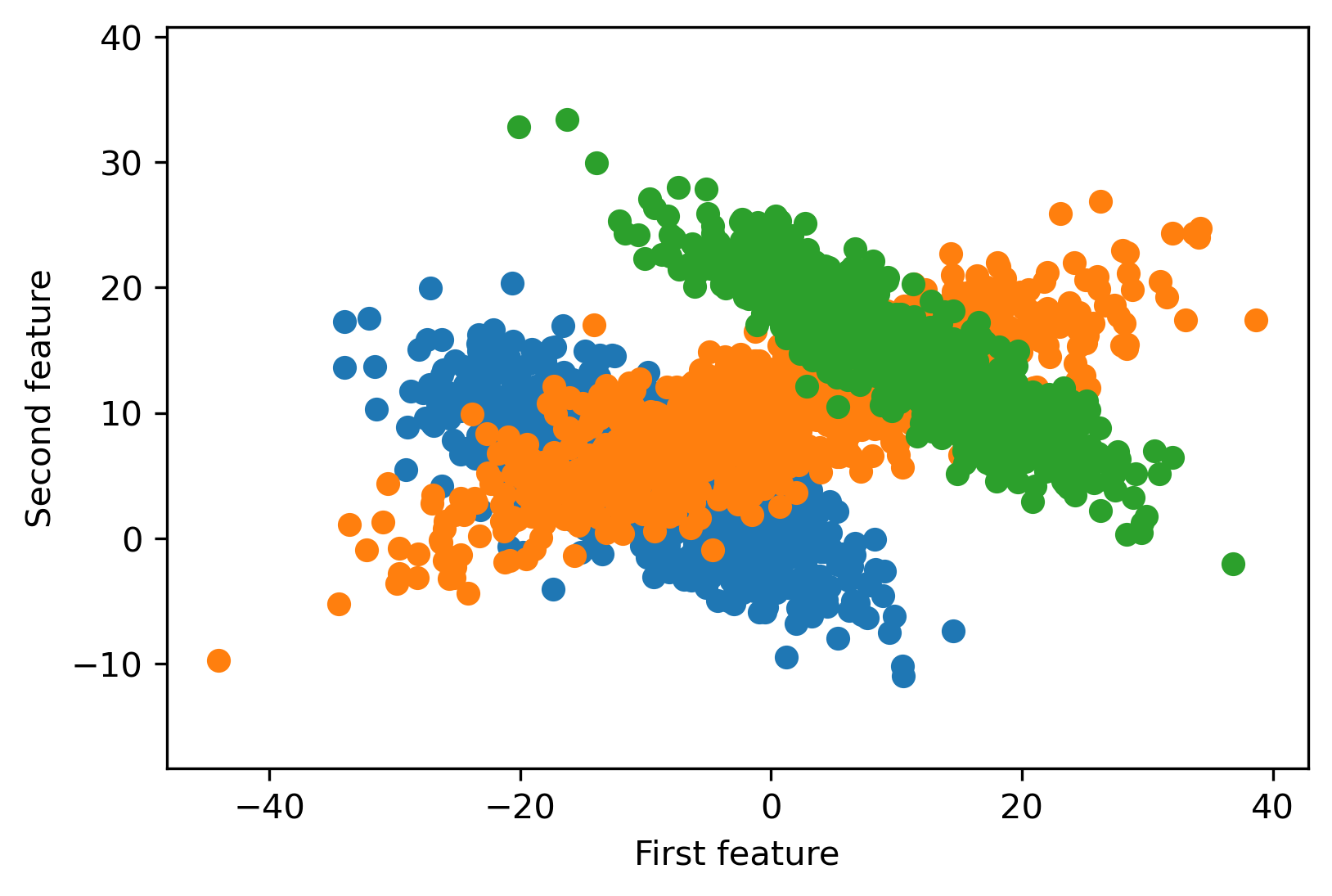}
    \caption{Particles sampled from $\pi$.}
  \end{subfigure}
    \hfill
  \begin{subfigure}[t]{.32\textwidth}
    \centering
    \includegraphics[width=\linewidth]{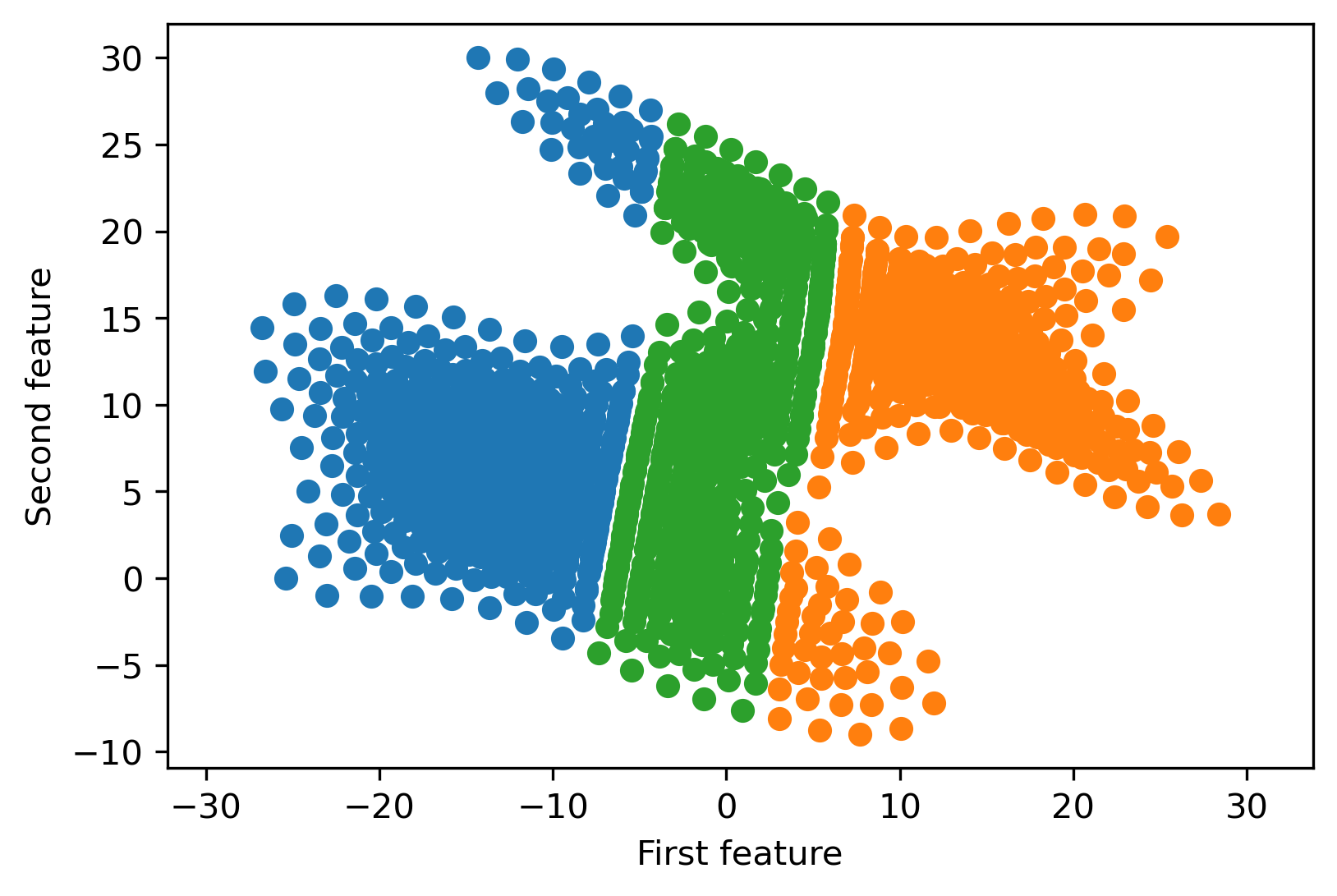}
    \caption{Decomposition with $W$.}
  \end{subfigure}
  \hfill
  \begin{subfigure}[t]{.32\textwidth}
    \centering
    \includegraphics[width=\linewidth]{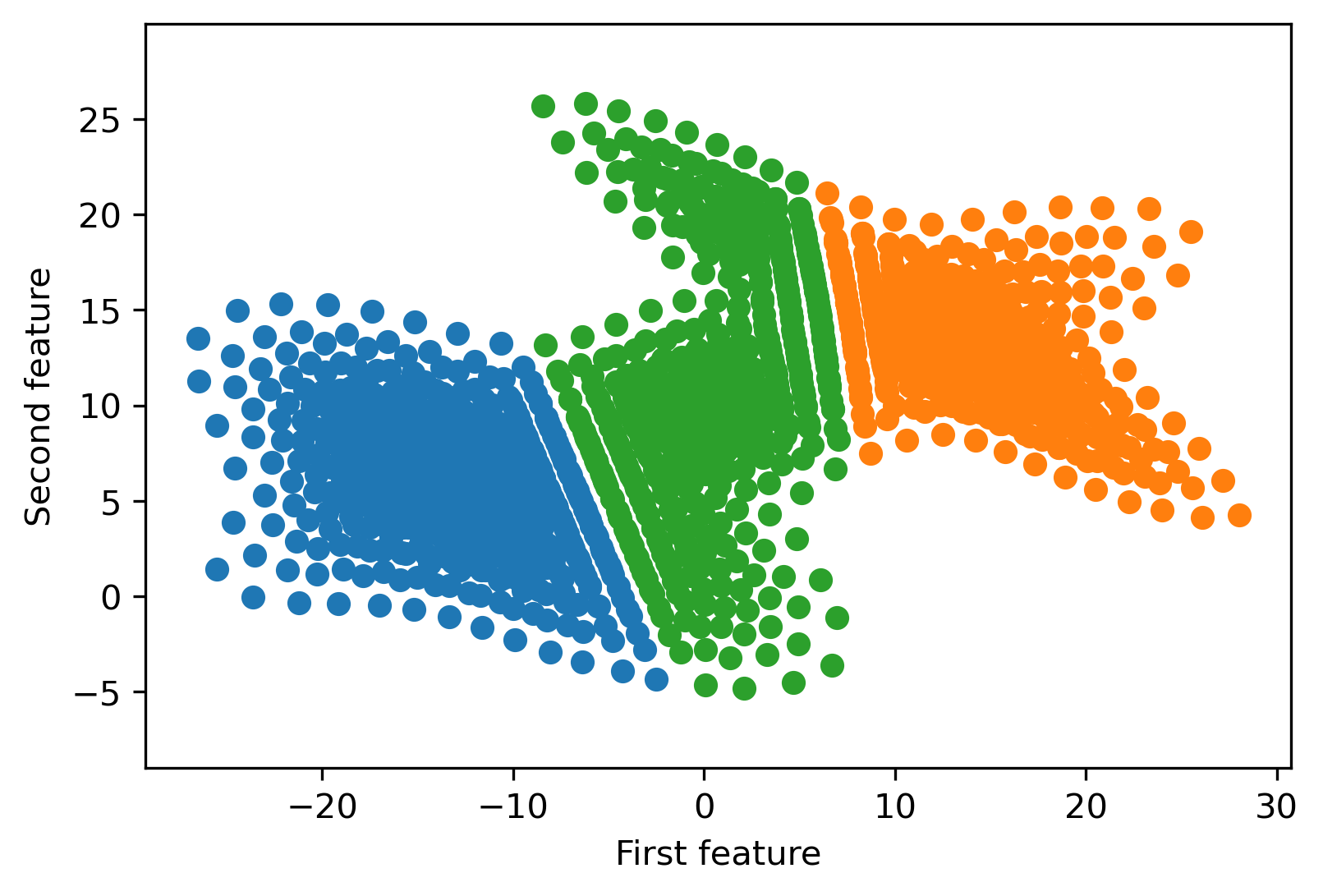}
    \caption{Decomposition when $W$ is the identity matrix.}
  \end{subfigure}
  \caption{Results of \cref{exp:multi-objective}}\label{fig:gaussian-2}
\end{figure}
\end{example}

\subsection{Case study: League design with Elo scores}\label{ss:case-study-elo}

In this subsection, we explore the league design problem in \cref{exp:elo}. Recall that the distribution loss function $L$ for this problem is defined in \cref{eq:elo-loss}.

\begin{example}[$\pi$ is uniform distribution]\label{exp:uniform}
    Let $\pi\in\calP_{2,ac}(\R)$ be the density function of a uniform distribution over interval $[10, 30]$. We implement \cref{alg:WCCGF} to solve this instance of \cref{prob:decomposition-2} with $N=200$ particles. Our benchmark is the ``Grand League" design, where all players are placed in the same league. This essentially means there is no specific league design in place. \cref{fig:uniform}(a) shows the final positions of particles given by \cref{alg:WCCGF}. Blue (resp. orange) particles in league 1 (resp. league 2) represent the density $\mu_1$ (resp. $\mu_2$). \cref{fig:uniform}(b) shows the histogram of $\mu_1$ (in blue) and $\mu_2$ (in orange), which are basically the densities of uniform distributions over $[10,20]$ and $[20,30]$, respectively. This means the optimal league design is to have a ``novice" league (for players with skill levels in $[10,20]$) and a ``veteran" league (for players with skill levels in $[20,30]$). This result is not surprising. Intuitively, since $\pi$ is uniform distribution and we require $p_1=p_2=\sfrac{1}{2}$, to minimize the Elo loss, $\mu_1$ and $\mu_2$ should also be uniform distributions supported on intervals with length $10$. \cref{fig:uniform}(c) shows the average win rate of each skill level. The green curve represents the average win rate of each skill level under the ``Grand League" design; that is, for players with skill level $x$, the average win rate is $\int_{y\in[10,30]}\frac{x}{x+y}d\pi(y)$. The blue curve represents the average win rate of each skill level in league 1 (veteran league), i.e., for players with skill level $x\in[20,30]$, the average win rate is $\int_{y\in[20,30]}\frac{x}{x+y}d\mu_1(y)$. The orange curve represents the average win rate of each skill level in league 2 (novice league), i.e., for players with skill level $x\in[10,20]$, the average win rate is $\int_{y\in[10,20]}\frac{x}{x+y}d\mu_2(y)$. Compared with the ``Grand League" design, players in the novice league have higher win rates, and players in the veteran league have lower win rates. This is generally considered desirable as it reduces frustration for novice players and increases the challenge for veteran players. Particularly, notice that players in the veteran league have roughly the same win rates as the ``high-end" players (in $[12.5, 20]$) in the novice league. 
\begin{figure}
  \begin{subfigure}[t]{.32\textwidth}
    \centering
    \includegraphics[width=\linewidth]{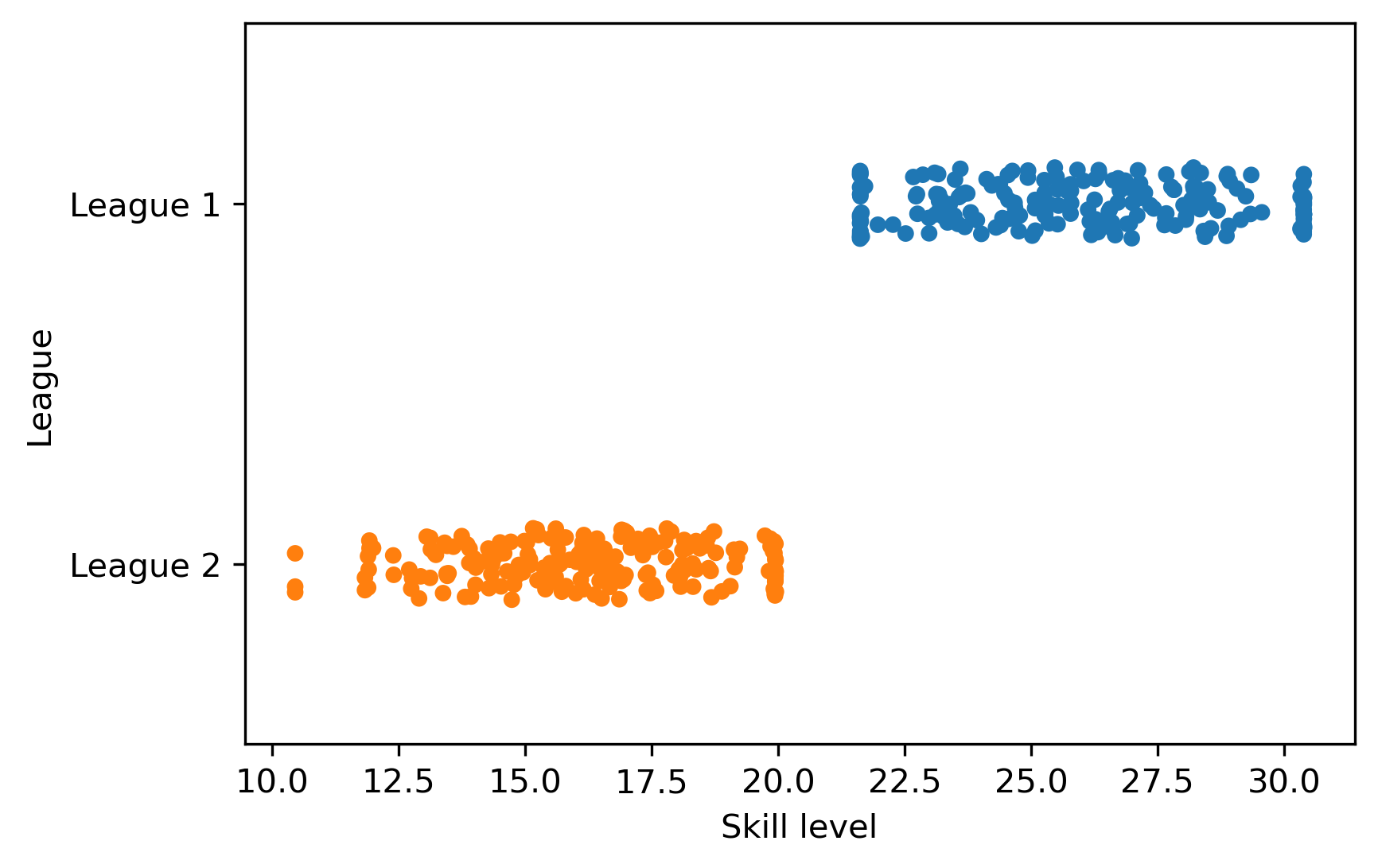}
    \caption{Final positions of particles}
  \end{subfigure}
\hfill
 \begin{subfigure}[t]{.32\textwidth}
    \centering
    \includegraphics[width=\linewidth]{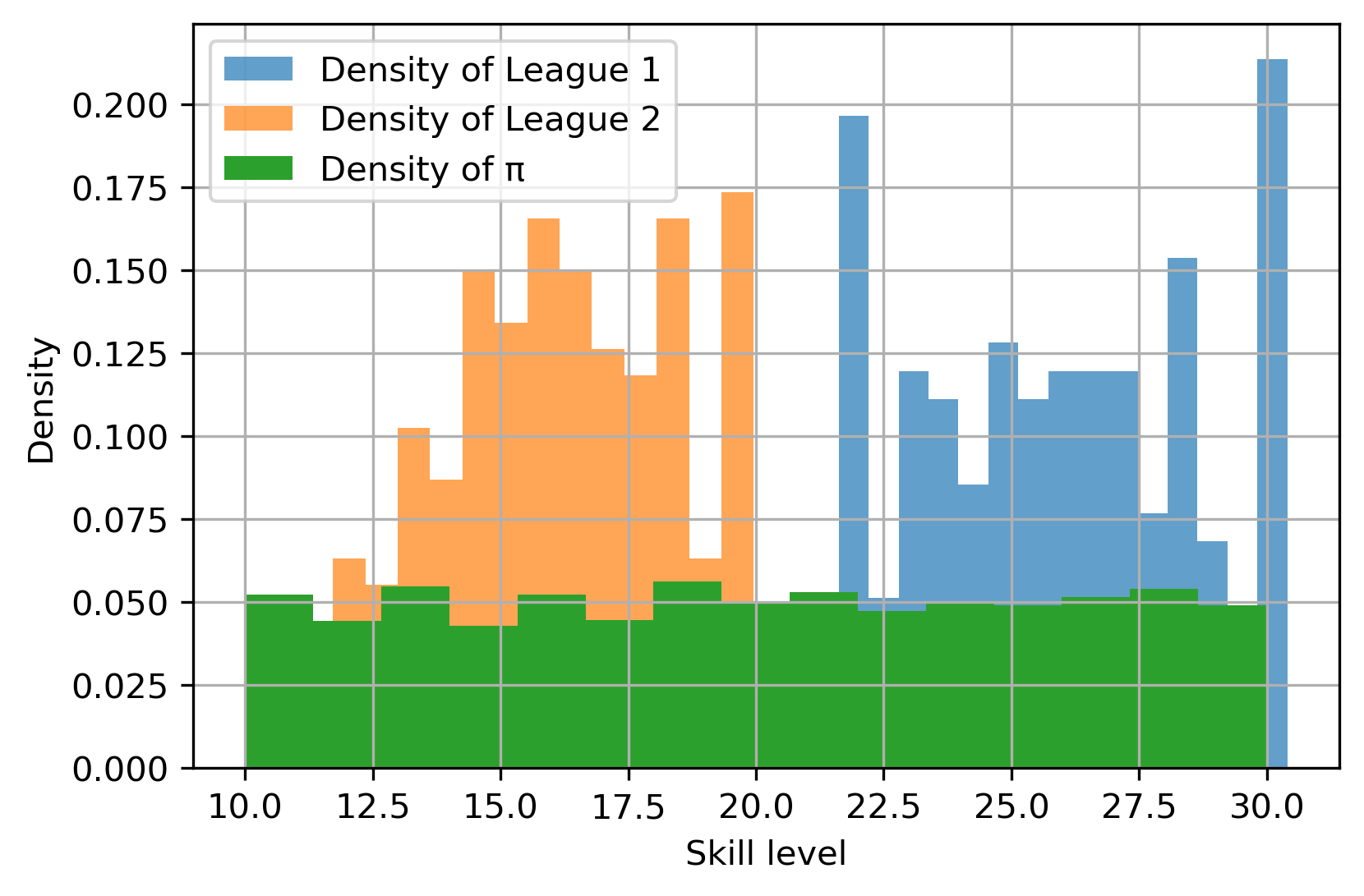}
    \caption{Histogram of $\mu_1,\mu_2$ and $\pi$}
  \end{subfigure}
\hfill
  \begin{subfigure}[t]{.32\textwidth}
    \centering
    \includegraphics[width=\linewidth]{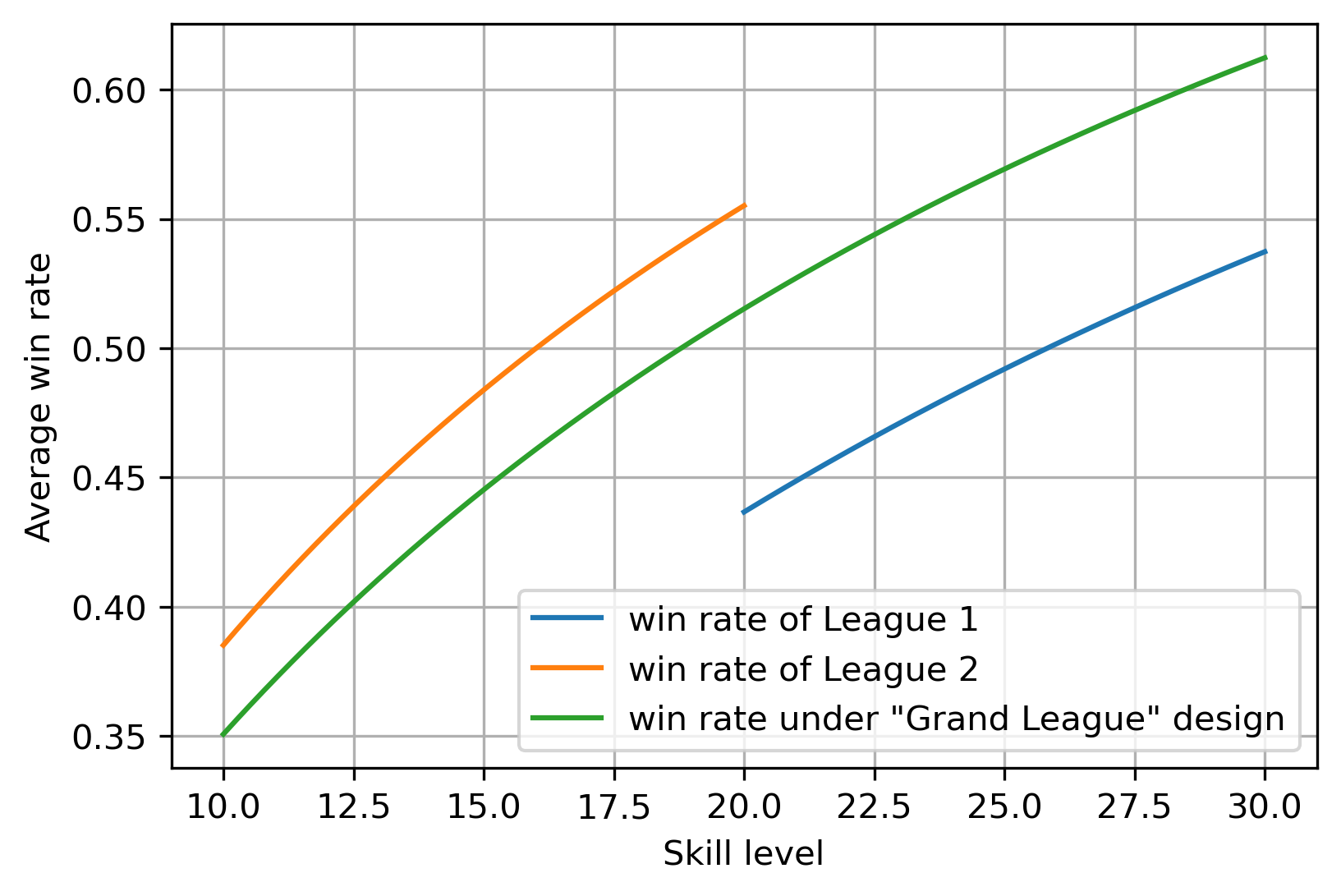}
    \caption{Average win rate}
  \end{subfigure}
  \caption{Results of \cref{exp:uniform}}\label{fig:uniform}
\end{figure}
\end{example}

\begin{example}[$\pi$ is lognormal distribution]\label{exp:lognormal}
    Let $\pi\in\calP_{2,ac}(\R)$ be the density function of a lognormal distribution with parameter $(4, 0.5^2)$. As shown in the \cref{fig:lognormal}(b), the density is single-peaked, indicated by the green color. We implement \cref{alg:WCCGF} to solve this instance of \cref{prob:decomposition-2} with $N=200$ particles. Again, we use the ``Grand League" design as our benchmark. \cref{fig:lognormal}(a) shows the final positions of particles given by \cref{alg:WCCGF}. Blue (resp. orange) particles in league 1 (resp. league 2) represent the density $\mu_1$ (resp. $\mu_2$). \cref{fig:lognormal}(b) shows the histogram of $\mu_1$ (in blue) and $\mu_2$ (in orange). The optimal league design still consists of a novice league (league 2) and a veteran league (league 1). \cref{fig:lognormal}(c) shows the average win rate of each skill level, which are computed in the same way as shown in \cref{exp:uniform}. Compared with the ``Grand League" design, players in the novice league have higher win rates, and players in the veteran league have lower win rates.
\begin{figure}
  \begin{subfigure}[t]{.32\textwidth}
    \centering
    \includegraphics[width=\linewidth]{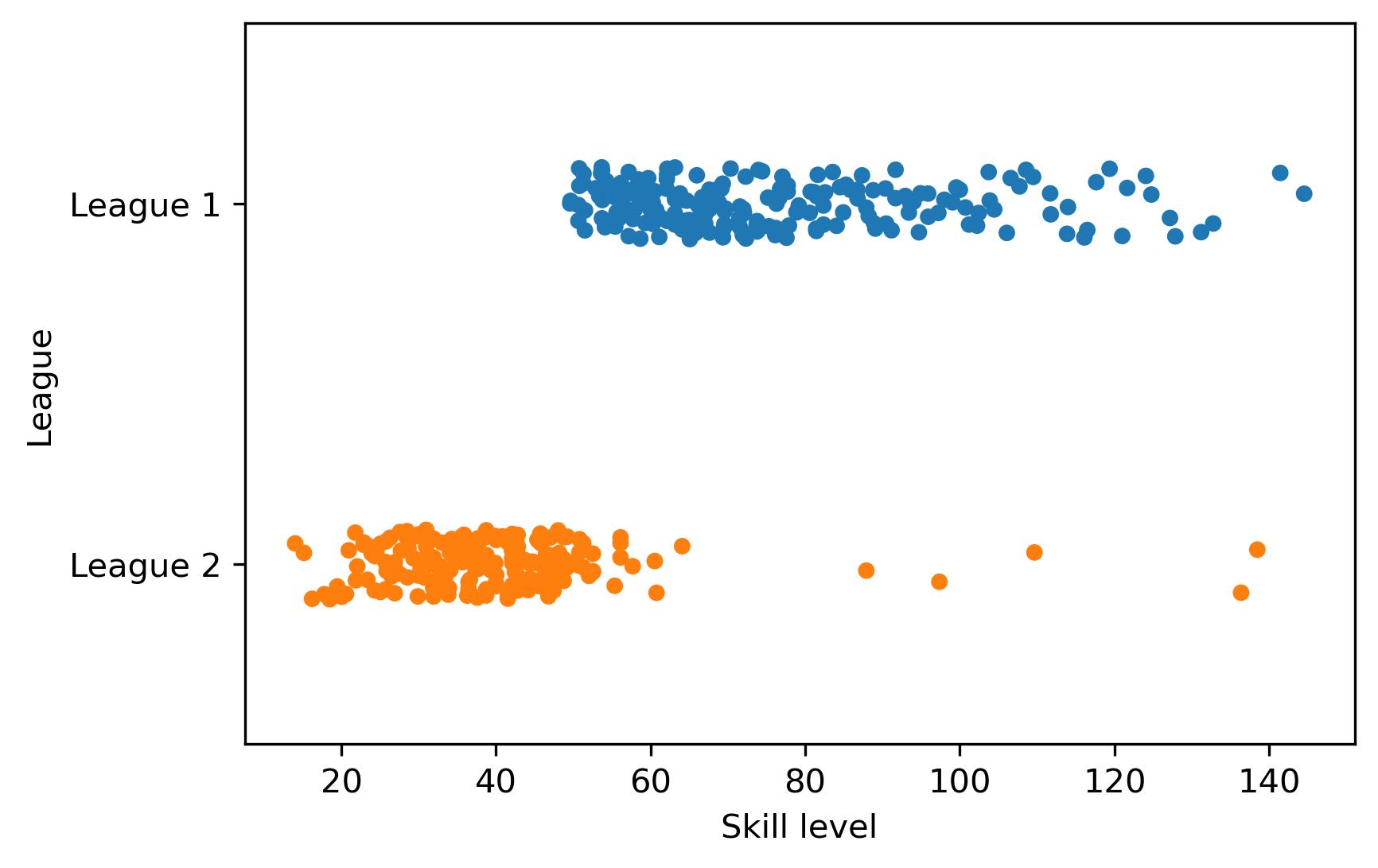}
     \caption{Final positions of particles}
  \end{subfigure}
    \hfill
      \begin{subfigure}[t]{.32\textwidth}
    \centering
    \includegraphics[width=\linewidth]{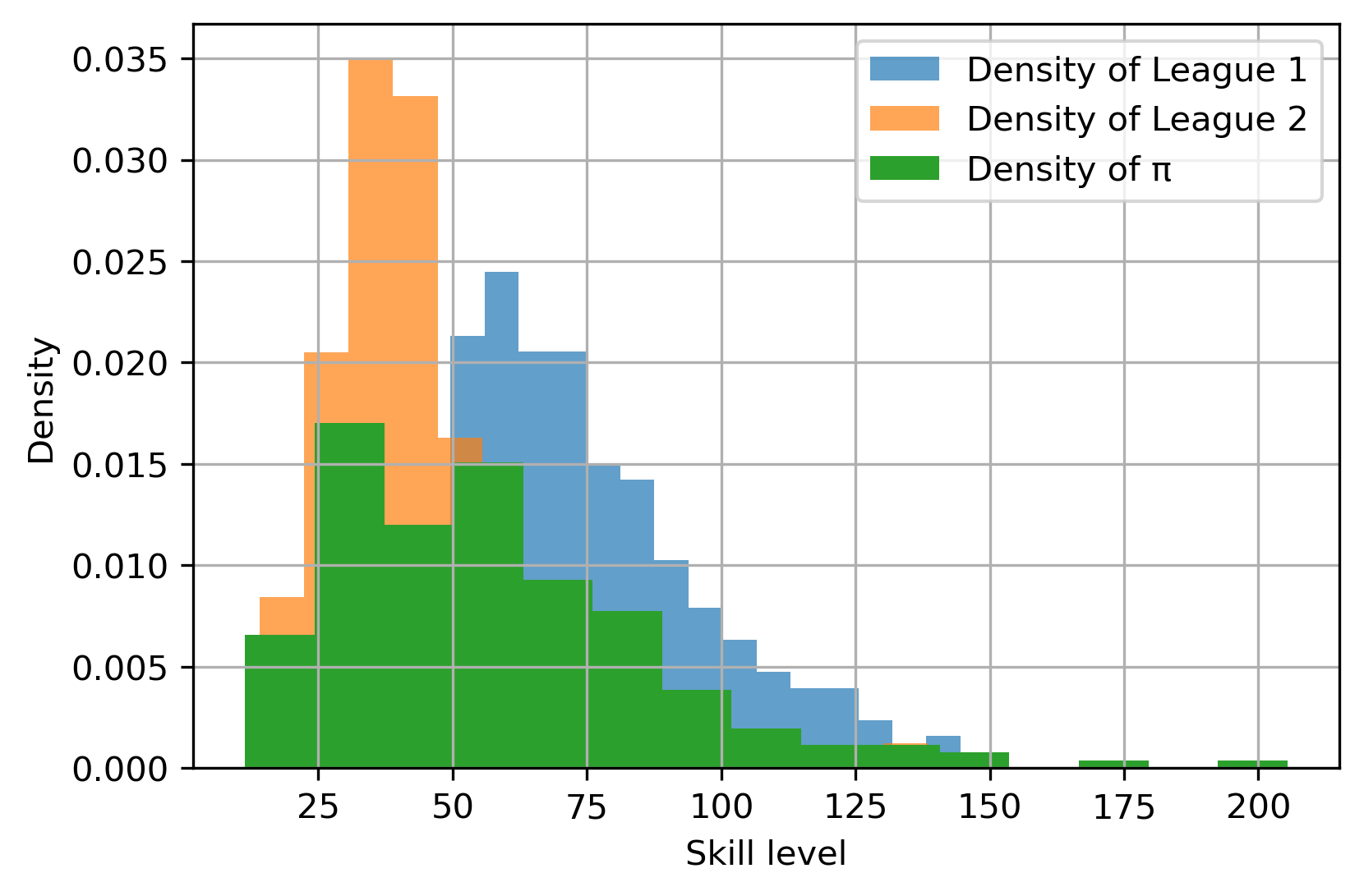}
    \caption{Histogram of $\mu_1,\mu_2$ and $\pi$}
  \end{subfigure}
  \hfill
  \begin{subfigure}[t]{.32\textwidth}
    \centering
    \includegraphics[width=\linewidth]{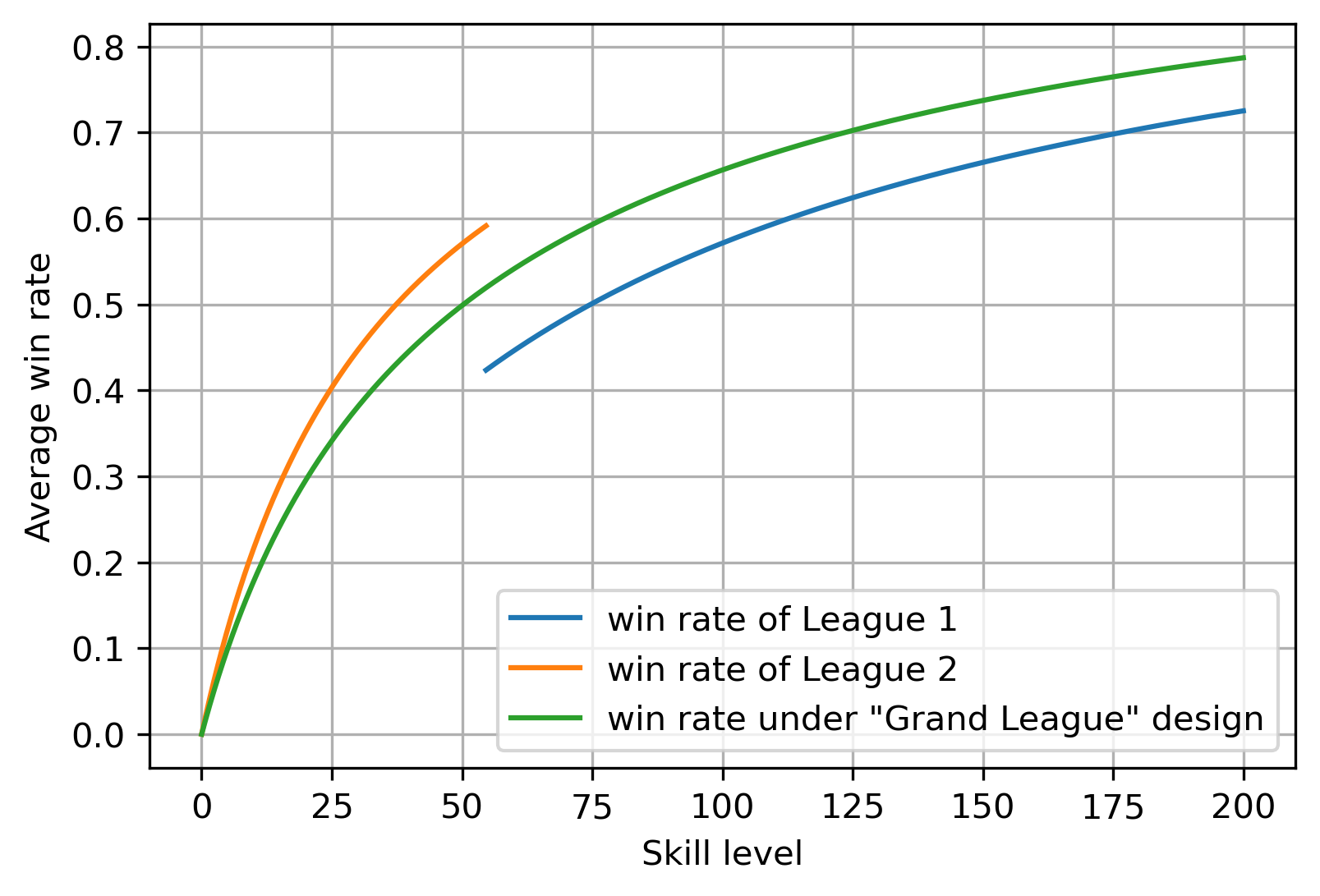}
     \caption{Average win rate before and after grouping}
  \end{subfigure}
  \caption{Results of \cref{exp:lognormal}}\label{fig:lognormal}
\end{figure}
\end{example}

However, the ``novice/veteran" league design is not always optimal. As we showed in \cref{exp:non-convex}, it is possible that the support of optimal solution is not convex in pairs. As shown below, this can happen if $\pi$ is a mixed lognormal distribution with three peaks. 

\begin{example}[$\pi$ is mixed lognormal distribution]\label{exp:mixedlognormal}
   Suppose $\pi$ is the density of a mixed lognormal distribution, i.e., $\pi=0.1\pi_1+0.6\pi_2+0.3\pi_3$, where $\pi_1,\pi_2,\pi_3$ are densities of lognormal distributions with parameters $(4.8, 0.05^2),(5.8,0.05^2),$ and $(7.5,0.05^2)$. In this example, we decompose $\pi$ into two probability densities $\mu_1$ and $\mu_2$ with weights $p_1=0.4$ and $p_2=0.6$. We use the ``novice/veteran" league design as our benchmark. The novice league (resp. veteran league) consists of players with skill levels in $[0,347]$ (resp. $[347,2000]$). This boundary is chosen to guarantee that the weight of the novice league is $0.4$ as required. \cref{fig:mixedlognormal}(a) shows one outcome of \cref{alg:WCCGF}. Blue particles form a ``tail league'' representing density $\mu_1$ while orange particles form a ``bulk league'' representing density $\mu_2$. \cref{fig:mixedlognormal}(b) shows the histogram of $\mu_1$ (in blue) and $\mu_2$ (in orange). Note that the optimal league design consists of a tail and bulk league. Compared with the ``novice/veteran" league design, the Elo loss (i.e., the objective value in \cref{prob:decomposition-2} with Elo loss function in \cref{exp:elo}) is reduced by $7.49\%$. \cref{fig:mixedlognormal}(c) shows the average win rate of each skill level, which are computed in the same way as shown in \cref{exp:uniform}. The green curve represents the win rate of each skill level under the ``novice/veteran" league design. The ``tail/bulk" league design significantly decreases the win rate of novice players (those with skill levels in the range of 0 to 347),  markedly increasing the win rate of veteran players (those with skill levels greater than 347). In practice, this outcome is often considered undesirable as it may lead to frustration among low-skill players because they are paired with high-skill players and are likely to lose. This is an illustration of the failure of the convex in pairs condition.

   \begin{figure}
  \begin{subfigure}[t]{.32\textwidth}
    \centering
    \includegraphics[width=\linewidth]{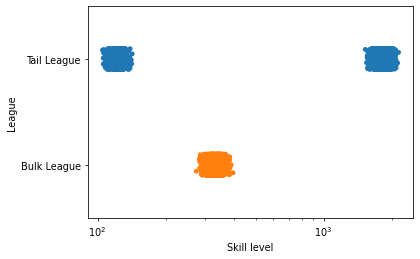}
     \caption{Final positions of particles}
  \end{subfigure}
    \hfill
      \begin{subfigure}[t]{.32\textwidth}
    \centering
    \includegraphics[width=\linewidth]{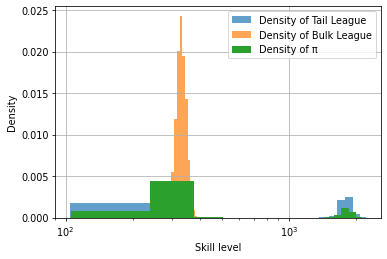}
    \caption{Histogram of $\mu_1,\mu_2$ and $\pi$}
  \end{subfigure}
  \hfill
  \begin{subfigure}[t]{.32\textwidth}
    \centering
    \includegraphics[width=\linewidth]{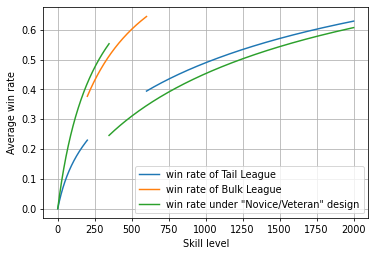}
     \caption{Average win rate under different designs}
  \end{subfigure}
  \caption{Resuts of \cref{exp:mixedlognormal}}\label{fig:mixedlognormal}
\end{figure}
\end{example}

\subsection{Implementation of Wasserstein flow in \cref{def:Wasserstein-CCGF2}}

In this subsection, we implement \cref{def:Wasserstein-CCGF2} to solve \cref{prob:opt-decomposition} in the setting described by \cref{exp:clustering,exp:elo}. In \cref{def:Wasserstein-CCGF2}, we also need to update the weight $p_k$ of each density $\mu_k$. Again, we discreterize the time and flow to get the following implementable algorithm (Algorithm 2). 

\begin{algorithm}
\caption{Wasserstein CCGF particles algorithm with dynamic weights}\label{alg:WCCGF2}
\begin{algorithmic}[1]
\Require Input number of iterations $T$, step size $\eta,\eta_2>0$, number $K$ of sub-populations, and initial distributions $\mu_k(0)$.
\State For each $k\in[K]$, randomly sample a population of $N$ particles $X_k(0)\doteq (\mathbf x_k^i(0))_{i=1}^N$ with $\mathbf x_k^i(0)\in\R^d$ from the initial distributions $\mu_k(0)$. Similarly, define $X_k(t)\doteq (\mathbf x_k^i(t))_{i=1}^N$ for all $t\geq 0$. 
\State Initialize $p_k(0)$ for $k\in[K]$.
\For {each time step $t=1,\ldots, T$}
\For {each sub-population $k=1,\ldots,K$}
    \For {each particle $i=1,\ldots, N$}
        \State Compute the velocity $\phi_k(t,\mathbf x_k^i(t))\in\R^d$ according to equation \cref{eq:WCCGF2} in \cref{def:Wasserstein-CCGF2}.
        \State Update the particle with: 
            $\mathbf x_k^i(t+1)=\mathbf x_k^i(t)+\eta\phi_k(t,\mathbf x_k^i(t)).$
    \EndFor
    \State Compute the velocity $v_k(t, p_k(t))\in\R$ according to equation \cref{eq:WCCGF2}.
    \State Update the weight with $p_k(t+1)=p_k(t)+\eta_2v_k(t,p_k(t))$.
\EndFor
\EndFor
\State Output $X_k(T)$ and $p_k(T)$ for all $k\in[K]$. 
\end{algorithmic}
\end{algorithm}

 In the following example, we show that it is not always optimal to have equal weights. 

    \begin{figure}[h]
      \begin{subfigure}[t]{.32\textwidth}
    \centering
    \includegraphics[width=\linewidth]{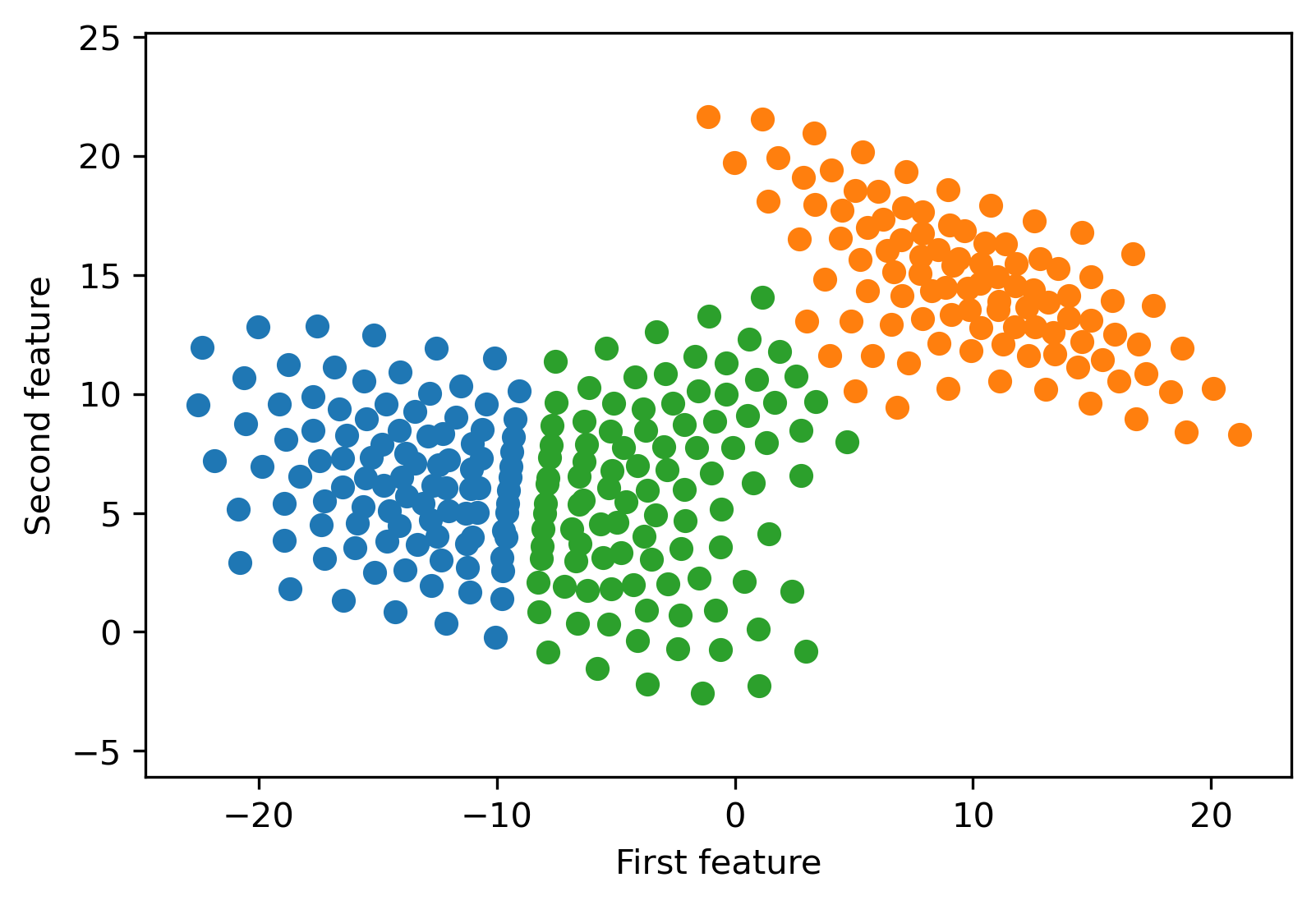}
     \caption{Final positions of particles}
  \end{subfigure}
    \hfill
      \begin{subfigure}[t]{.32\textwidth}
    \centering
    \includegraphics[width=\linewidth]{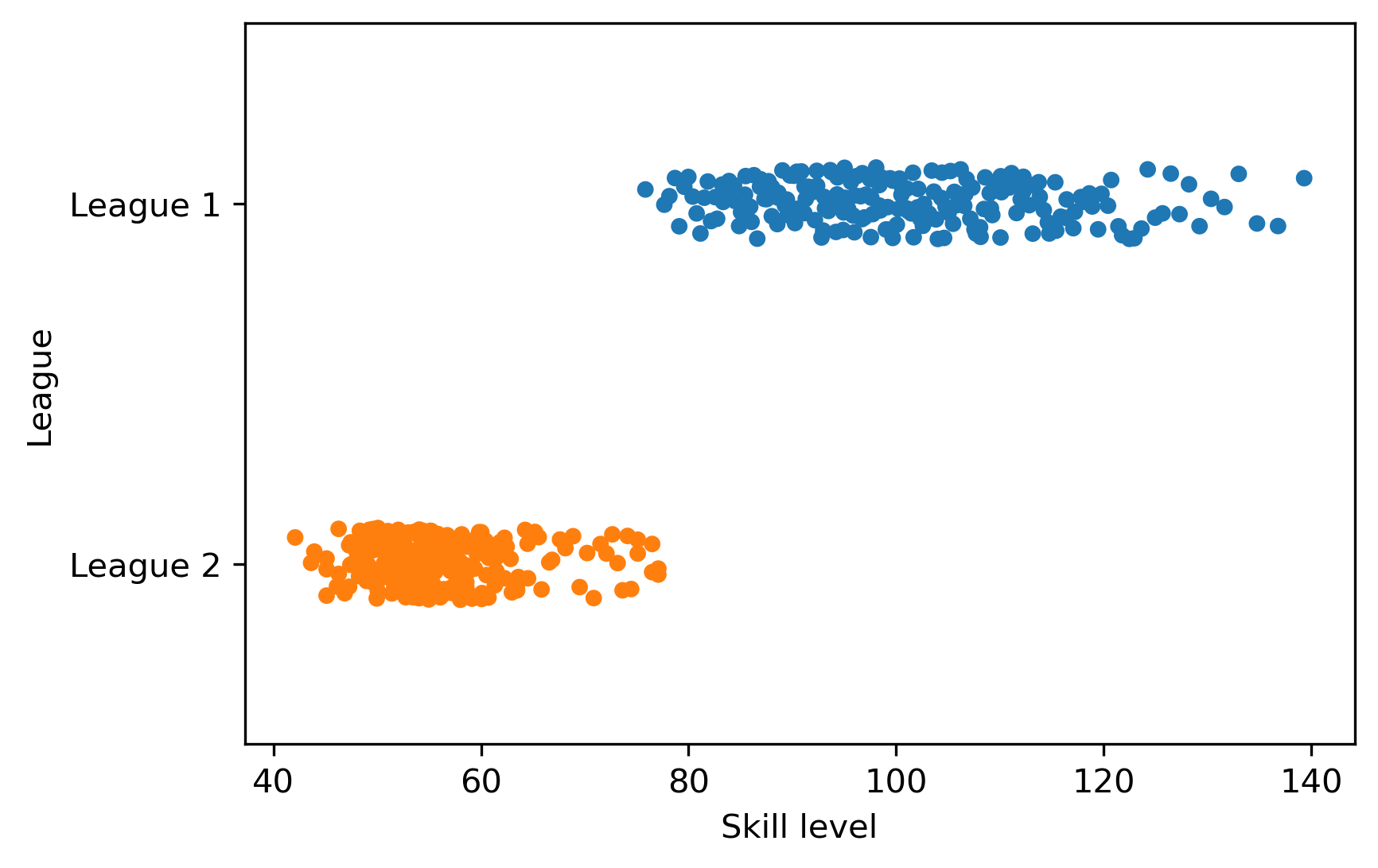}
    \caption{Final positions of particles}
  \end{subfigure}
  \hfill
  \begin{subfigure}[t]{.32\textwidth}
    \centering
    \includegraphics[width=\linewidth]{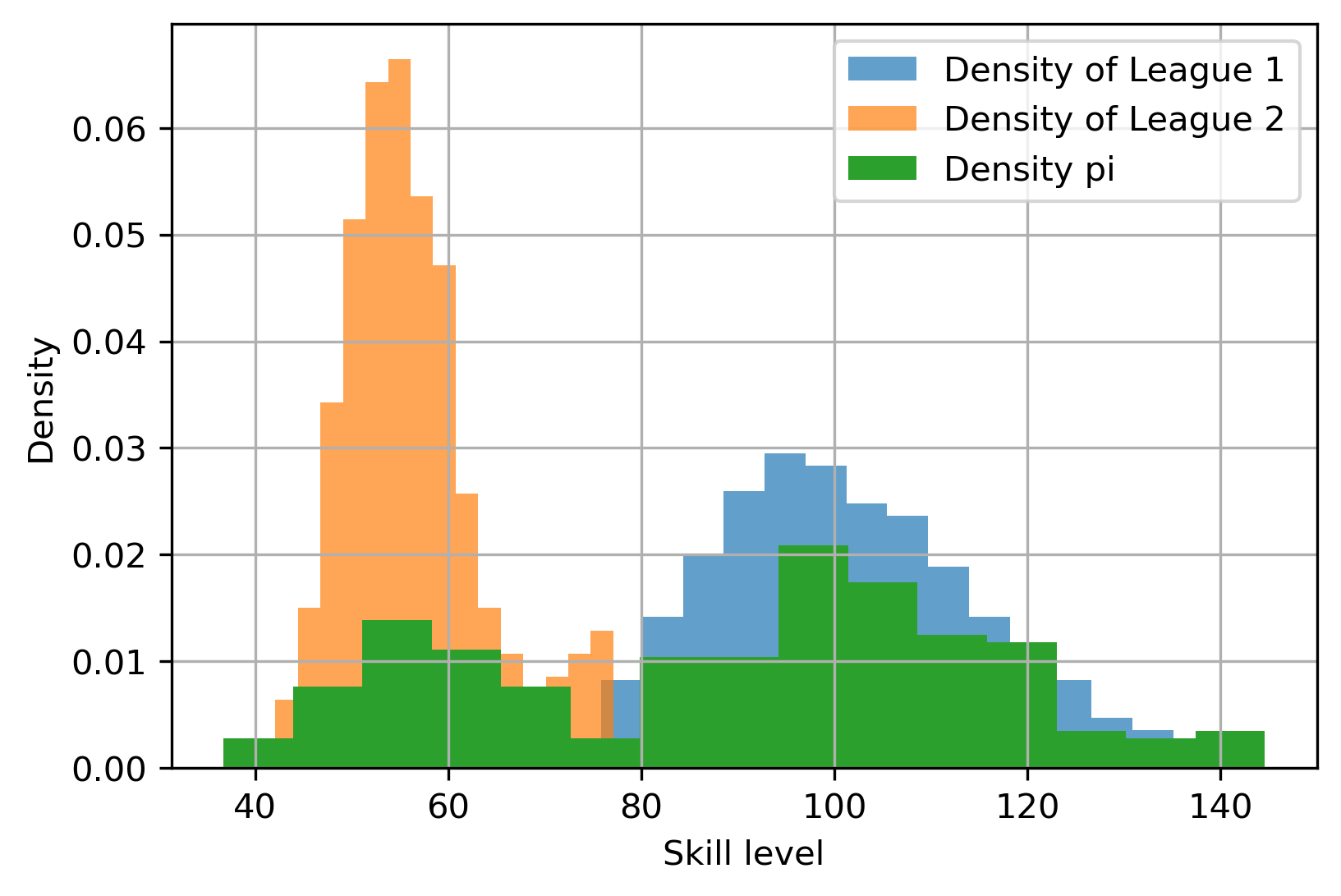}
     \caption{Histogram of $\mu_1,\mu_2$ and $\pi$}
  \end{subfigure}
  \caption{Resuts of \cref{exp:dynamic-weights}}\label{fig:dynamic}
    \end{figure}
    
\begin{example}[WCCGF with dynamic weights]\label{exp:dynamic-weights}
\noindent\textbf{(i)} We revisit the setting described by \cref{exp:multi-objective} with $W$ equal to the identity matrix where the result is shown in \cref{fig:gaussian-2}(c). The output of \cref{alg:WCCGF2} is shown in \cref{fig:dynamic}(a) with weight of blue particle-population being $0.2881$, green particle-population being $0.3087$, and orange particle-population being $0.4322$. In comparison, the weights in \cref{exp:multi-objective} are fixed at $\sfrac{1}{3}$ for all populations. Compared with the result in \cref{exp:multi-objective}, solution of \cref{alg:WCCGF2} reduces the objective value (computed by \cref{eqn:genR}) by $21\%$, reduces the weighted variance loss (computed by \cref{eqn:fixpR}) by $37.15\%$, and increases the total weight loss (computed by $\sum_{k=1}^K\sfrac{1}{p_k}$) by $2.12\%$. Thus, a slight increase in the weight loss leads to a notable reduction in the weighted distribution loss and objective value. We can observe that, compared with \cref{fig:gaussian-2}(c), orange particle population in \cref{fig:dynamic}(a) occupied a larger area, while the blue particle population shrinks. These observations are also reflected in the change of their corresponding weights.

\noindent\textbf{(ii)} Suppose $\pi$ is the density of a mixed lognoraml distribution defined by $\pi=0.3\pi_1+0.7\pi_2$, where $\pi_1$ and $\pi_2$ are densities of lognormal distributions with (scale, shape) parameters given by $(4, 0.1)$ and $(4.6,0.15)$. The density $\pi$ is depicted by green color in \cref{fig:dynamic}(c). The output of \cref{alg:WCCGF2} is shown in \cref{fig:dynamic}(b) with weight of orange particle-population being $0.3303$ and blue particle-population being $0.6697$. \cref{alg:WCCGF2} creates a ``novice/veteran" league design with $33.03\%$ players being categorized as ``novice", and $66.97\%$ players being categorized as ``veteran". \cref{fig:dynamic}(c) shows the histogram of $\mu_1$ and $\mu_2$. Compared with the $50\%$-$50\%$ design, league design of \cref{alg:WCCGF2} reduces the objective value (computed by \cref{eqn:genR}) by $32.85\%$, reduces the weighted Elo loss (computed by \cref{eqn:fixpR}) by $37.25\%$, and increases the total weight loss (computed by $0.0001\sum_{k=1}^K\sfrac{1}{p_k}$) by $13.02\%$. This $0.0001$ guarantees that the weight loss is in the same scale as the distribution loss.

\end{example}

\section{Conclusion}\label{s:conclusion}

In the paper we have provided approaches to solving the optimal probability measure decomposition problem. The approaches use fresh ideas from Wasserstein gradient flow to handle a constrained problem. The numerical implementation uses a particle-based method that provides insights into our motivating examples.

Of course, there are opportunities for future work to extend our approach. On the application side, one could study the league design problem while also considering the spending habits of players as another dimension to designing leagues. For example, players who are big spenders in a game could be given higher priority for winning to keep them highly engaged with the game. This would require the loss function and applying our methodology to this new setting. Many similar extended applications are possible.

On the theory side, one could consider adapting our algorithm to tackle additional constraints for the nature of sub-population decomposition. For example, one could consider requiring sub-populations to have similar moments. This would require developing new optimality conditions and, adapting new optimality conditions and adapting the gradient flow algorithms accordingly.

\newpage

\bibliographystyle{plainnat} 
\bibliography{reference}

\newpage

\begin{appendices}

\pagenumbering{arabic}
\renewcommand*{\thepage}{A.\arabic{page}}
\renewcommand{\thelemma}{A.\arabic{lemma}}
\renewcommand{\thesection}{A.\arabic{section}}
\renewcommand{\theproposition}{A.\arabic{proposition}}
\renewcommand{\thecorollary}{A.\arabic{corollary}}
\renewcommand{\theequation}{A.\arabic{equation}}
\renewcommand{\theremark}{A.\arabic{remark}}
\setcounter{lemma}{0}
\setcounter{section}{0}
\setcounter{proposition}{0}
\setcounter{equation}{0}
\setcounter{corollary}{0}

\begin{center}
    \section*{Appendix: Omitted proofs}
\end{center}

\section{Proofs of \cref{s:structure}}


\subsection{Proof of \cref{prop:disjoint-support}}

\begin{proof}
Consider a $\mathbf x_0$ and a unit norm vector $\mathbf v\in\mathbb{R}^d$ such that $\langle \mathbf v, \nabla^2_{1,2}\ell(\mathbf x_0,\mathbf x_0)\rangle<0$. For the sake of contradiction, suppose for some $\delta,c>0$ and $i,j\in[K]$, $\mathbf x_0\in S_i(\delta,c)\cap S_j(\delta,c)$. Without loss of generality, we assume $x_0=0$, $i=1$, and $j=2$. Let $B_0$ be the ball of radius $r<\delta$ centered at $x_0$.

Define two half-balls
\[
B_1=\{\mathbf x\in B_0: \langle \mathbf v, \mathbf x\rangle>0\},\quad B_2=\{\mathbf x\in B_0:\langle \mathbf v, \mathbf x\rangle<0\}.
\]
Define $q(\mathbf x)=1_{B_1}(\mathbf x)-1_{B_2}(\mathbf x)$.

For $\epsilon>0$ such that $c-\frac{\epsilon}{p_1}>0$ and $c-\frac{\epsilon}{p_2}>0$, construct the following densities
\begin{equation*}
    \begin{split}
\mu^+_1(\mathbf x)&=\mu_1(\mathbf x)-\frac{\epsilon}{p_1} q(\mathbf x),
\quad 
\mu^+_2(\mathbf x)=\mu_2(x)+\frac{\epsilon}{p_2} q(\mathbf x),\quad \mu^+_k(\mathbf x)=\mu_k(\mathbf x),\quad k\geq 3. \\
        \mu^-_1(\mathbf x)&=\mu_1(\mathbf x)+\frac{\epsilon}{p_1} q(\mathbf x),
\quad 
\mu^-_2(\mathbf x)=\mu_2(\mathbf x)-\frac{\epsilon}{p_2} q(\mathbf x),\quad \mu^-_k(\mathbf x)=\mu_k(\mathbf x),\quad k\geq 3. 
    \end{split}
\end{equation*}
Note that both $\{\mu^+_k\}_{k\in[K]}$ and $\{\mu^-_k\}_{k\in[K]}$ are feasible to \cref{prob:opt-decomposition}. To see this, by the choice of $\epsilon$, $\mu_1^+(\mathbf x)>0$ and $\mu_2^-(\mathbf x)>0$ for any $\mathbf x\in\mathbb{R}^d$. Also, by the construction, $\{\mu^+_k\}_{k\in[K]}$ and $\{\mu^-_k\}_{k\in[K]}$ are two decomposition of $\pi$.

We have
\begin{equation*}
    \begin{split}
        p_1(L(\mu_1^+)-L(\mu_1))&=-\epsilon\int\int_{B_0\times B_0} \ell(\mathbf x,\mathbf y)(q(\mathbf x)\mu_1(\mathbf y)+q(\mathbf y)\mu_1(\mathbf x))d\mathbf xd\mathbf y +\frac{\epsilon^2}{p_1}\int\int_{B_0\times B_0} \ell(\mathbf x,\mathbf y)q(\mathbf x)q(\mathbf y)d\mathbf xd\mathbf y\\
        &=-2\epsilon\int\int_{B_0\times B_0} \ell(\mathbf x,\mathbf y)q(\mathbf x)\mu_1(\mathbf y)d\mathbf xd\mathbf y +\frac{\epsilon^2}{p_1}\int\int_{B_0\times B_0} \ell(\mathbf x,\mathbf y)q(\mathbf x)q(\mathbf y)d\mathbf xd\mathbf y,
    \end{split}
\end{equation*}
where the second equality is because $\ell:\R^d\times\R^d\to\R$ is symmetric.
Likewise, 
\begin{equation*}
        p_2(L(\mu_2^+)-L(\mu_2))=2\epsilon\int\int_{B_0\times B_0} \ell(\mathbf x,\mathbf y)q(\mathbf x)\mu_2(\mathbf y)d\mathbf xd\mathbf y +\frac{\epsilon^2}{p_2}\int\int_{B_0\times B_0} \ell(\mathbf x,\mathbf y)q(\mathbf x)q(\mathbf y)d\mathbf xd\mathbf y.
\end{equation*}
So we have 
\begin{align*}
         \sum_{k\in[K]}p_k(L(\mu_k^+)-L(\mu_k))=2\epsilon&\int\int_{B_0\times B_0} \ell(\mathbf x,\mathbf y)(\mu_2(\mathbf y)-\mu_1(\mathbf y))q(\mathbf x)d\mathbf xd\mathbf y \\
         &+ \epsilon^2(\frac{1}{p_1}+\frac{1}{p_2})\int\int_{B_0\times B_0} \ell(\mathbf x,\mathbf y)q(\mathbf x)q(\mathbf y)d\mathbf xd\mathbf y.
\end{align*}

Likewise, we have
\begin{align*}
         \sum_{k\in[K]}p_k(L(\mu_k^-)-L(\mu_k))=-2\epsilon&\int\int_{B_0\times B_0} \ell(\mathbf x,\mathbf y)(\mu_2(\mathbf y)-\mu_1(\mathbf y))q(\mathbf x)d\mathbf xd\mathbf y \\ 
         &+ \epsilon^2(\frac{1}{p_1}+\frac{1}{p_2})\int\int_{B_0\times B_0} \ell(\mathbf x,\mathbf y)q(\mathbf x)q(\mathbf y)d\mathbf xd\mathbf y.
\end{align*}

Suppose $\int\int_{B_0\times B_0} \ell(\mathbf x,\mathbf y)(\mu_2(\mathbf y)-\mu_1(\mathbf y))q(\mathbf x)d\mathbf xd\mathbf y\neq 0$. Then, for $\epsilon$ small enough, either $\sum_{k\in[K]}p_kL(\mu_k^+)$ or $\sum_{k\in[K]}p_kL(\mu_k^-)$ is smaller than $\sum_{k\in[K]}p_kL(\mu_k)$, which contradicts the optimality of $\{\mu_k\}_{k\in[K]}$.

In the case that \[
\int\int_{B_0\times B_0} \ell(\mathbf x,\mathbf y)(\mu_2(\mathbf y)-\mu_1(\mathbf y))q(\mathbf x)d\mathbf xd\mathbf y= 0,
\]
we check the term 
\begin{equation}\label{eq:integral}
    \int_{B_0\times B_0} \ell(\mathbf x,\mathbf y)q(\mathbf x)q(\mathbf y)d\mathbf xd\mathbf y.
\end{equation}
We apply Talyor expansion as follows.
\begin{equation*}
    \begin{split}
        \ell(\mathbf x,\mathbf y)&=\ell(0,0)+(\mathbf x,\mathbf y)^\top \nabla \ell(0,0) +\frac12(\mathbf x,\mathbf y)^\top\nabla^2 \ell(0,0) (\mathbf x,\mathbf y)
+O(\|(\mathbf x,\mathbf y)\|^3).
    \end{split}
\end{equation*}

Therefore, in the integral \cref{eq:integral}, the constant term becomes
 \begin{equation*}
           \int \int_{B_0\times B_0} \ell(0,0)q(\mathbf x)q(\mathbf y)d\mathbf xd\mathbf y = \ell(0,0)\int_{B_0}q(\mathbf x)d\mathbf x\int_{B_0}q(\mathbf y)d\mathbf y=0
    \end{equation*}
    where $\int_{B_0}q(\mathbf x)d\mathbf x=0$.


Note that $f(\mathbf x,\mathbf y)\doteq (\mathbf x,\mathbf y)^\top\nabla \ell(0,0)q(\mathbf x)q(\mathbf y)$ is an odd function on $\R^d\times\R^d$, i.e., $f(\mathbf x,\mathbf y)=-f(-\mathbf x,-\mathbf y)$. In the integral \cref{eq:integral}, the first order term is 
\[
\int\int_{B_0\times B_0}(\mathbf x,\mathbf y)^\top\nabla \ell(0,0)q(\mathbf x)q(\mathbf y)d\mathbf xd\mathbf y=0.
\]


Next, we deal with the second-order term. For each $\mathbf x\in B_0$, define $\mathbf x^-\doteq 2\langle \mathbf x, \mathbf v\rangle \mathbf v-\mathbf x$, which is mirrow image of $\mathbf x$ with respect to the hyperplane $\{\mathbf z\in\R^d:\mathbf z^\top \mathbf v=0\}$. Note that the distribution of $\mathbf x^-$ and $\mathbf x$ are identical under $q(\mathbf x)$. Denote
\[
\nabla^2 \ell(0,0)=\begin{bmatrix} H_1 &H_2\\
H_2^T &H_3\end{bmatrix}.
\]
In the integral \cref{eq:integral}, the second order term is
\begin{equation*}
        \begin{split}
            &\int\int_{B_0\times B_0} (\mathbf x,\mathbf y)^\top\nabla^2 \ell(0,0)
            (\mathbf x,\mathbf y)q(\mathbf x)q(\mathbf y)d\mathbf xd\mathbf y\\
            &= \int\int_{B_0\times B_0} x^\top H_1 x q(\mathbf x) q(\mathbf y)d\mathbf xd\mathbf y+
            2\int\int_{B_0\times B_0} \mathbf x^\top H_2 \mathbf y q(\mathbf x)q(\mathbf y)d\mathbf xd\mathbf y+\int\int_{B_0\times B_0} \mathbf y^\top H_3 \mathbf y q(\mathbf y)q(\mathbf x)d\mathbf xd\mathbf y\\
            &= 2\int\int_{B_0\times B_0} \mathbf x^\top H_2 \mathbf y q(\mathbf x)q(\mathbf y)d\mathbf xd\mathbf y\\
            &=\int\int_{B_0\times B_0} (\mathbf x+\mathbf x^-)^\top H_2 y q(\mathbf x)q(\mathbf y)d\mathbf xd\mathbf y\\
            &=\frac12 \int\int_{B_0\times B_0} (\mathbf x+\mathbf x^-)^\top H_2 (y+y^{-}) q(\mathbf x)q(\mathbf y)d\mathbf xd\mathbf y\\
            &=2\int\int_{B_0\times B_0} \mathbf v^\top \mathbf x   (\mathbf v^\top H_2 \mathbf v) \mathbf v^\top \mathbf y q(\mathbf x)q(\mathbf y)d\mathbf xd\mathbf y\\
            &=2(\mathbf v^\top H_2 \mathbf v) \int\int_{B_0\times B_0} (\mathbf v^\top \mathbf x)(\mathbf v^\top \mathbf y) (1_{\mathbf v^\top \mathbf x>0}-1_{\mathbf v^\top \mathbf x<0})(1_{\mathbf v^\top \mathbf y>0}-1_{\mathbf v^\top\mathbf  y<0})dxdy\\
        &=2(\mathbf v^\top H_2 \mathbf v)(\int_{B_0}(\mathbf v^\top \mathbf x)(1_{\mathbf v^\top \mathbf x>0}-1_{\mathbf v^\top \mathbf x<0})d\mathbf x)^2  =-O((rm(B_0))^2),
        \end{split}
    \end{equation*}
where the second equality is because 
\[
\int\int_{B_0\times B_0} \mathbf x^\top H_1 \mathbf x q(\mathbf x) q(\mathbf y)d\mathbf xd\mathbf y=(\int_{B_0}\mathbf x^\top H_1\mathbf xq(\mathbf x)d\mathbf x)(\int_{B_0}q(\mathbf y)d\mathbf y)=0,
\]
since $\int_{B_0}q(\mathbf y)d\mathbf y=0$; the third equality is because 
\[
\int\int_{B_0\times B_0} \mathbf x^\top H_2 \mathbf y q(\mathbf x)q(\mathbf y)d\mathbf xd\mathbf y=
\int\int_{B_0\times B_0} (\mathbf x^-)^\top H_2 \mathbf y q(\mathbf x)q(\mathbf y)d\mathbf xd\mathbf y;
\]
and the last equality is because $v^\top H_2v<0$ by our condition.

Plug in this estimate and symmetry of $B_1,B_2$,  we find 
\begin{equation*}
        \int\int_{B_0\times B_0} \ell(\mathbf x,\mathbf y)q(\mathbf x)q(\mathbf y)d\mathbf xd\mathbf y=-|O(r^2)m(B_0)^2|+O(r^3)m(B_0)^2.
\end{equation*}


Thus, for $r$ small enough, 
    \begin{equation*}
        \sum_{k\in[K]}p_k(L(\mu_k^-)-L(\mu_k))=\int\int_{B_0\times B_0} \ell(\mathbf x,\mathbf y)q(\mathbf x)q(\mathbf y)d\mathbf xd\mathbf y<0,
    \end{equation*}
    which again contradicts the optimality of $\{\mu_i\}_{i \in [K]}$. 
\end{proof}

\subsection{Proof of \cref{cor:disjoint-support}}

\begin{proof}
\textbf{Variance loss:}
Let $\ell(\mathbf x,\mathbf y)=\langle \mathbf x-\mathbf y,W(\mathbf x-\mathbf y)\rangle$. The gradient and Hessian matrix are given by  
\[
\nabla \ell(\mathbf x,\mathbf y)=2W[\mathbf x-\mathbf y,\mathbf y-\mathbf x],\quad \nabla^2 \ell(\mathbf x,\mathbf y)=2\begin{bmatrix}
W &-W\\
-W &W
\end{bmatrix}.
\]
Thus, for any $\mathbf x_0$ and unit norm vector $\mathbf v$, 
$\langle \mathbf v, \nabla^2_{1,2}\ell(\mathbf x_0,\mathbf x_0)\mathbf v\rangle=-2W$. By \cref{prop:disjoint-support}, as long as not all entries in $W$ are non-positive (equivalently, $W$ is not negative semi-definite), the optimal densities have no overlap for their $(\delta,c)$-support. In particular, if $W=\text{id}$, i.e., $\ell(\mathbf x,\mathbf y)=\|\mathbf x-\mathbf y\|^2$, then the optimal densities have no overlap for their $(\delta,c)$-interior densities. 

\textbf{Elo loss:}
Let $\ell(x,y)=\frac{x^2+y^2}{(x+y)^2}$. We note that 
\[
\partial_1 \ell=\frac{2xy-2y^2}{(x+y)^3},\quad \partial_2 \ell=\frac{2xy-2x^2}{(x+y)^3},
\]
\[
\nabla^2 \ell=
\frac{1}{(x+y)^4}\begin{bmatrix}8y^2-4xy & 2(x^2+y^2-4xy)\\
2(x^2+y^2-4xy)&  8x^2-4xy
\end{bmatrix}.
\]
In this case, $d=1$ and the only unit norm vector is $v=1$, thus
    \[
    \nabla^2\ell(x,x)=\frac{1}{4x^2}\begin{bmatrix}
        1 & -1\\
        -1 & 1
    \end{bmatrix}.
    \]
Hence, $\langle v, \nabla^2_{1,2}\ell(x,x)v\rangle<0$. By \cref{prop:disjoint-support}, the optimal densities have no overlap for their $(\delta,c)$- support.
\end{proof}

\subsection{Proof of claim in \cref{exp:non-convex}}

\begin{proof}
  Suppose we decompose $\pi$ into two disbritutions $\mu_1$ and $\mu_2$, i.e., $K=2$ in \cref{prob:decomposition-2}. Consider the following decomposition:
\[
\mu_1=\frac{1}{2}\delta_{a-1}+\frac12\delta_a,\quad 
\mu_2=\frac{1}{2}\delta_{a+1}+\frac12\delta_a.
\]

Omitting the constant $-\frac{1}{4}$, the Elo loss of this decomposition is given by 

\begin{equation*}
    \begin{split}
        &\frac{1}{2}(\frac{1}{2}\cdot P(\text{equal match})+
\frac{2a^2-2a+1}{(2a-1)^2}\cdot P(\text{unequal match}))+
\frac{1}{2}(\frac{1}{2}\cdot P(\text{equal match})\\&+
\frac{2a^2+2a+1}{(2a+1)^2}\cdot P(\text{unequal match}))\\
&=\frac{1}{4}(\frac{1}{2}+\frac{1}{2}+
\frac{1}{2(2a-1)^2})+
\frac{1}{4}(\frac{1}{2}+\frac{1}{2}+
\frac{1}{2(2a+1)^2})\\
&=\frac12(1+\frac{1}{4(2a-1)^2}+\frac{1}{4(2a+1)^2}),
    \end{split}
\end{equation*}
where the second equality is due to $P(\text{equal match})=P(\text{unequal match})=\frac{1}{2}$ and $\frac{2a^2-2a+1}{(2a-1)^2}=\frac{1}{2}+\frac{1}{2(2a-1)^2} $.

Meanwhile, consider the following bulk/tail decomposition:
\[
\mu_1=\frac{1}{2}\delta_{a-1}+\frac12\delta_{a+1},\quad 
\mu_2=\delta_{a}.
\]

Omitting the constant $-\frac{1}{4}$, its Elo loss is
\begin{equation*}
    \begin{split}
        &\frac{1}{2}(\frac{1}{2}\cdot P(\text{equal match})+
\frac{2a^2+2}{(2a)^2}\cdot P(\text{unequal match}))+
\frac{1}{4}\\=&\frac{1}{4}(\frac 12+\frac 12+\frac{2}{4a^2})+\frac 14\\
&=\frac 12(1+\frac{1}{4a^2}).
    \end{split}
\end{equation*}

Then, when $a>1$ is close to $1$,
\[
\frac{1}{4(2a-1)^2}+\frac{1}{4(2a+1)^2}\geq \frac{1}{4a^2}. 
\]

So, the second plan, which is nonconvex, yields a lower Elo loss. 
\end{proof}

\subsection{Proof of \cref{lem:variance-convex}}

\begin{proof}
Let $(\mu_1,\ldots,\mu_K)$ be the optimal solution to \cref{prob:decomposition-2} with variance loss. For the sake of contradiction, suppose, for some $\delta, c>0$, there exists $\mathbf x_0,\mathbf x_1\in S_1(\delta,c)$, $\lambda\in(0,1)$ such that $\mathbf z\doteq\lambda \mathbf x_0+(1-\lambda)\mathbf x_1\in S_2(\delta,c)$. Without loss of generality, assume $\mathbf x_0=0$. Let $B_0, B_1, B_{\mathbf z}$ be open ball of radius $\delta$ centered at $\mathbf x_0,\mathbf x_1$, and $\mathbf z$, respectively. By the definition of $S_1(\delta,c)$ and $S_2(\delta,c)$, we know that $\mu_1(\mathbf x)>c$ $\forall x\in B_0\cup B_1$ and $\mu_2(\mathbf x)>c$ for all $\mathbf x\in B_{\mathbf z}$.

Consider the following alternative solution $(\mu'_1,\ldots,\mu'_K)$ obtained by 
\[
\mu'_1(\mathbf x)=\mu_1(\mathbf x)-\frac{c}{p_1} (1_{B_0}(\mathbf x)-1_{B_z}(\mathbf x)),
\quad 
\mu'_2(\mathbf x)=\mu_2(\mathbf x)-\frac{c}{p_2}(1_{B_{\mathbf z}}(\mathbf x)-1_{B_0}(\mathbf x)),
\]
and $\mu'_k=\mu_k$ for all $k>2$.

Let $\mathbf m_i$ (resp. $\mathbf m_i'$) be the mean of probability measure $\mu_i$ (resp. $\mu_i'$) consisting of means in $d$-dimensions of the underlying features. Since $\mathbf m_1'\neq \mathbf m_1$ is the mean of $\mu_1'$, we have
\[
L(\mu_1')=\int_{\R^d} (\mathbf x-m_1')^\top W(\mathbf x-\mathbf m_1')d\mu_1'(\mathbf x)<\int_{\R^d} (\mathbf x-\mathbf m_1)^\top W(\mathbf x-\mathbf m_1)d\mu_1'(\mathbf x).
\]
Therefore,
\begin{equation*}
    \begin{split}
         p_1(L(\mu_1)-L(\mu_1'))&> c\int_{B_0} (\mathbf x-\mathbf m_1)^\top W(\mathbf x-\mathbf m_1)dx-c\int_{B_z} (\mathbf x-\mathbf m_1)^\top W(\mathbf x-\mathbf m_1)dx\\
&= c\int_{B_0} ((\mathbf x-\mathbf m_1)^\top W(\mathbf x-\mathbf m_1)-(\mathbf x+\mathbf z-\mathbf m_1)^\top W(\mathbf x+\mathbf z-\mathbf m_1))dx\\
&= c\int_{B_0} (-2(\mathbf x-\mathbf m_1)^\top W\mathbf z-\mathbf z^\top W\mathbf z)d\mathbf x\\
&=cm(B_0)(2\mathbf m_1-\mathbf z)^\top W\mathbf z
    \end{split}
\end{equation*}
Similarly,
\[
p_2(L(\mu_2)-L(\mu_2'))> -c\int_{B_0} (\mathbf x-\mathbf m_2)^\top W(\mathbf x-m_2)d\mathbf x+c\int_{B_{\mathbf z}} (\mathbf x-\mathbf m_2)^\top W(\mathbf x-\mathbf m_2)dx=-c m(B_0)(2\mathbf m_2-\mathbf z)^\top W \mathbf z.
\]

The change of objective value is given by
\[
\Delta_1\doteq\sum_{k=1}^Kp_k(L(\mu_k)-L(\mu'_k))=p_1(L(\mu_1)-L(\mu_1'))+p_2(L(\mu_2)-L(\mu_2'))
> 2cm(B_0) (\mathbf m_1-\mathbf m_2)^\top W\mathbf z. 
\]

Next, we consider another swapping 
\[
\mu''_1(\mathbf x)=\mu_1(\mathbf x)-\frac{c}{p_1} (1_{B_1}(\mathbf x)-1_{B_z}(\mathbf x)),
\quad 
\mu''_2(\mathbf x)=\mu_2(\mathbf x)-\frac{c}{p_2} (1_{B_z}(\mathbf x)-1_{B_1}(\mathbf x)),
\]
and $\mu''_k=\mu_k$ for all $k>2$. We note that 
\begin{align*}
p_1(L(\mu_1)-L(\mu_1'')) &> c\int_{B_1} (\mathbf x-m_1)^\top W(\mathbf x-\mathbf m_1)d\mathbf x-c\int_{B_z}  (\mathbf x-\mathbf m_1)^\top W(\mathbf \mathbf x-\mathbf m_1)d\mathbf x \\
&=cm(B_0)(2\mathbf m_1-\mathbf z+\mathbf x_1)^\top W(\mathbf z-\mathbf x_1),
\end{align*}
and 
\begin{align*}
p_2(L(\mu_2)-L(\mu_2'')) &> -c\int_{B_1} (\mathbf x-\mathbf m_2)^\top W(\mathbf x-\mathbf m_2)d\mathbf x+c\int_{B_z} (\mathbf x-\mathbf m_2)^\top W(\mathbf x-\mathbf m_2)d\mathbf x \\
&=cm(B_0)(-2\mathbf m_2+\mathbf z-\mathbf x_1)^\top W(\mathbf z-\mathbf x_1).
\end{align*}
The change of objective value is given by
\[
\Delta_2\doteq\sum_{k=1}^Kp_k(L(\mu_k)-L(\mu''_k))=p_1(L(\mu_1)-L(\mu_1''))+p_2(L(\mu_2)-L(\mu_2''))> 2cm(B_0)(\mathbf z-\mathbf x_1)^\top W(\mathbf m_1-\mathbf m_2). 
\]
Note that $\lambda \Delta_1+(1-\lambda)\Delta_2>2cm(B_0)(\mathbf m_1-\mathbf m_2)^\top W(\mathbf z-(1-\lambda)\mathbf x_1)=0$ since $\mathbf z-(1-\lambda)\mathbf x_1=\lambda \mathbf x_0=0$, which indicates that either $\Delta_1>0$ or $\Delta_2>0$. In either case $(\mu_k)_{k\in[K]}$ is not optimal. 
\end{proof}

\section{Proofs for \cref{s:preliminary}}

\subsection{Proof of \cref{lem:U}}

\begin{proof}
Consider a curve $\mu(\cdot):[0,1]\rightarrow\calP_2(\R^d)$ with $\mu(0)=\mu$ and velocity $\phi:\R^d\to\R^d$ at time $t=0$; that is, $\frac{d}{dt}\mu(t,\mathbf x)|_{t=0}=-\nabla\cdot (\phi(\mathbf x)\mu(0,\mathbf x))$.
Then
\begin{equation*}
    \begin{split}
\frac{d}{dt} L(\mu(t))|_{t=0}&=\tfrac{d}{dt} \int_{\R^d}\int_{\R^d} \ell(\mathbf x,\mathbf y)d\mu(t,\mathbf x)d\mu(t,\mathbf y)|_{t=0}\\
&=\int_{\R^d}\frac{d}{dt}\int_{\R^d} \ell(\mathbf x,\mathbf y) d\mu(t,\mathbf x))d\mu(t, y)+
\int_{\R^d}\frac{d}{dt}\int_{\R^d} \ell(\mathbf x,\mathbf y) d\mu(t,\mathbf y))d\mu(t,\mathbf x)|_{t=0}\\
&=-\int_{\R^d}\int_{\R^d} \ell(\mathbf x,\mathbf y)\nabla\cdot(\phi(\mathbf x)\mu(\mathbf x))d(\mathbf x)d\mu(\mathbf y)-\int_{\R^d}\int_{\R^d} \ell(\mathbf x,\mathbf y) \nabla\cdot(\phi(\mathbf y)\mu(\mathbf y))d\mathbf yd\mu(\mathbf x)\\
&=\int_{\R^d}\int_{\R^d} \phi(\mathbf x)^\top\nabla_1 \ell(\mathbf x,\mathbf y)d\mu(\mathbf x)d\mu(\mathbf y)+\int_{\R^d}\int_{\R^d} \phi(\mathbf y)^\top\nabla_2 \ell(\mathbf x,\mathbf y) d\mu(\mathbf x)d\mu(\mathbf y)\\
&=\int_{\R^d} d\mathbf x\phi(\mathbf x)^\top\left(\int_{\R^d} (\nabla_1 \ell(\mathbf x,\mathbf z)+\nabla_2 \ell(\mathbf z,\mathbf x))d\mu(\mathbf z)\right).   \qquad \qquad \qedhere
\end{split} 
\end{equation*} 
\end{proof}

\subsection{Proof of \cref{lem:KLder}}
\begin{proof}
    Let $\mubar=\sum_{k=1}^Kp_k\mu_k$. For each $k\in[K]$, consider a curve $\mu_k(\cdot):[0,1]\rightarrow\calP_2(\R^d)$ with $\mu_k(0)=\mu_k$ and velocity $\phi_k:\R^d\to\R^d$ at time $t=0$; that is, $\frac{d}{dt}\mu(t,\mathbf x)_k|_{t=0}=-\nabla\cdot (\phi_k(\mathbf x)\mu_k(0,\mathbf x))$. We have
\begin{equation*}
    \begin{split}
        \tfrac{d}{dt}\text{KL}(\mubar(t)\|\pi)|_{t=0}&= \tfrac{d}{dt}\int_{\R^d}\mubar(t)\log\tfrac{\mubar(t)}{\pi}d\mathbf x|_{t=0}\\
&=\int(\tfrac{d}{dt} \mubar(t))\log\tfrac{\mubar(t)}{\pi}d\mathbf x|_{t=0}\\
&=\sum_{k\in [K]}p_k\int_{\R^d}(\tfrac{d}{dt} \mu_k(t))\log\tfrac{\mubar(t)}{\pi}d\mathbf x|_{t=0}\\
&=-\sum_{k\in[K]}p_k\int_{\R^d} \nabla\cdot(\phi_k\mu_k)\log\tfrac{\mubar}{\pi}d\mathbf x\\
&=\sum_{k\in [K]}p_k\int_{\R^d} \langle\phi_k(\mathbf x), s_{\mubar}(x)-s_{\pi}(\mathbf x)\rangle \mu_k (\mathbf x)d\mathbf x\\
&=\sum_{k\in [K]}\int_{\R^d} \langle\phi_k(\mathbf x), p_k(s_{\mubar}(x)-s_{\pi}(x))\rangle \mu_k(\mathbf x)d\mathbf x\\
&=\sum_{k\in [K]}\langle\phi_k(\mathbf x), p_k(s_{\mubar}(\mathbf x)-s_{\pi}(\mathbf x))\rangle_{\mu_k},
    \end{split}
\end{equation*}
which satisfies the definition of the Wasserstein gradient in equation \cref{eq:Wasserstein-gradient-K}.
\end{proof}

\section{Proofs for \cref{s:optimality-condition}}

\subsection{Proof of \cref{prop:opt-condition-1}}

We start with several useful lemmas. 

\begin{lemma}
\label{lem:TP}
For any $\bfmu\in\calP_{\pi,\bfp}$, 
\[
\mathcal{T} \calP_{\pi,p}(\bfmu)\subseteq \mathcal{T}' \calP_{\pi,p}(\bfmu)\doteq\left\{\bfphi:
\sum_{k\in [K]} p_k\nabla \cdot(\mu_k \phi_k)=0\right\}. 
\]    
\end{lemma}
\begin{proof}
For any $\bfphi\in\calT\calP_{\pi,p}(\bfmu)$, there exists a curve $\bfmu(t)$ in $\calP_{\pi,p}$ such that $\frac{d}{dt}\mu_k(t)\large|_{t=0}=-\nabla\cdot(\mu_k(0) \phi_k)$ and $\mu(0)=\bfmu$. Since $\bfmu(t)\in\calP_{\pi,p}$, $\sum_{k\in[K]}p_k\mu_k(t)=\pi$ $\forall t\in[0,1]$. This equality holds in the sense that for all $t\in[0,1]$,
    \begin{equation*}
        \sum_{k\in[K]}p_k\int_{\R^d} f(\mathbf x)d\mu_k(t,\mathbf x)=\int_{\R^d} f(\mathbf x)d\pi(\mathbf x),\quad\forall f\in C_c^\infty(\R^d).
    \end{equation*}
Note that the right-hand side does not depend on $t$. Hence, $\sum_{k\in[K]}p_k\frac{d}{dt}\int_{\R^d} f(\mathbf x)d\mu_k(t,\mathbf x)=0$, $\forall f\in C^\infty_c(\R^d)$ for all $t\in[0,1]$, which means $\frac{d}{dt}\sum_{k\in[K]}p_k\mu_k(t)=0$ for all $t\in[0,1]$. Therefore, we have 
\[
\sum_{k\in [K]} p_k\nabla \cdot(\mu_k(0) \phi_k)=-\frac{d}{dt}\sum_{k\in [K]} p_k \mu_k(t)|_{t=0}=0,
\]
which means $\bfphi\in\calT'\calP_{\pi,p}(\bfmu)$, thus completing the proof. 
\end{proof}

\begin{lemma}
\label{lem:diffeo}
Given a $\phi=\nabla \Psi$ with $\Psi\in C^\infty_c(\R^d)$, denote $T_{t}(\mathbf x)=\mathbf x+\phi(\mathbf x)t$. Then for $t < t_0\doteq \sfrac{1}{\|\nabla \phi\|_\infty}$,  $T_t$ is a diffeomorphism of $\reals^d$.  Moreover, $\|T_t^{-1}(\mathbf x)-\mathbf x\|\leq \|\phi\|_\infty t$ (for all $t$, not just $t < t_0$).
\end{lemma}
\begin{proof}
First, we show $T_t$ is one-to-one. If 
$T_t(\mathbf x)=T_t(\mathbf y)$ for some given $\mathbf x$ and $\mathbf y$, then we have 
\[
\|\mathbf x-\mathbf y\|=t\|\phi(\mathbf x)-\phi(\mathbf y)\|\leq t\|\nabla \phi\|_\infty\|\mathbf x-\mathbf y\|,
\]
by the mean value theorem and the Cauchy-Schwartz inequality. Thus, if $t<1/\|\nabla \phi\|_\infty$, we have $\mathbf x=\mathbf y$.

Next, we show that $T_t$ is onto. It suffices to show for any $\mathbf x$, there exists a $\mathbf z$ such that $\mathbf x = \mathbf z+\phi(\mathbf z)t$. Since $\Phi=\nabla\Psi$ and $\psi$ is compactly supported, $T_t=\text{id}$ outside the compact support. Within the compact support, define $f_{\mathbf x,\mathbf t}(\mathbf z)=\mathbf x-\phi(\mathbf z)t$. We show that $f_{\mathbf x,t}$ is a contraction mapping. To see this, for any $\mathbf z_1\neq\mathbf z_2$ in the compact support, $\|f_{\mathbf x,t}(\mathbf z_1)-f_{\mathbf x,t}(\mathbf z_2)\|=\|\phi(\mathbf z_1)-\phi(\mathbf z_2)\|t$. By the mean value theorem, we have $\|\phi(\mathbf z_1)-\phi(\mathbf z_2)\|\leq \|\nabla \phi\|_{\infty}\|\mathbf z_1-\mathbf z_2\|$. This implies $\|f_{\mathbf x,t}(\mathbf z_1)-f_{\mathbf x,t}(\mathbf z_2)\|\leq \|\nabla \phi\|_{\infty}t\|\mathbf z_1-\mathbf z_2\|< \|\mathbf z_1-\mathbf z_2\|$ for $t< t_0$. Hence, $f_{\mathbf x,t}$ is a contraction mapping. By the Banach fixed point theorem, $f_{\mathbf x,t}$ admits a fixed point, i.e., there exists a $\mathbf z$ such that $\mathbf x-\phi(\mathbf z)t=\mathbf z$. Hence, $T_t$ is also onto. 

Then, we show that $T_t$ is differentiable with a differentiable inverse. The differential of $T_t$ is given by 
\[
\nabla T_t(\mathbf x)=\text{id}+t \nabla \phi(\mathbf x),
\]
where $\text{id}$ is the identity matrix. Then, using the inverse map rule, the differential of the inverse $T^{-1}_t$ is given by
\[
\nabla T_t^{-1}(\mathbf x)=(\text{id}+t \nabla \phi(T_t^{-1}(\mathbf x)))^{-1}.
\]
Note that $\text{id}$ and $\nabla \phi(T_t^{-1}(x))$ are both invertible. Recall that the set of invertible matrices is a convex cone. Hence, since $t>0$, $\text{id}+t \nabla \phi(T_t^{-1}(\mathbf x))$ is also invertible. This means that $\nabla T_t^{-1}(\mathbf x)$ exists. 

Thus, $T_t$ is a diffeomorphism. Also, note that if we let $\mathbf z=T^{-1}_t(\mathbf x)$, i.e.,  $\mathbf x=\mathbf z+\phi(\mathbf z)t$, then $\|\mathbf x-\mathbf z\|\leq t\|\phi\|_\infty$. 
\end{proof}

\begin{proof}[Proof of \cref{prop:opt-condition-1}]
Recall that, by definition, $\nabla F_\bfp$ is the unique element in $\calT\calP^{\otimes K}_2(\bfmu^*)=\bigotimes_{k\in[K]}\overline{\{\nabla \Psi: \Psi\in C^\infty_c(\reals^d)\}}^{\calL^2_{\mu_k}}$ such that for all $\bfmu(\cdot):[0,\infty)\to\calP_2^{\otimes K}(\R^d)$ with $\bfmu(0)=\bfmu^*$ and velocity $ \bfphi\in \calT\calP^{\otimes K}_2(\bfmu^*)$, 
\begin{align*}
\frac{d}{dt}F_{\bfp}(\bfmu(t))|_{t=0} 
&=\langle\nabla F_\bfp,\phi\rangle_{\bfmu^*} \\
&=\sum_{k\in[K]}\langle\nabla_{\mu_k} F_\bfp,\phi_k\rangle_{\mu^*_k} \\
&=\sum_{k\in[K]}\int\phi_k(x)^\top\nabla_{\mu_k}F_\bfp(x)d\mu^*_k(x).
\end{align*}
Since $F_\bfp(\bfmu)=\sum_{k\in[K]}p_kL(\mu_k)$, $\nabla_{\mu_k} F_\bfp(\bfmu)=p_k\nabla L(\mu_k)$.

For the sake of contradiction, suppose $\|\nabla F_\bfp(\bfmu^*)\|_{\mathcal{T}\calP_{\pi,\bfp}(\bfmu^*)}>0$. Then, there exists a $\bfphi \in \calT\calP^{\otimes K}_2(\bfmu^*)$ such that
\begin{enumerate}[label=(\roman*)]
    \item $\bfphi=\nabla\Phi$ for some $\Phi\in C_c^\infty(\R^d)$;
    \item $ A\doteq \langle \nabla F_\bfp,\bfphi\rangle_{\bfmu^*}= \sum_{k\in [K]} p_k\langle \phi_k, \nabla L\rangle_{\mu^*_k}<0$;
    \item  $\sum_{k\in [K]} p_k\nabla \cdot(\mu^*_k \phi_k)=0$,
\end{enumerate}
where the equality in claim (iii) holds in the sense that, for any $f\in C^\infty_c(\R^d)$, 
    \[
    \sum_{k\in[K]}p_k\int_{\R^d}\langle\nabla f(x),\phi_k(x)\rangle d\mu^*_k(x)=0.
    \]
Claim (i) and (ii) hold because $\nabla F_\bfp$ is in the closure of $\{\nabla \Psi: \Psi\in C^\infty_c(\reals^d)\}$.  Claim (iii) holds by \cref{lem:TP} since $\bfphi\in \calT\calP^{\otimes K}_2(\bfmu^*)$. 

When such a $\bfphi$ exists, we show that we can construct a solution to \cref{prob:decomposition-2} which is strictly better than $\bfmu^*$. Consider the following time-independent velocity field
\[
\phi_k(t)=\phi_k,\quad\forall  t\in[0,1].
\]
Then, the map $T^\phi_{k,t}(\mathbf x)=\mathbf x+\phi_k(\mathbf x)t$ defines a curve $\mu_k(t)\doteq T^\phi_{k,t}\sharp\mu^*_k$ for $t\in[0,1]$ in $\calP_2(\R^d)$. Define $\bar{\mu}(t)=\sum_{k\in[K]}p_k\mu_k(t)$ in the sense that, for all $t\in[0,1]$, $\bar{\mu}(t)$ is a probability measure in $\calP_2(\R^d)$ such that for all $ f\in C^\infty_c(\R^d)$, $\int_{\R^d} fd\bar{\mu}(t)=\sum_{k\in[K]}p_k\int_{\R^d} fd\mu_k(t)$. In particular, $\bar{\mu}(0)=\pi$ since $\sum_{k\in[K]}p_k\mu_k(0)=\pi$. 

Note that the curve $\bfmu\doteq(\mu_1(t),\ldots,\mu_K(t))$ may not be in the feasible set $\calP_{\pi,\bfp}$, i.e., $\Bar{\mu}(t)$ may not equal to $\pi$. Hence, to get the ``better solution" we are looking for, we need to ``project" this $\bfmu(t)$ back to the feasible set $\calP_{\pi,\bfp}$. However, this ``projection'' is subtle since the underlying space is not a Hilbert space.

Below, we show the following facts: there exists a $C>0$ such that for all $ t\in[0,1]$,
\begin{enumerate}[label=(\roman*)]
    \item $F_\bfp(\bfmu(t))\leq F_\bfp(\bfmu^*)+At+\frac12 Ct^2$;
    \item there exists an optimal transport map $S_t$ from $\mubar(t)$ to $\pi$ where
    \[
    F_\bfp(S_t\sharp \bfmu(t))-F_\bfp(\bfmu(t))\leq  W_1(\mubar(t),\pi)\leq \frac12Ct^2.
    \]
\end{enumerate}
Given these two facts, note that $S_t\sharp \bfmu(t)=(S_t\sharp\mu_1(t),\ldots,S_t\sharp\mu_K(t))$ is feasible, since 
    \[
\sum_{k\in[K]}p_kS_t\sharp\mu_k(t)=S_t\sharp(\sum_{k\in[K]}p_k\mu_k(t))=S_t\sharp\mubar(t)=\pi.
    \]
Combining these two inequalities, we have
\[
F_\bfp(S\sharp \bfmu(t))\leq F_\bfp(\bfmu^*)+At+Ct^2.
\]
Since $A<0$ and $C>0$, by setting $t$ small enough, we show that $F_\bfp(S\sharp\bfmu(t))< F_\bfp(\bfmu^*)$, which contradicts the optimality of $\bfmu^*$. 


Thus, it suffices to prove Claim (i) and (ii). 

\textbf{Claim (i):} $F_\bfp(\bfmu(t))\leq F_\bfp(\bfmu^*)+At+\frac12 Ct^2$.

Observe that
\begin{equation*}
    \begin{split}
        L(\mu_k(t))-L(\mu_k(0))&=\int^t_0 \int_{\R^d}\langle \nabla L(\mu_k(s)), \phi_k\rangle d\mu_k(s,\mathbf x)ds\\
&=\int^t_0 \int_{\R^d}\int_{\R^d}\langle \nabla_1 \ell(\mathbf x,\mathbf z)+\nabla_2 \ell(\mathbf z,\mathbf x), \phi_k(\mathbf x)\rangle d\mu_k(s,\mathbf x)d\mu_k(s,\mathbf z)ds\\
&=\int^t_0 \int_{\R^d}\int_{\R^d}\langle \nabla_x \ell(T_{k,s}^{-1} x,T_{k,s}^{-1}\mathbf z)+\nabla_2 \ell(T_{k,s}^{-1} \mathbf z,T_{k,s}^{-1} \mathbf x), \phi_k(T_{k,s}^{-1} \mathbf x)\rangle d\mu^*_k(\mathbf x)d\mu^*_k(\mathbf z)ds,
    \end{split}
\end{equation*}
where the first equality is due to the fundamental theorem of calculus, the second equality is due to \cref{lem:U}, and the third equality is due to the definition of the pushforward measure. 

Denote
\[
Q_{k,s} \doteq \int_{\R^d}\int_{\R^d}\langle \nabla_1 \ell(T_{k,s}^{-1} \mathbf x ,T_{k,s}^{-1} \mathbf z)+\nabla_2 \ell(T_{k,s}^{-1} \mathbf z,T_{k,s}^{-1} \mathbf x), \phi_k(T_{k,s}^{-1} \mathbf x)\rangle d\mu^*_k(\mathbf x)d\mu^*_k(\mathbf z)
\]

In particular, when $s=0$,
\[
Q_{k,0}=\int_{\R^d}\langle \nabla_1 \ell(\mathbf x,\mathbf z)+\nabla_2 \ell(\mathbf z,\mathbf x), \phi_k(\mathbf x)\rangle d\mu^*_k(\mathbf x)d\mu^*_k(\mathbf z).
\]
Then note that $\|T_{k,s}^{-1} \mathbf x-\mathbf x\|\leq Cs$, we have that
$\|\nabla_2 \ell(T_{k,s}^{-1} \mathbf x,T_{k,s}^{-1} \mathbf z)-\nabla_1 \ell(T_{k,s}^{-1} \mathbf x,T_{k,s}^{-1}\mathbf z)\|\leq Cs (\|\mathbf x\|+\|\mathbf z\|+1), \phi_k(T_{k,s}^{-1} \mathbf x)-\phi_k(\mathbf x)\leq Cs$. Therefore 
\[
Q_{k,s}-Q_{k,0}\leq Cs.
\]
In summary, we have 
\[
L(\mu_k(t))-L(\mu^*_k)-tQ_{k,0}\leq \frac12 Ct^2.
\]
Hence, 
\begin{equation*}
\begin{split}
    F_\bfp(\bfmu(t)) &=\sum_{k\in[K]} p_k L(\mu_k(t))\\
&\leq \sum_{k\in[K]}p_kL(\mu^*_k)+ t\sum_{k\in[K]} p_k Q_{k,0}+\frac12 C t^2\\
&=F_\bfp(\bfmu^*)+At+\frac12 Ct^2.
\end{split}
\end{equation*}

\textbf{Claim (ii): There exists an optimal transport map $S_t$ from $\mubar(t)$ to $\pi$ and \[
F_\bfp(S_t\sharp \bfmu(t))-F_\bfp(\bfmu(t))\leq  W_1(\mubar(t),\pi)\leq \frac12Ct^2.
\] }

 By the duality formula of Kantorovich-Rubinstein distance \citep{villani2009optimal}, there is a Lipschitz continuous function $f$ with Lipschitz constant no larger than $1$ such that 
 \begin{equation*}
    \begin{split}
        W_1(\mubar(t),\pi)&=\sup_{f'\in \text{Lip-1}}\E_{\mubar(t)}[f']-\E_{\pi}[f']=\E_{\mubar(t)}[f]-\E_{\pi}[f]=
     \E_{\mubar(t)}[f]-\E_{\mubar(0)}[f]\\
     &=\int^t_0 \sum_{k\in [K]} p_k \int_{\R^d} \langle \nabla f(\mathbf x), \phi_k(\mathbf x)\rangle d\mu_k(s,\mathbf x)ds,
    \end{split}
 \end{equation*}
 where the last equality is by the definition of the distributional solution of the continuity equation with boundary conditions. 

Define $R_{k,s}\doteq\int\langle\nabla T_{k,s}^{-1}(\mathbf x)\nabla f(T_{k,s}^{-1}\mathbf x), \phi_k(\mathbf x)\rangle d\mu^*_k(\mathbf x) 
$. There exists a constant $C>0$ such that for any $s\in[0,t]$,
\begin{equation*}
    \begin{split}
          \int_{\R^d} \langle \nabla f(\mathbf x), \phi_k(\mathbf x) \rangle d\mu_k(s,\mathbf x) 
&= \int_{\R^d} \langle\nabla f(T_{k,s}^{-1}\mathbf x), \phi_k(T_{k,s}^{-1}\mathbf x) \rangle d\mu^*_k(\mathbf x) \\
&=\int_{\R^d} \langle\nabla f(T_{k,s}^{-1}\mathbf x), \phi_k(T_{k,s}^{-1}x)-\nabla T_{k,s}^{-1}(\mathbf x)\phi_k(\mathbf x)\rangle d\mu^*_k(\mathbf x) +R_{k,s}\\
&\leq \int_{\R^d} \|\phi_k(T_{k,s}^{-1}\mathbf x)-\nabla T_{k,s}^{-1}(\mathbf x)\phi_k(\mathbf x)\| d\mu^*_k(\mathbf x)+R_{k,s}\\
&\leq \int_{\R^d} (\|\phi_k(T_{k,s}^{-1}\mathbf x)-\phi_k(\mathbf x)\|+ \|T_{k,s}^{-1}(\mathbf x)-I\|\|\phi_k(\mathbf x)\|) d\mu^*_k(\mathbf x)+R_{k,s}\\
&\leq C s+R_{k,s}.
    \end{split}
\end{equation*}

Because $\sum_{k\in[K]}p_k\nabla\cdot(\mu^*_k\phi_k)=0$, we have
\begin{equation*}
    \begin{split}
        \sum p_k R_{k,s}&=\int_{\R^d} \sum_{k} p_k\langle\nabla T_{k,s}^{-1}(\mathbf x)\nabla f(T_{k,s}^{-1}\mathbf x), \phi_k(\mathbf x)\rangle d\mu_k(\mathbf x) \\
&=\int_{\R^d} \sum_{k} p_k\langle\nabla g(\mathbf x), \phi_k(\mathbf x)\rangle d\mu^*_k(\mathbf x) =0,
    \end{split}
\end{equation*}
where $g(\mathbf x)=f(T_{k,s}^{-1}\mathbf x)$.

Plugging in these results, we have
\begin{equation*}
    W_1(\mubar(t),\pi)=\E_{\mubar(t)}[f]-\E_{\mubar(0)}[f]=\int^t_0 \sum p_k \int_{\R^d} \langle \nabla f_k(\mathbf x), \phi_k(\mathbf x)\rangle d\mu_k(s,\mathbf x) ds
\leq \int^t_0 Cs ds=\frac12 Ct^2.
\end{equation*}

Since $\mubar(t)$ is absolutely continuous (as we will show later), there exists an optimal transport map $S_t:\R^d\rightarrow\R^d$ from $\mubar(t)$ to $\pi=\mubar$ for the 1-Wasserstein distance $W_1(\mubar(t),\pi)$; that is, 
    \begin{equation*}
        S_t\sharp\mubar(t)=\pi,\quad W_1(\mubar(t),\pi)=\int_{\R^d}\|S_t(\mathbf x)-\mathbf x\|d\mubar(t,\mathbf x).
    \end{equation*}
We show that $\mubar(t)$ is absolutely continuous. Recall $\mu_k(t)=T_{k,t}\sharp\mu^*_k$, where $T_{k,t}:\R^d\rightarrow\R^d$ is a diffeomorphism for small $t$ by \cref{lem:diffeo}. Since $T^{-1}_{k,t}$ is continuously differentiable, $T^{-1}_{k,t}$ satisfies the Luzin N property \citep{evans2018measure}, which claims that, $\lambda(T^{-1}_{k,t}(B))=0$ if $\lambda(B)=0$ for any measurable set $B$, where $\lambda$ is the Lebesgue measure. Hence, for any measurable set $B$ with $0$ measure, $\mu_k(t)(B)=\mu^*_k(T^{-1}_{k,t}(B))=0$, since $\mu^*_k$ is absolutely continuous. 

Recall that, $L(\mu)=\int_{\R^d}\int_{\R^d} \ell(\mathbf x,\mathbf y)d\mu(\mathbf x)d\mu(\mathbf y)$ for some function $\ell:\R^d\times\R^d\rightarrow\R$ such that 
    \begin{equation*}
        |\ell(\mathbf z,\mathbf x)-\ell(\mathbf z,\mathbf y)|\leq\|\mathbf x-\mathbf y\|
    \end{equation*}
for all $\mathbf z \in \R^d$. Define $h(\mathbf y)=\int_{\R^d} \ell(\mathbf x,\mathbf y)d\mu(\mathbf x)$. Then, $h$ is Lipschitz continuous with Lipschitz constant less than one because for every $\mathbf x,\mathbf y\in\R^d$,
    \begin{equation*}
            \left|h(\mathbf x)-h(\mathbf y)\right|=\left|\int_{\R^d} \ell(\mathbf z,\mathbf x)-\ell(\mathbf z,\mathbf y)d\mu(\mathbf z)\right|\leq\int_{\R^d} \left|\ell(\mathbf z,\mathbf x)-\ell(\mathbf z,\mathbf y)\right|d\mu(\mathbf z)\leq\|\mathbf x-\mathbf y\|.
    \end{equation*}
This implies
    \begin{equation*}
        \begin{split}
            F_\bfp(S_t\sharp\bfmu(t))-F_\bfp(\bfmu(t))&=\sum_{k\in[K]}p_k\int_{\R^d} hdS_t\sharp\mu_k(t)-\sum_{k\in[K]}p_k\int_{\R^d} hd\mu_k(t)\\
            &=\int_{\R^d} hd(\sum_{k\in[K]}p_kS_t\sharp\mu_k(t))-\int_{\R^d} hd(\sum_{k\in[K]}p_k\mu_k(t))\\
            &=\int_{\R^d} hd\pi - \int_{\R^d} hd\mubar(t)\\&
            \leq \sup_{f'\in \text{Lip-1}}\E_{\mubar(t)}[f']-\E_{\pi}[f']=W_1(\pi,\mubar)\leq\frac 12Ct^2,
        \end{split}
    \end{equation*}
    where the first equality is due to the definition of objective function.
\end{proof}

\subsection{Proof of \cref{prop:opt-condition-2}}

The proof is similar to that of \cref{prop:opt-condition-1} in the previous section. We need the following lemma.

\begin{lemma}\label{lem:TP2}
  For any $(\bfmu,\bfp)\in\calP_\pi$,
\[
\mathcal{T} \calP_{\pi}(\bfmu,\bfp)\subseteq \mathcal{T}' \calP_{\pi}(\bfmu,\bfp)\doteq\left\{(\bfphi,\bfv):
\sum_{k\in [K]} p_k\nabla \cdot(\mu_k \phi_k)=\sum_{k\in[K]}v_k\mu_k, \sum_{k\in[K]}v_k=0\right\}. 
\]    
\end{lemma}

\begin{proof}
    For any $(\bfphi,\bfv)\in \calT\calP_\pi(\bfmu,\bfp)$, there exists a curve $(\bfmu(t),\bfp(t))\in\calP_\pi$ such that, $\bfmu(0)=\bfmu$, $\bfp(0)=\bfp$, $\bfmu(t)\in\calP_{\pi,\bfp(t)}$, $\frac{d}{dt}\mu_k(t)|_{t=0}=-\nabla\cdot( \mu_k\phi_k)$ and $\frac{d}{dt}p_k(t)|_{t=0}=v_k$ $\forall k\in[K]$. Particularly, we have for all $t\in[0,1]$, $\sum_{k\in[K]}p_k(t)\mu_k(t)=\pi$ in the sense that, for all $f\in C^\infty_c(\R^d)$,
        \begin{equation*}
            \sum_{k\in[K]}p_k(t)\int_{\R^d} f(\mathbf x)d\mu_k(t,\mathbf x)=\int_{\R^d} f(\mathbf x)d\pi(x).
        \end{equation*}
Since the right hand side does not depend on $t$, we have for all $f\in C^\infty_c(\R^d)$,
    \begin{equation*}
        \frac{d}{dt}\sum_{k\in[K]}p_k(t)\int_{\R^d} f(\mathbf x)d\mu_k(t,\mathbf x)=0.
    \end{equation*}
Hence, for all $f\in C^\infty_c(\R^d)$,
    \begin{equation*}
        \begin{split}
            \frac{d}{dt}\sum_{k\in[K]}p_k(t)\int_{\R^d} f(\mathbf x)d\mu_k(t,\mathbf x)|_{t=0}&=\sum_{k\in[K]}(\tfrac{d}{dt}p_k(t))\int_{\R^d} f(\mathbf x)d\mu_k(t,\mathbf x)|_{t=0}+p_k(t)\frac{d}{dt}\int_{\R^d} f(\mathbf x)d\mu_k(t,\mathbf x)|_{t=0}\\
            &=\sum_{k\in[K]}v_k\int_{\R^d} f(\mathbf x)d\mu_k(\mathbf x)+p_k\int_{\R^d} \langle\nabla f(\mathbf x),\phi_k(\mathbf x)\rangle d\mu_k(\mathbf x)=0.
        \end{split}
    \end{equation*}
 where the second equality holds because $\frac{d}{dt}\mu_k(t)|_{t=0}=-\nabla\cdot (\mu_k\phi_k)$ and $\frac{d}{dt}p_k(t)|_{t=0}=v_k$. 
 
 This is the definition of the following equation   \begin{equation*}
        \sum_{k\in[K]}v_k\mu_k=\sum_{k\in[K]}p_k\nabla\cdot(\mu_k\phi_k).
    \end{equation*}

Similarly, since $\sum_{k\in[K]}p_k(t)=1$ for all $t$, 
    \[
    \frac{d}{dt}(\sum_{k\in[K]}p_k(t))=0,
    \]
which implies, at $t=0$,
    \[
    \sum_{k\in[K]}\frac{d}{dt} p_k(0)=\sum_{k\in[K]}v_k=0. \qedhere
    \]
\end{proof}

\begin{proof}[Proof of \cref{prop:opt-condition-2}]
    Consider the following objective function $F(\bfmu,\bfp)=\sum_{k\in[K]}(p_kL(\mu_k)+\frac{\theta}{p_k^\beta})$. The Wasserstein gradient is given by 
       \[
       \nabla F(\bfmu,\bfp)=(\nabla_{\bfmu} F(\bfmu,\bfp),\nabla_\bfp F(\bfmu,\bfp))=(p_1\nabla L(\mu_1),\ldots,p_K\nabla L(\mu_k),L(\mu_1)-\frac{\beta\theta}{p_1^{\beta+1}},\ldots,L(\mu_K)-\frac{\beta\theta}{p_K^{\beta+1}}).
       \]

    For any tangent vector $(\bfphi,\bfv)\in\calT\calP_\pi(\bfmu,\bfp)$, 
\[
\langle\nabla F(\bfmu,\bfp),(\bfphi,\bfv)\rangle_{\bfmu}=\sum_{k\in[K]}p_k\langle\nabla L(\mu_k),\phi_k\rangle_{\mu_k}+\sum_{k\in [K]}v_k(L(\mu_k)-\frac{\beta\theta}{p_k^{\beta+1}}),
\]
    where $\langle\nabla L(\mu_k),\phi_k\rangle_{\mu_k}=\int_{\R^d}\nabla L(\mu_k)(\mathbf x )^\top\phi_k(\mathbf x)d\mu_k(\mathbf x)$.

    For the sake of contradiction, assume $\|\nabla F\|_{\mathcal{T}\calP_{\pi}(\bfmu^*,\bfp^*)}>0$. Then, with the same reason as in the proof of \cref{prop:opt-condition-1}, there exists tangent vector $(\bfphi,\bfv)\in \calT\calP_\pi(\bfmu^*,\bfp^*)$ such that 
        \begin{enumerate}[label=(\roman*)]
            \item $\phi_k\in C^\infty_c$;
            \item $\sum_{k\in[K]}p^*_k\nabla\cdot(\mu^*_k\phi_k)=\sum_{k\in[K]}v_k\mu^*_k$;
            \item $A\doteq\sum_{k\in[K]}v_k( L(\mu^*_k)-\frac{\beta\theta}{(p_k^*)^{\beta+1}})+p^*_k\langle\nabla L(\mu^*_k),\phi_k\rangle_{\mu^*_k}<0$.
        \end{enumerate}

We construct a new feasible solution with this tangent vector as follows. Consider the following time-independent velocity field
    \[
    \phi_k(t)=\phi_k,\quad\forall t\in[0,1].
    \]
Then, the map $T^{\bfphi}_{k,t}=x+\phi_k(x)t$ and $p_k(t)=p^*_k+v_kt$ define a curve $(\bfmu(t),\bfp(t))$ by $\mu_k(t)\doteq T^{\bfphi}_{k,t}\sharp \mu^*_k$ for $t\in[0,1]$. Define $\bar{\mu}(t)=\sum_{k\in[K]}p_k(t)\mu_k(t)$ in the sense that, $\forall t\in[0,1]$, $\bar{\mu}(t)$ is a probability measure in $\calP_2(\R^d)$ such that for all $f\in C^\infty_c(\R^d)$, $\int_{\R^d} f(x)d\bar{\mu}(t)=\sum_{k\in[K]}p_k(t)\int_{\R^d} f(x)d\mu_k(t)$. In particular, $\bar{\mu}(0)=\pi$ since $\sum_{k\in[K]}p_k(0)\mu_k(0)=\pi$. 

Note that the curve $(\bfmu(t),\bfp(t))$ may not be in the feasible set $\calP_{\pi}$, i.e., $\Bar{\mu}(t)$ may not equal to $\pi$. Hence, to get the ``better solution" we are looking for, we need to ``project" this $(\bfmu(t),\bfp(t))$ back to the feasible set $\calP_{\pi}$. 

We show the following facts: there exists a $C>0$ such that for all  $t\in[0,1]$,
\begin{enumerate}[label=(\roman*)]
    \item $F(\bfmu(t),\bfp(t))\leq F(\bfmu^*,\bfp^*)+At+\frac12 Ct^2$;
    \item there exists an optimal transport map $S_t$ from $\mubar(t)$ to $\pi$ and
    \[
    F(S_t\sharp \bfmu(t),\bfp(t))-F(\bfmu(t),\bfp(t))\leq  W_1(\mubar(t),\pi)\leq \frac12Ct^2.
    \]
\end{enumerate}
Assuming that (i) and (ii) are true, note that $S_t\sharp \bfmu(t)=(S_t\sharp\mu_1(t),\ldots,S_t\sharp\mu_K(t))$ is feasible, since 
\[
\sum_{k\in[K]}p_k(t)S_t\sharp\mu_k(t)=S_t\sharp(\sum_{k\in[K]}p_k(t)\mu_k(t))=S_t\sharp\mubar(t)=\pi.
\]
Moreover, since $\sum_{k\in[K]}v_k=0$, we have $\sum_{k\in[K]}p_k(t)=\sum_{k\in[K]}p^*_k+v_k=1$. Since $p^*_k>0$, $p_k(t)>0$ for small $t$. Hence, $\bfp(t)$ is also feasible. 

Combining inequalities in statements (i) and (ii), we have
\[
F(S\sharp \bfmu(t),\bfp(t))\leq F(\bfmu^*,\bfp^*)+At+Ct^2.
\]
Since $A<0$ and $C>0$, by setting small enough $t$, we show that $F(S\sharp\bfmu(t),\bfp(t))< F(\bfmu^*,\bfp^*)$, which contradicts the optimality of $\bfmu^*$.   

We prove Claim (i) and (ii) as follows. 

\textbf{Claim (i): $F(\bfmu(t),\bfp(t))\leq F(\bfmu^*,\bfp^*)+At+\frac12 Ct^2$.}

With the same computation as in \cref{prop:opt-condition-1},
\begin{equation*}
        L(\mu_k(t))-L(\mu_k(0))=\int^t_0 \int_{\R^d}\int_{\R^d}\langle \nabla_1 \ell(T_{k,s}^{-1} \mathbf x,T_{k,s}^{-1} \mathbf z)+\nabla_2 \ell(T_{k,s}^{-1} \mathbf  z,T_{k,s}^{-1} \mathbf x), \phi_k(T_{k,s}^{-1}\mathbf x)\rangle d\mu^*_k(\mathbf x)d\mu^*_k(\mathbf z)ds
\end{equation*}
and we define
\begin{align*}
Q_{k,s}\doteq\int\langle \nabla_1 \ell(T_{k,s}^{-1} \mathbf x,T_{k,s}^{-1} \mathbf z)+\nabla_2 \ell(T_{k,s}^{-1} z,T_{k,s}^{-1} \mathbf x), \phi_k(T_{k,s}^{-1}\mathbf x)\rangle d\mu^*_k(\mathbf x)d\mu^*_k(\mathbf z).
\end{align*}
In particular, when $s=0$,
\[
Q_{k,0}=\int\langle \nabla_1 \ell(\mathbf x,\mathbf z)+\nabla_2 \ell(\mathbf z,\mathbf x), \phi_k(\mathbf x)\rangle d\mu^*_k(\mathbf x)d\mu^*_k(\mathbf z).
\]
Note that $Q_{k,0}=\langle\nabla L(\mu^*_k),\phi_k\rangle_{\mu^*_k}$. Also, note that $\|T_{k,s}^{-1} \mathbf x-\mathbf x\|\leq C_1s$, we have that
$\|\nabla_1 \ell(T_{k,s}^{-1} \mathbf x,T_{k,s}^{-1} \mathbf z)-\nabla_2 \ell(T_{k,s}^{-1} \mathbf x,T_{k,s}^{-1} \mathbf z)\|\leq C_2s (\|\mathbf x\|+\|\mathbf z\|+1), \phi_k(T_{k,s}^{-1} \mathbf x)-\phi_k(\mathbf x)\leq C_3s$. Therefore 
\[
Q_{k,s}-Q_{k,0}\leq C_4s.
\]
In summary, we have 
    \begin{equation}\label{eq:bound1}
        L(\mu_k(t))-L(\mu^*_k)-tQ_{k,0}\leq \frac12 C_5t^2.
    \end{equation}

Also, note that 
    \begin{equation}\label{eq:bound2}
           \frac{\theta}{p^\beta_k(t)}\leq \frac{\theta}{(p_k^*)^\beta}-\frac{\theta\beta}{(p_k^*)^{\beta+1}}v_kt+\frac{1}{2}C_6t^2.
    \end{equation}
This is because $\frac{\theta}{p^\beta_k(t)}$ has Lipschitz continuous derivative on $(0,1)$.

Hence, 
\begin{equation*}
\begin{split}
    F&(\bfmu(t),\bfp(t)) =\sum_{k\in[K]} (p_k(t) L(\mu_k(t))+\frac{\theta}{p_k(t)^\beta})\\
    &\leq \sum_{k\in[K]}p_k(t)L(\mu^*_k)+t\sum_{k\in[K]}p_k(t)Q_{k,0}+\sum_{k\in[K]}(\frac{\theta}{(p_k^*)^\beta}-\frac{\theta\beta}{(p_k^*)^{\beta+1}}v_kt)+\frac{1}{2}C_7t^2 \\
&=\sum_{k\in[K]}(p^*_kL(\mu^*_k) +\frac{\theta}{(p_k^*)^\beta})+ t(\sum_{k\in[K]}v_k(L(\mu^*_k)- \frac{\theta\beta}{(p_k^*)^{\beta+1}})+\sum_{k\in[K]} p^*_k Q_{k,0})+t^2\sum_{k\in[K]}v_kQ_{k,0}+ \frac12 C_7 t^2\\
&= F(\bfmu^*,\bfp^*)+At+\frac12 Ct^2,
\end{split}
\end{equation*}
where the second line (or the first inequality) is because of the bounds \cref{eq:bound1,eq:bound2}, the third line is from the definition $p_k(t)=p_k^*+v_kt$, and the last line is by definition of $A$. 

\textbf{Claim (ii): there exists an optimal transport map $S_t$ from $\mubar(t)$ to $\pi$ and 
\[
F(S_t\sharp \bfmu(t),\bfp(t))-F(\bfmu(t),\bfp(t))\leq  W_1(\mubar(t),\pi)\leq \frac12Ct^2.
\] }

By the duality formula of Kantorovich-Rubinstein distance \citep{villani2009optimal}, there is a Lipschitz continuous function $f$ with Lipschitz constant no larger than $1$ such that 
 \begin{equation*}
    \begin{split}
        W_1(\mubar(t),\pi)&=\sup_{f'\in \text{Lip-1}}\E_{\mubar(t)}[f']-\E_{\pi}[f']=\E_{\mubar(t)}[f]-\E_{\pi}[f]=
     \E_{\mubar(t)}[f]-\E_{\mubar(0)}[f]\\
     &=\int_0^t\sum_{k\in[K]}v_k\int_{\R^d} f(\mathbf x)d\mu_k(s,\mathbf x)ds+\int^t_0 \sum_{k\in [K]} p_k(s) \int_{\R^d} \langle \nabla f(\mathbf x), \phi_k(\mathbf x)\rangle d\mu_k(s,\mathbf x)ds.
    \end{split}
 \end{equation*}

Define $R_{k,s}\doteq\int_{\R^d}\langle\nabla T_{k,s}^{-1}(\mathbf x)\nabla f(T_{k,s}^{-1}\mathbf x), \phi_k(\mathbf x)\rangle d\mu^*_k(\mathbf x) 
$. There exists a constant $C>0$ such that for any $s\in[0,t]$,
\begin{equation*}
    \begin{split}
          \int_{\R^d} \langle \nabla f(\mathbf x), \phi_k(\mathbf x) \rangle d\mu_k(s,\mathbf x) 
&= \int_{\R^d} \langle\nabla f(T_{k,s}^{-1}\mathbf x), \phi_k(T_{k,s}^{-1}\mathbf x) \rangle d\mu^*_k(\mathbf x) \\
&=\int_{\R^d} \langle\nabla f(T_{k,s}^{-1}\mathbf x), \phi_k(T_{k,s}^{-1}\mathbf x)-\nabla T_{k,s}^{-1}(\mathbf x)\phi_k(\mathbf x)\rangle d\mu^*_k(\mathbf x) +R_{k,s}\\
&\leq \int_{\R^d} \|\phi_k(T_{k,s}^{-1}\mathbf x)-\nabla T_{k,s}^{-1}(\mathbf x)\phi_k(\mathbf x)\| d\mu^*_k(x)+R_{k,s}\\
&\leq \int_{\R^d} (\|\phi_k(T_{k,s}^{-1}\mathbf x)-\phi_k(\mathbf x)\|+ \|T_{k,s}^{-1}(\mathbf x)-I\|\|\phi_k(\mathbf x)\|) d\mu^*_k(\mathbf x)+R_{k,s}\\
&\leq C s+R_{k,s}.
    \end{split}
\end{equation*}

Because $\sum_{k\in[K]}p^*_k\nabla\cdot(\mu^*_k\phi_k)=\sum_{k\in[K]}v_k\mu^*_k$, we have
\begin{equation*}
    \begin{split}
    &\sum_{k\in[K]}v_k\int_{\R^d} f(\mathbf x)d\mu_k(s,\mathbf x)+  \sum p_k(s) R_{k,s}\\&=\sum_{k\in[K]}\int_{\R^d} v_kf(T^{-1}_{k,s}(\mathbf x))d\mu^*_k(\mathbf x)+\int_{\R^d} \sum_{k} p_k(s)\langle\nabla T_{k,s}^{-1}(\mathbf x)\nabla f(T_{k,s}^{-1}\mathbf x), \phi_k(\mathbf x)\rangle d\mu^*_k(\mathbf x) \\
&=\int_{\R^d} \sum_{k\in[K]} v_kg(\mathbf x)+p_k(s) \langle\nabla g(\mathbf x),\phi_k(\mathbf x)\rangle d\mu^*_k(\mathbf x)\\
&=\int_{\R^d} \sum_{k\in[K]} v_kg(\mathbf x)+p_k^* \langle\nabla g(\mathbf x),\phi_k(\mathbf x)\rangle d\mu^*_k(\mathbf x)+s\int_{\R^d} v_k\langle\nabla g(\mathbf x),\phi_k(\mathbf x)\rangle d\mu^*_k(\mathbf x)\leq Cs ,
    \end{split}
\end{equation*}
where $g(\mathbf x)=f(T_{k,s}^{-1}\mathbf x)$.

Plugging in these results, we have
\begin{equation*}
    W_1(\mubar(t),\pi)=\E_{\mubar(t)}[f]-\E_{\mubar(0)}[f]=\int^t_0 \sum_{k \in [K]} p_k(s) \int_{\R^d} \langle \nabla f_k(\mathbf x), \phi_k(\mathbf x)\rangle d\mu_k(s,\mathbf x) ds
\leq \int^t_0 Cs ds=\frac12 Ct^2.
\end{equation*}
Since $\mubar(t)$ is absolutely continuous, there exists an optimal transport map $S_t:\R^d\rightarrow\R^d$ from $\mubar(t)$ to $\pi=\mubar$ for the 1-Wasserstein distance $W_1(\mubar(t),\pi)$; that is, 
    \begin{equation*}
        S_t\sharp\mubar(t)=\pi,\quad W_1(\mubar(t),\pi)=\int_{\R^d}\|S_t(\mathbf x)-\mathbf x\|d\mubar(t,\mathbf x).
    \end{equation*}
With the same reason as in the proof of \cref{prop:opt-condition-1}, 
\begin{equation*}
            F(S_t\sharp\bfmu(t),\bfp(t))-F(\bfmu(t),\bfp(t))
            \leq W_1(\pi,\mubar)\leq\frac 12Ct^2. \qedhere
    \end{equation*}
\end{proof}

\section{Proofs of \cref{s:algorithms}}

We will frequently apply the following Gr\"onwall's inequality in the proofs. This inequality is key to obtaining the exponential decrease of the Kullback-Leibler divergence in \cref{s:algorithms}.
\begin{lemma}[Gr\"{o}nwall's inequality]
        Let $\beta$ and $u$ be real-valued continuous functions defined on the closed interval $[a,b]$. If $u$ is differentiable on $(a,b)$ and satisfies $u'(t)\leq\beta(t)u(t)$ for all  $t\in(a,b)$ then $u(t)\leq u(a)e^{\int_a^t\beta(s)ds}$  for all $t\in [a,b]$.
    \end{lemma}

\subsection{Proof of \cref{lem:discrete-iteration}}

\begin{proof}
Denote the optimal value of \cref{eq:discrete-scheme} by $V^*$,i.e., 
     \begin{equation*}
        \begin{split}
         V^*\doteq   \min_{\mathbf x \in \R^n}&\quad f(\mathbf x^\tau_k)+\langle \nabla f(\mathbf x^\tau_k),\mathbf x-\mathbf x^\tau_k\rangle+\frac{1}{2\tau}d(\mathbf x,\mathbf x^\tau_k)^2\\
            \text{s.t.}&\quad \Tilde{g}(\mathbf x)= (1-\alpha h)g(\mathbf x^\tau_k).
        \end{split}
    \end{equation*}
Note that this is a convex optimization problem with one affine inequality constraint. Thus, it satisfies the weak Slater's condition, which implies strong duality; that is, 
\begin{equation*}
    \begin{split}
        V^*=\max_{\lambda\in\R}\min_{\mathbf x \in \R^n} f(\mathbf x^\tau_k)+\langle  \nabla f(\mathbf x^\tau_k)+\lambda \nabla g(\mathbf x^\tau_k),x-\mathbf x^\tau_k\rangle+\frac{1}{2\tau}\|\mathbf x-\mathbf x_k\|^2+\lambda\alpha \tau g(\mathbf x_k^\tau),
    \end{split}
\end{equation*}
where the right-hand side is the optimal value of the Lagrangian dual problem.

Since the objective function of the right-hand side is a quadratic function in $x$, the solution $(\mathbf x^\tau_{k+1},\lambda^*)$ is given by 
 \begin{equation*}
         \mathbf x^\tau_{k+1}=\mathbf x^\tau_k-\tau (\nabla f(\mathbf x^\tau_k)+\lambda^* \nabla g(\mathbf x^\tau_k)),\quad 
\lambda^*=\frac{-\langle \nabla g(\mathbf x^\tau_k), \nabla f(\mathbf x^\tau_k)\rangle +\alpha g(\mathbf x^\tau_k)}{\|\nabla g(\mathbf x^\tau_k)\|^2}. \qedhere
    \end{equation*} 
\end{proof}

\subsection{Proof of \cref{thm:convergence-Euclidean}}

\begin{proof}
\noindent\textbf{(i):}
\begin{equation*}
    \begin{split}
        \tfrac{d}{dt}g(\mathbf x(t))&=\langle \mathbf x'(t),\nabla g(\mathbf x(t))\rangle = \langle\phi(\mathbf x(t)),\nabla g(\mathbf x(t))\rangle\\
        &=-\langle\nabla f(\mathbf x(t)),\nabla g(\mathbf x(t))\rangle - \lambda(t)\|\nabla g(\mathbf x(t))\|^2\\
        &=-\alpha g(\mathbf x(t)).
    \end{split}
\end{equation*}
By Gr\"onwall's inequality, we prove statement (i).



\noindent\textbf{(ii):} We start with the derivative of $f(\bfx(t))$.
\begin{equation}\label{eq:df}
    \begin{split}
    \tfrac{d}{dt}f(\mathbf x(t))&=\langle \mathbf x'(t),\nabla f(\mathbf x(t))\rangle = \langle\phi(\mathbf x(t)),\nabla f(\mathbf x(t))\rangle\\
        &=-\|\nabla f(\mathbf x(t))\|^2-\lambda(t)\langle\nabla g(\mathbf x(t)),\nabla f(\mathbf x(t))\rangle,
    \end{split}
\end{equation}
where the first line is by chain's rule and the second line is due to the definition of $\phi$.

Integrating both parts of \cref{eq:df} and using the fundamental theorem of calculus, we have 
    \begin{equation*}
        f(\mathbf x(T))-f(\mathbf x(0))=\int_0^T\frac{d}{dt}f(\mathbf x(t))dt=-\int_0^T\|\nabla f(\mathbf x(t))+\lambda(t)\nabla g(\mathbf x(t))\|^2+\alpha\lambda(t)g(\mathbf x(t))dt.
    \end{equation*}
This implies
    \begin{equation*}
        \begin{split}
            \int_0^T\|\nabla f(\mathbf x(t))+\lambda(t)\nabla g(\mathbf x(t))\|^2dt &= f(\mathbf x(0))-f(\mathbf x(T))+\int_0^T\alpha\lambda(t)g(\mathbf x(t))dt\\
            &\leq f(\mathbf x(0))-f_{\text{min}}+\int^T_0\frac{|\langle\nabla g(\mathbf x(t)),\nabla f(\mathbf x(t))\rangle|g(\mathbf x(t))}{\|\nabla g(\mathbf x(t))\|^2}+\frac{\alpha g(\mathbf x(t))^2}{\|\nabla g(\mathbf x(t))\|^2}dt\\
            &\leq f(\mathbf x(0))-f_{\text{min}}+\int^T_0L\sqrt{\frac{g(\mathbf x(t))}{\kappa}}+\frac{\alpha g(\mathbf x(t))}{\kappa}dt\\
            &\leq f(\mathbf x(0))-f_{\text{min}}+\frac{2g(\mathbf x(0))}{\kappa}+\frac{L}{\alpha\sqrt{\kappa}}\sqrt{g(\mathbf x(0))},
        \end{split}
    \end{equation*}
where the first inequality is because 
    \[
    \lambda(t)=\frac{-\langle\nabla f(\mathbf x(t)),\nabla g(\mathbf x(t))\rangle+\alpha g(\mathbf x(t))}{\|\nabla g(\mathbf x(t))\|^2}\leq \frac{|\langle\nabla f(\mathbf x(t)),\nabla g(\mathbf x(t))\rangle|+\alpha g(\mathbf x(t))}{\|\nabla g(\mathbf x(t))\|^2};
    \]
the second inequality is by \cref{as:CCGFEuclidean}, and the last inequality is by statement (i), i.e., $g(\mathbf x(t))\leq e^{-\alpha t}g(\mathbf x(0))$.

Since the upper bound we have just shown does not depend on the time index $T$, we have for all $T > 0$ that
    \begin{equation*}
        \min_{t\leq T}\|\nabla f(\mathbf x(t))+\lambda(\mathbf x(t))\nabla g(\mathbf x(t))\|^2\leq\frac{C}{T},
    \end{equation*}
    for some constant $C>0$.
Define $M\doteq\{\mathbf x\in\R^d:g(\mathbf x)=0\}$ and $\calT M(x)\doteq\{\mathbf v\in\R^d:\mathbf v^\top\nabla g(\mathbf x)=0\}$. Then, for all $ \mathbf v\in\calT M(\mathbf x)$, $\mathbf v^\top\nabla f(\mathbf x(t))=\mathbf v^\top(f(\mathbf x(t))+\lambda(t)\nabla g(\mathbf x(t)))$. Hence,
    \begin{equation*}
        \begin{split}
               \min_{t\leq T}\|\nabla f(\mathbf x(t))\|^2_{\cT M_{g(\mathbf x(t))}}&=\min_{t\leq T}\left (\sup_{v\in\calT M(\mathbf x)}\frac{\mathbf v^\top\nabla f(\mathbf x(t))}{\|\mathbf v\|}\right )^2=\min_{t\leq T}\left (\sup_{\mathbf v\in\calT M(\mathbf x)}\frac{v^\top(f(\mathbf x(t))+\lambda(t)\nabla g(\mathbf x(t)))}{\|\mathbf v\|}\right )^2 \\
               &\leq \min_{t\leq T}\|\nabla f(\mathbf x(t))+\lambda(\mathbf x(t))\nabla g(\mathbf x(t))\|^2\leq\frac{C}{T}. 
        \end{split}
    \end{equation*}
\end{proof}

\subsection{Proof of \cref{thm:convergence-1}}

We first show the following direct consequences of \cref{as:WCCGF} that will be useful later. 
\begin{lemma}\label{lem:technical}
    Suppose \cref{as:WCCGF} holds. Then, 
        \begin{enumerate}[label=(\roman*)]
            \item $\|\nabla_{\mu_k}\text{KL}(\mubar\|\pi)\|^2_{\mubar}>0$ if $\text{KL}(\mubar\|\pi)>0$;
            \item $\|L(\mu)\|\leq \ell_{\text{max}}$, $\forall\mu\in\calP_{2,ac}(\R^d)$;
            \item $F_{\bfp}(\bfmu)>-\ell_{\text{max}}$, $\forall\bfmu\in\calP^{\otimes K}_{2,ac}$;
            \item $\|\nabla L(\mu)\|_{\mu}\leq 2L_{\text{max}}\sqrt{K}$, $\forall\mu\in\calP_{2,ac}(\R^d)$;
            \item $\tfrac{\text{KL}(\mubar\|\pi)}{\|\nabla_{\bfmu}\text{KL}(\mubar\|\pi)\|_{\bfmu}}\leq \frac{\sqrt{\text{KL}(\mubar\|\pi)}}{p_{\text{min}}\sqrt{\kappa}}$.
        \end{enumerate}
\end{lemma}

\begin{proof}
\noindent\textbf{(i)} This is because $\|\nabla_{\mu_k}\text{KL}(\mubar\|\pi)\|^2_{\mubar}=p_k^2\|s_{\mubar}-s_\pi\|^2_{\mubar}>0$ if $\mubar\neq\pi$.

\noindent\textbf{(ii)} This is because $|L(\mu)|=|\int_{\R^d}\int_{\R^d}\ell(\mathbf x,\mathbf y)d\mu(\mathbf x)d\mu(\mathbf y)\leq\int_{\R^d}\int_{\R^d}|\ell(\mathbf x,\mathbf y)|d\mu(\mathbf x)d\mu(\mathbf y)\leq \ell_{\text{max}}$ for all $\mu\in\calP_{2,ac}(\R^d)$.

\noindent\textbf{(iii)} $F_{\bfp}(\bfmu)=\sum_{k=1}^Kp_kL(\mu_k)>-\ell_{\text{max}}$.

\noindent\textbf{(iv)} For all $\mathbf x\in\R^d$,
\begin{equation*}
    \begin{split}
        \|\nabla L(\mu)(\mathbf x)\|^2 &= \|\int_{\R^d} \nabla_1\ell(\mathbf x,\mathbf z)+\nabla_2\ell(\mathbf z,\mathbf x)d\mu(\mathbf z)\|^2\\
    &=\sum_{i=1}^{2d}(\int_{\R^d}\nabla_1\ell(\mathbf x,\mathbf z)_i+\nabla_2\ell(\mathbf z,\mathbf x)_id\mu(\mathbf z))^2\\
    &\leq \sum_{i=1}^{2d}\int_{\R^d}(\nabla_1\ell(\mathbf x,\mathbf z)_i+\nabla_2\ell(\mathbf z,\mathbf x)_i)^2d\mu(\mathbf z)\\
    &\leq \sum_{i=1}^{2d}\int(|\nabla_1\ell(\mathbf x,\mathbf z)_i|+|\nabla_2\ell(\mathbf z,\mathbf x)_i|)^2d\mu(z)\\
    &\leq 4L_{\text{max}}^2K.
    \end{split}
\end{equation*}
Hence, $\|\nabla L(\mu)\|^2_\mu=\int_{\R^d}\|\nabla L(\mu)(\mathbf x)\|^2d\mu(\mathbf x)\leq 4L_{\text{max}}^2K$.

\noindent\textbf{(v)}
\begin{equation}
        \begin{split}
            & \frac{\text{KL}(\mubar\|\pi)}{\sqrt{\sum_{k\in[K]}\|p_k(s_{\mubar}-s_\pi)\|^2_{\mu_k}}}\\
            &\leq  \frac{\text{KL}(\mubar(t)\|\pi)}{p_{\text{min}}\sqrt{\sum_{k\in[K]}p_k\|(s_{\mubar}-s_\pi)\|^2_{\mu_k}}}
              \\
    &=  \frac{\text{KL}(\mubar\|\pi)}{p_{\text{min}}\sqrt{\|(s_{\mubar}-s_\pi)\|^2_{\mubar}}}
              \\
              &\leq\tfrac{1}{p_{\text{min}}\sqrt{\kappa}}\sqrt{\text{KL}(\mubar\|\pi)},
        \end{split}
    \end{equation}
    where $p_{\text{min}}=\min_{k\in[K]}p_k$; the equality is because 
\[
\sum_{k\in[K]}p_k\|(s_{\mubar}-s_\pi)\|^2_{\mu_k} = \sum_{k\in[K]}p_k\int_{\R^d}\|(s_{\mubar}-s_\pi)\|^2d\mu_k(\mathbf x)=\int_{\R^d}\|(s_{\mubar}-s_\pi)\|^2d\left(\sum_{k\in[K]}p_k\mu_k(\mathbf x)\right) = \|(s_{\mubar}-s_\pi)\|^2_{\mubar};
\]
and the last line is by the statement (iii) in \cref{as:WCCGF}.
\end{proof}

\begin{proof}[Proof of \cref{thm:convergence-1}]
    \noindent\textbf{(i)} We show that $\frac{d}{dt}\text{KL}(\mubar(t)\|\pi)\leq -\alpha\text{KL}(\mubar(t)\|\pi)$. Then, by Gr\"{o}nwall's inequality, we have 
 \begin{equation*}
        \text{KL}(\mubar(t)\|\pi)\leq e^{-\alpha t}\text{KL}(\mu(0)\|\pi).
    \end{equation*}

    First note that 
        \begin{equation*}
            \begin{split}
                \frac{d}{dt}\text{KL}(\mubar(t)\|\pi)&=\frac{d}{dt}\int_{\R^d}\mubar(t,\mathbf x)\log\frac{\mubar(t,\mathbf x)}{\pi(\mathbf x)}dx\\
                &=\int_{\R^d}\frac{d}{dt}[ \mubar(t,\mathbf x)\log\frac{\mubar(t,\mathbf x)}{\pi(\mathbf x)}]dx\\
                &=\int_{\R^d} [\frac{d}{dt}\mubar(t,\mathbf x)]\cdot\log\frac{\mubar(t,\mathbf x)}{\pi(\mathbf x)}+\mubar(t,\mathbf x)\cdot\frac{d}{dt}[\log\frac{\mubar(t,\mathbf x)}{\pi(\mathbf x)}]dx\\
                &=\int_{\R^d}[\frac{d}{dt}\mubar(t,\mathbf x)]\cdot\log\frac{\mubar(t,\mathbf x)}{\pi(\mathbf x)}+\mubar(t,\mathbf x)\cdot\frac{\frac{d}{dt}\mubar(t,\mathbf x)}{\mubar(t,\mathbf x)}d\mathbf x\\
                &=\int_{\R^d}[-\sum_{k\in[K]}p_k\nabla\cdot (\phi_k\mu_k(t))]\cdot\log\frac{\mubar(t,\mathbf x)}{\pi(\mathbf x)}-\sum_{k\in[K]}p_k\nabla\cdot (\phi_k\mu_k(t))d\mathbf x\\
                &=\sum_{k\in[K]}p_k\int_{\R^d}\langle \nabla\log\frac{\mubar(t,\mathbf x)}{\pi(\mathbf x)},\phi_k(\mathbf x)\rangle d\mu_k(t,\mathbf x)\\
         &=\sum_{k\in[K]}p_k\langle s_{\mubar(t)}-s_\pi,\phi_k\rangle_{\mu_k(t)}\\
         &=-\lambda(\bfmu(t))\sum_{k\in[K]}\|p_k(s_{\mubar(t)}-s_\pi)\|^2_{\mu_k(t)} - p_k\langle p_k(s_{\mubar(t)}-s_\pi),\nabla L(\mu_k(t))\rangle_{\mu_k(t)}\\
         &=-\alpha\text{KL}(\mubar(t)\|\pi).
            \end{split}
        \end{equation*}
        The sixth equality follows by integration by parts and the fact that the integral of the divergence of a vector field vanishing at infinity is 0 by the divergence theorem. The seventh equality is by the definition of $s_{\mubar}$ and $s_\pi$. The eighth equality is due to equation \cref{eq:WCCGF}. We get the last equality by inserting the formula of $\lambda(\bfmu(t))$ into equation \cref{eq:WCCGF}.
    

    \noindent\textbf{(ii)}
    \begin{equation*}
        \begin{split}
            \frac{d}{dt}F_\bfp(\bfmu(t))&= \sum_{k\in[K]}p_k\langle \nabla L(\mu_k(t)),\phi_k\rangle_{\mu_k(t)}\\
            &=\sum_{k\in[K]}p_k\langle\nabla L(\mu_k(t)),-p_k(\nabla L(\mu_k(t))+\lambda(\bfmu(t))(s_{\mubar(t)}-s_\pi))\rangle_{\mu_k(t)}\\
            &=\sum_{k\in[K]}p_k\langle\nabla L(\mu_k(t))+\lambda(\bfmu(t))(s_{\mubar(t)}-s_{\pi}),-p_k(\nabla L(\mu_k(t))+\lambda(\bfmu(t))(s_{\mubar(t)}-s_\pi))\rangle_{\mu_k(t)}\\
            &-\sum_{k\in[K]}p_k\langle \lambda(\bfmu(t))(s_{\mubar(t)}-s_\pi)),-p_k(\nabla L(\mu_k(t))+\lambda(\bfmu(t))(s_{\mubar(t)}-s_\pi))\rangle_{\mu_k(t)}\\
            &=-\sum_{k\in[K]}p^2_k\|\nabla L(\mu_k(t))+\lambda(\bfmu(t))(s_{\mubar(t)}-s_{\pi})\|^2_{\mu_k(t)}-\lambda(\bfmu(t))\sum_{k\in[K]}p_k\langle s_{\mubar(t)}-s_\pi, \phi_k(\bfmu(t))\rangle_{\mu_k(t)}\\
            &=-\sum_{k\in[K]}\|\phi_k(\bfmu(t))\|^2_{\mu_k(t)}-\lambda(\bfmu(t))\sum_{k\in[K]}p_k\langle s_{\mubar(t)}-s_\pi, \phi_k(\bfmu(t))\rangle_{\mu_k(t)}
        \end{split}
    \end{equation*}
The first equality follows from the chain rule. The second equality is by the formula \cref{eq:WCCGF} of $\bfphi$. The last two lines are due to the formula \cref{eq:WCCGF} of $\bfphi$.

By the formula for $\lambda(\bfmu(t))$ in equation \cref{eq:WCCGF}, we have 
    \[
    \sum_{k\in[K]}p_k\langle s_{\mubar(t)}-s_\pi, \phi_k(\bfmu(t))\rangle_{\mu_k(t)} = -\alpha\text{KL}(\mubar(t)\|\pi).
    \]
Hence, 
\[
\frac{d}{dt}F_\bfp(\bfmu(t))=-\sum_{k\in[K]}\|\phi_k(\bfmu(t))\|^2_{\mu_k(t)}+\lambda(\bfmu(t))\alpha\text{KL}(\mubar(t)\|\pi).
\]

Then, using the fundamental theorem of calculus, we have for any $T>0$,
    \begin{equation*}
        \begin{split}
            F_\bfp(\bfmu(T))-F_\bfp(\bfmu(0)) &= \int_0^T\frac{d}{dt}F_\bfp(\bfmu(t))dt \\&= \int_0^T-\sum_{k\in[K]}\|\phi_k(\bfmu(t))\|^2_{\mu_k(t)}+\lambda(\bfmu(t))\alpha\text{KL}(\mubar(t)\|\pi)dt.
        \end{split}
    \end{equation*}
This implies
    \begin{equation*}
        \begin{split}
        \int_0^T\sum_{k\in[K]}\|\phi_k(\bfmu(t))\|^2_{\mu_k(t)}dt&= F_\bfp(\bfmu(0))-F_\bfp(\bfmu(T))+\int^T_0\lambda(\bfmu(t))\alpha\text{KL}(\mubar(t)\|\pi)dt\\
            &\leq F_\bfp(\bfmu(0)) + \ell_{\text{max}}+\int^T_0\lambda(\bfmu(t))\alpha\text{KL}(\mubar(t)\|\pi)dt,
        \end{split}
    \end{equation*}
where the inequality is due to statement (iii) in \cref{lem:technical}, i.e., $F_{\bfp}(\mu)\geq -l_{\text{max}}$.

By inserting the formula for $\lambda(\bfmu(t))$, we have 
    \begin{equation}\label{eq:first-second}
        \begin{split}
            \int^T_0\lambda(\bfmu(t))\alpha\text{KL}(\mubar(t)\|\pi)dt &= -\alpha \int_0^T\frac{\sum_{k\in[K]}\langle p_k\nabla L(\mu_k),p_k(s_{\mubar(t)}-s_\pi)\rangle_{\mu_k(t)}}{\sum_{k\in[K]}\|p_k(s_{\mubar(t)}-s_\pi)\|^2_{\mu_k(t)}}\text{KL}(\mubar(t)\|\pi)dt\\
            &+\alpha^2\int^T_0\frac{\text{KL}(\mubar(t)\|\pi)^2}{\sum_{k\in[K]}\|p_k(s_{\mubar(t)}-s_\pi)\|^2_{\mu_k(t)}}dt
        \end{split}
    \end{equation}

We bound these two terms in this equation separately. For the first term, we have 

\begin{equation}\label{eq:first-term}
        \begin{split}
             & -\alpha \int_0^T\frac{\sum_{k\in[K]}\langle p_k\nabla L(\mu_k),p_k(s_{\mubar(t)}-s_\pi)\rangle_{\mu_k(t)}}{\sum_{k\in[K]}\|p_k(s_{\mubar(t)}-s_\pi)\|^2_{\mu_k(t)}}\text{KL}(\mubar(t)\|\pi)dt \\
             &\leq \alpha\int^T_0 \frac{\sqrt{\sum_{k\in[K]}\|p_k\nabla L(\mu_k(t))\|^2_{\mu_k(t)}}}{\sqrt{\sum_{k\in[K]}\|p_k(s_{\mubar(t)}-s_\pi)\|^2_{\mu_k(t)}}}\text{KL}(\mubar(t)\|\pi)dt,
        \end{split}
    \end{equation}
    by Cauchy-Schwartz inequality.

With statement (iv) in \cref{lem:technical}, we can bound the numerator as 
    \begin{equation}\label{eq:numerator}
          \sqrt{\sum_{k\in[K]}\|p_k\nabla L(\mu_k(t))\|^2_{\mu_k(t)}}\leq 2L_{\text{max}}\sqrt{K}.
    \end{equation}

By statement (v) in \cref{lem:technical}, we can bound the denominator as 
    \begin{equation}\label{eq:deno}
        \begin{split}
            \int^T_0 \frac{1}{\sqrt{\sum_{k\in[K]}\|p_k(s_{\mubar(t)}-s_\pi)\|^2_{\mu_k(t)}}}\text{KL}(\mubar(t)\|\pi)dt
             \leq\tfrac{1}{p_{\text{min}}\sqrt{\kappa}}\int^T_0\sqrt{\text{KL}(\mubar(t)\|\pi)}dt,
        \end{split}
    \end{equation}

Hence, equation \cref{eq:first-term} becomes
\begin{equation*}
    \begin{split}
         & -\alpha \int_0^T\frac{\sum_{k\in[K]}\langle p_k\nabla L(\mu_k),p_k(s_{\mubar(t)}-s_\pi)\rangle_{\mu_k(t)}}{\sum_{k\in[K]}\|p_k(s_{\mubar(t)}-s_\pi)\|^2_{\mu_k(t)}}\text{KL}(\mubar(t)\|\pi)dt \\
             &\leq \alpha\int^T_0 \frac{\sqrt{\sum_{k\in[K]}\|p_k\nabla L(\mu_k(t))\|^2_{\mu_k(t)}}}{\sqrt{\sum_{k\in[K]}\|p_k(s_{\mubar(t)}-s_\pi)\|^2_{\mu_k(t)}}}\text{KL}(\mubar(t)\|\pi)dt\\
             &\leq\frac{\alpha 2L_{\text{max}}\sqrt{K}}{ p_{\text{min}}\sqrt{\kappa}} \int^T_0\sqrt{\text{KL}(\mubar(t)\|\pi)}dt\\
              &\leq \frac{\alpha 2L_{\text{max}}\sqrt{K}}{p_{\text{min}}\sqrt{\kappa}} \int^T_0\sqrt{e^{-\alpha t}\text{KL}(\mubar(0)\|\pi)}dt\\
              &\leq \frac{4\alpha L_{\text{max}}\sqrt{K}}{p_{\text{min}}\sqrt{\kappa}}\sqrt{\text{KL}(\mubar(0)\|\pi)},
    \end{split}
\end{equation*}
where the second inequality is by \cref{eq:numerator,eq:deno}, the third inequality is by statement (i) of \cref{thm:convergence-1}, and the last inequality is by integration. 

Next, we bound the second term in \cref{eq:first-second} as follows. 
    \begin{equation*}
        \begin{split}
              &\alpha^2\int^T_0\frac{\text{KL}(\mubar(t)\|\pi)^2}{\sum_{k\in[K]}\|p_k(s_{\mubar(t)}-s_\pi)\|^2_{\mu_k(t)}}dt\\
              &\leq \alpha^2\int^T_0\frac{\text{KL}(\mubar(t)\|\pi)^2}{p_{\text{min}}\|p_k(s_{\mubar(t)}-s_\pi)\|^2_{\mubar(t)}}dt
              \\&\leq\alpha^2\frac{1}{p_{\text{min}}\kappa} \int^T_0\text{KL}(\mubar(t)\|\pi)dt\\
              &\leq \alpha^2\frac{1}{\kappa}\text{KL}(\mubar(0)\|\pi)\int^T_0e^{-\alpha t}dt\\
              &\leq  \alpha\frac{1}{p_{\text{min}}\kappa}\text{KL}(\mubar(0)\|\pi),
        \end{split}
    \end{equation*}
where the first inequality is obtained by the same way as in \cref{eq:deno}, the second inequality is due to statement (iii) of \cref{as:WCCGF}, and the third inequality is due to part (i) of \cref{thm:convergence-1}.

To conclude, 
\[
\int_0^T\sum_{k\in[K]}\|\phi_k(\bfmu(t))\|^2_{\mu_k(t)}dt\leq F_\bfp(\bfmu(0))+\ell_{\text{max}}+\frac{4\alpha L_{\text{max}}\sqrt{K}}{p_{\text{min}}\sqrt{\kappa}}\sqrt{\text{KL}(\mubar(0)\|\pi)} + \alpha\frac{1}{p_{\text{min}}\kappa}\text{KL}(\mubar(0)\|\pi).
\]
Let $C$ denote the constant $ F_\bfp(\bfmu(0))+\ell_{\text{max}}+\frac{4\alpha L_{\text{max}}\sqrt{K}}{p_{\text{min}}\sqrt{\kappa}}\sqrt{\text{KL}(\mubar(0)\|\pi)} + \alpha\frac{1}{p_{\text{min}}\kappa}\text{KL}(\mubar(0)\|\pi)$.
This implies
    \[
    \min_{t\leq T}\sum_{k\in[K]}\|\phi_k(\bfmu(t))\|^2_{\mu_k(t)}\leq\frac{C}{T}.
    \]
For any unit-norm $v\in\calT\calP_{\pi,\bfp}(\bfmu(t))$, we have 
    \begin{equation*}
        \begin{split}
            \langle\nabla F_\bfp(\bfmu(t)),v\rangle=\sum_{k\in[K]}\langle p_k\nabla L(\mu_k(t)),v_k\rangle_{\mu_k(t)}=&\sum_{k\in[K]}\langle p_k\nabla L(\mu_k(t))+p_k\lambda(\bfmu(t))(s_{\mubar(t)}-s_\pi),v_k\rangle_{\mu_k(t)}\\
            &-\sum_{k\in[K]}\langle p_k\lambda(\bfmu(t))(s_{\mubar(t)}-s_\pi),v_k\rangle_{\mu_k(t)}.
        \end{split}
    \end{equation*}

By \cref{lem:TP}, $\sum_{k\in[K]}p_k\nabla\cdot(v_k\mu_k(t))=0$. Thus,
\begin{equation*}
        \begin{split}
            \sum_{k\in[K]}\langle p_k(s_{\mubar(t)}-s_\pi),v_k\rangle_{\mu_k(t)} &= \sum_{k\in[K]}p_k\int_{\R^d} \langle \nabla\log\frac{\mubar(t,x)}{\pi(x)}, v_k\rangle d\mu_k(t,x)=0,
            \end{split}
    \end{equation*}
where the second equality is due to $\sum_{k\in[K]}p_k\nabla\cdot(v_k\mu_k(t))=0$.

To conclude, 
\begin{equation*}
        \begin{split}
            \langle\nabla F_\bfp(\bfmu(t)),v\rangle=&\sum_{k\in[K]}\langle p_k\nabla L(\mu_k(t))+p_k\lambda(\bfmu(t))(s_{\mubar(t)}-s_\pi),v_k\rangle_{\mu_k(t)}\\
            &\leq \sum_{k\in[K]}\|p_k\nabla L(\mu_k(t))+p_k\lambda(\bfmu(t))(s_{\mubar(t)}-s_\pi)\|_{\mu_k(t)}\|v_k\|_{\mu_k(t)}\\
            &\leq \frac{1}{\sqrt{K}}\sum_{k\in[K]}\|\phi_k(\bfmu(t))\|_{\bfmu(t)}.
        \end{split}
    \end{equation*}
Therefore, 
 \[
 \min_{t\leq T}\|\nabla F_\bfp(\bfmu(t))\|_{\mathcal{T}\calP_{\pi,\bfp}(\bfmu(t))}\leq \frac{C}{\sqrt{KT}}. \qedhere
    \]
\end{proof}

\subsection{Proof of \cref{thm:convergence-2}}

We first show the following bounds that will be useful later. 
\begin{lemma}\label{lem:technical2}
    Suppose \cref{as:WCCGF} holds. Then for all $\bfmu\in\calP_{2,ac}(\R^d)^{\otimes K},\bfp\in\R^{Kd}$,
\[
\|\nabla F(\bfmu,\bfp)\|_{\bfmu,P}\leq 2\sqrt{L_{\text{max}}^2K^2+K(\ell_{\text{max}}+\tfrac{\theta\beta}{p_{\text{min}}^{\beta+1}})^2}.
\]
\end{lemma}
\begin{proof}
By definition, 
    \[
    \|\nabla F(\bfmu,\bfp)\|_{\bfmu,P}^2 = \|\nabla_{\bfmu} F(\bfmu,\bfp)\|_{\bfmu}^2+\|P\nabla_{\bfp} F(\bfmu,\bfp)\|^2.
    \]
    We bound these two terms separately. 

\begin{equation*}
        \|\nabla_{\bfmu} F(\bfmu,\bfp)\|_{\bfmu}^2 =\sum_{k=1}^K\|\nabla_{\mu_k} F(\bfmu,\bfp)\|_{\mu_k}^2=\sum_{k=1}^K\|p_k\nabla L(\mu_k)\|_{\mu_k}^2\leq 4L_{\text{max}}^2K^2,
\end{equation*}
where the last inequality is by statement (iv) in \cref{lem:technical}.

\begin{equation*}
    \begin{split}
           \|P\nabla_{\bfp} F(\bfmu,\bfp)\|^2 &=\sum_{k=1}^K[P(L(\mu_k)-\tfrac{\theta\beta}{p_k^{\beta+1}})]^2=\sum_{k=1}^K(L(\mu_k)-\tfrac{\theta\beta}{p_k^{\beta+1}}-\tfrac{1}{K}\sum_{j=1}^K(L(\mu_j)-\tfrac{\theta\beta}{p_j^{\beta+1}}))^2
           \\ &\leq\sum_{k=1}^K(|L(\mu_k)|+|\tfrac{\theta\beta}{p_k^{\beta+1}}|+\tfrac{1}{K}\sum_{j=1}^K(|L(\mu_j)|+|\tfrac{\theta\beta}{p_j^{\beta+1}}|))^2\\
           &\leq \sum_{k=1}^K(\ell_{\text{max}}+\tfrac{\theta\beta}{p_{\text{min}}^{\beta+1}}+\ell_{\text{max}}+\tfrac{\theta\beta}{p_{\text{min}}^{\beta+1}})^2=4K(\ell_{\text{max}}+\tfrac{\theta\beta}{p_{\text{min}}^{\beta+1}})^2,
    \end{split}
\end{equation*}
where the last line is by \cref{lem:technical}. Combining these two bounds, we prove the lemma. 
\end{proof}

\begin{proof}[Proof of \cref{thm:convergence-2}]
\noindent\textbf{(i):}
For any $k$, 
\begin{equation*}
    \begin{split}
        v_k & = -P(\nabla_{\bfp}F+\lambda\nabla_{\bfp}\text{KL}_k)\\
        &=-[L(\mu_k)-\tfrac{\beta}{p_k^{\beta+1}}+\lambda(\int\mu_k\log\tfrac{\mubar}{\pi}dx+1)-\tfrac{1}{K}\sum_{j}(L(\mu_j)-\tfrac{\beta}{p_j^{\beta+1}}+\lambda(\int\mu_j\log\tfrac{\mubar}{\pi}dx+1))]\\
        &=\tfrac{1}{K}\sum_{j}(L(\mu_j)-L(\mu_k))-\tfrac{1}{K}\sum_j(\tfrac{\beta}{p_j^{\beta+1}}-\tfrac{\beta}{p_k^{\beta+1}})+\tfrac{1}{K}\lambda\sum_j\int(\mu_k-\mu_j)\log\tfrac{\mubar}{\pi}dx).
    \end{split}
\end{equation*}

Now suppose $t_0$ is the first $t$ that  $p_{i}(t)=p_{\min}$ for some $i$. Without loss of generality, assume $i=1$. Then $\tilde{p}(t)=[p_2,\ldots,p_K]\in S:=\{\sum_{j=2}^K p_j=1-p_{\min},p_j\geq p_{\min}\}$, which is a simplex. Observe that 
$f(x)=\sum x_{j}^{-(\beta+1)}$ is a convex function; therefore its maximal value is reached at a vertice of $S$, that is 
\[
\sum_{j=2}^K \frac{1}{p_j^{\beta+1}}\leq 
\frac{1}{(1-(K-1)p_{\min})^{\beta+1}}+\frac{K-2}{p_{\min}^{\beta+1}}.
\]
Next note that
\begin{equation*}
            \frac{d}{dt}\text{KL}(\mubar(t)\|\pi) = \langle\nabla_{\bfmu}\text{KL}(\mubar(t)\|\pi),\bfphi\rangle_{\bfmu}+\langle\nabla_{\bfp}\text{KL}(\mubar(t)\|\pi),\bfv\rangle=  -\alpha\text{KL}(\mubar(t)\|\pi)\leq 0.
    \end{equation*}

Moreover, we note that $\sum p_j \mu_j(x)=\mubar(x)$, so $0\leq \mu_j(x)\leq \frac{1}{p_{\min}}\mubar(x)$. Therefore by the convexity of $g(y)=y\log \frac{y}{\pi(x)}$, $-\frac{\pi(x)}{e}\leq g(y)\leq \max\{\frac{\mubar(x)}{p_{\min}}\log \frac{\mubar(x)}{p_{\min}\pi(x)},0\}$
\[
-\frac{\pi(x)}{e}\leq \mu_j(x)\log \frac{\mu_j(x)}{\pi(x)}\leq \max\{\frac{\mubar(x)}{p_{\min}}\log \frac{\mubar(x)}{p_{\min}\pi(x)},0\}\leq 
\frac{\mubar(x)}{p_{\min}}\log \frac{\mubar(x)}{p_{\min}\pi(x)}+\frac{\pi(x)}{e}. 
\]
\[
 \mu_j(x)\log \frac{\mubar(x)}{\pi(x)}\leq \max\{\frac{\mubar(x)}{p_{\min}}\log \frac{\mubar(x)}{\pi(x)},0\}\leq 
\frac{\mubar(x)}{p_{\min}}\log \frac{\mubar(x)}{\pi(x)}+\frac{\pi(x)}{p_{\min}e}. 
\]
Since $-\frac{\mubar(x)}{p_{\min}}\log \frac{\mubar(x)}{\pi(x)}\leq \frac{\pi(x)}{e p_{\min}}$. Likewise
\[
 \mu_j(x)\log \frac{\mubar(x)}{\pi(x)}\geq \min\{\frac{\mubar(x)}{p_{\min}}\log \frac{\mubar(x)}{\pi(x)},0\}\geq 
-\frac{\pi(x)}{e p_{\min}}. 
\]

Assuming $|L|\leq L_{\max}$, we find
\begin{equation*}
    \begin{split}
        v_{1}  
        &=\tfrac{1}{K}\sum_{j}(L(\mu_j)-L(\mu_{1}))-\tfrac{1}{K}\sum_j(\tfrac{\beta}{p_j^{\beta+1}}-\tfrac{\beta}{p_{1}^{\beta+1}})+\tfrac{1}{K}\lambda\sum_j\int_{\R^d}(\mu_{1}-\mu_j)\log\tfrac{\mubar}{\pi}dx\\
        &\geq -2L_{\max}-\frac{\beta}{K(1-(K-1)p_{\min})^{\beta+1}}+\frac{\beta}{Kp_{\min}^{\beta+1}}+\tfrac{1}{K}\lambda\sum_j\int_{\R^d}(\mu_{1}-\mu_j)\log\tfrac{\mubar}{\pi}dx\\
        &\geq  -2L_{\max}-\frac{\beta}{K(1-(K-1)p_{\min})^{\beta+1}}+\frac{\beta}{Kp_{\min}^{\beta+1}}-\tfrac{1}{K}\lambda\sum_j\int_{\R^d} (\frac{\mubar(x)}{p_{\min}}\log\tfrac{\mubar}{\pi}+\frac{\pi(x)}{e p_{\min}})dx\\
        &\geq -C-C (1+KL(\mubar(t)\|\pi))/p_{\min}+\frac{\beta}{K p_{\min}^{\beta+1}}.
    \end{split}
\end{equation*}
Hence, if $p_{\min}$ is small enough, we have $v_1>0$. This implies $p_1(t)<p_{\min}$ for some $t<t_0$. This contradicts the fact that $t_0$ is the first time that $p_i(t)=p_{\text{min}}$ for some $i\in[K]$.

To show $\sum_{k=1}^Kp_k(t)=1$ for any $t>0$, we note that 
\[
\sum_{k=1}^K\int^T_0\tfrac{d}{dt}p_k(t)dt = \int^T_0\sum_{k=1}^Kv_k(\bfmu(t),\bfp(t))dt = 0,
\]
where the first equality is by the construction of $v_k$ and the second equality is due to the projection operator $P$ we defined. Hence, 
\[
\sum_{k=1}p_k(T)-\sum_{k=1}^Kp_k(0)=\sum_{k=1}^K\int^T_0\tfrac{d}{dt}p_k(t)dt=0
\]
for any $T>0$. This implies $\sum_{k=1}p_k(T)=1$ as long as $\sum_{k=1}p_k(0)=1$.

 \noindent\textbf{(ii):} Note that 
    \begin{equation*}
            \frac{d}{dt}\text{KL}(\mubar(t)\|\pi) = \langle\nabla_{\bfmu}\text{KL}(\mubar(t)\|\pi),\bfphi\rangle_{\bfmu}+\langle\nabla_{\bfp}\text{KL}(\mubar(t)\|\pi),\bfv\rangle=  -\alpha\text{KL}(\mubar(t)\|\pi),
    \end{equation*}
    where the first equality is by chain's rule and the second one is due to the formula of $\bfphi,\bfv$ in equation \cref{eq:WCCGF2}. Then, by the Gr\"onwall's inequality, we proved statement (ii).

\noindent\textbf{(iii):}
For abrivation, in this proof, we use $F(t)$, $\text{KL}(t)$, and $\lambda(t)$ to denote $F(\bfmu(t),\bfp(t))$, $\text{KL}(\mubar(t)\|\pi)$, and $\lambda(\bfmu(t),\bfp(t))$, repsectively. In this proof, with a little abuse of notation, we define the inner product $\langle \nabla F(\bfmu,\bfp),\nabla\text{KL}(\mubar\|\pi)\rangle$ as 
\[\langle\nabla_{\bfmu}F(\bfmu,\bfp),\nabla_{\bfmu}\text{KL}(\mubar\|\pi))\rangle_{\bfmu}+\langle\nabla_{\bfp}F(\bfmu,\bfp),\nabla_{\bfp}\text{KL}(\mubar\|\pi)\rangle,\]
where 
    \[
    \langle\nabla_{\bfmu}F(\bfmu,\bfp),\nabla_{\bfmu}\text{KL}(\mubar\|\pi))\rangle_{\bfmu} = \sum_{k\in[K]}\int_{\R^d}\langle\nabla_{\mu_k} F(\bfmu,\bfp),\nabla_{\mu_k}\text{KL}(\mubar\|\pi))\rangle d\mu_k.
    \]

Now we calculate $\frac{d}{dt}F(t)$ as follows.
    \begin{equation*}
        \begin{split}
            \frac{d}{dt}F(t) &= \langle\nabla F(t),(\bfphi,\bfv)\rangle \\
            &= -\langle\nabla F(t), \nabla F(t)+\lambda(t)\nabla\text{KL}(t)\rangle + \eta\langle\nabla_{\bfp}F(t),\mathbf{1}\rangle\\
            &= -\langle\nabla F(t), \nabla F(t)+\lambda(t)\nabla\text{KL}(t)\rangle + \frac{1}{K} \langle\nabla_{\bfp}F(t)+\lambda(t)\nabla_{\bfp}\text{KL}(t), \mathbf{1}\rangle \langle\nabla_{\bfp}F(t),\mathbf{1}\rangle\\
           &=- \langle\nabla_{\bfmu} F(t), \nabla_{\bfmu} F(t)+\lambda(t)\nabla_{\bfmu}\text{KL}(t)\rangle_{\bfmu(t)}\\ 
           &- (\langle\nabla_{\bfp} F(t), \nabla_{\bfp} F(t)+\lambda(t)\nabla_{\bfp}\text{KL}(t)\rangle - \frac{1}{K} \langle\nabla_{\bfp}F(t)+\lambda(t)\nabla_{\bfp}\text{KL}(t), \mathbf{1}\rangle\langle\nabla_{\bfp}F(t),\mathbf{1}\rangle )\\
           &=- (\langle\nabla_{\bfmu} F(t), \nabla_{\bfmu} F(t)+\lambda(t)\nabla_{\bfmu}\text{KL}(t)\rangle_{\bfmu(t)} +  \langle P\nabla_{\bfp} F(t), P(\nabla_{\bfp} F(t)+\lambda(t)\nabla_{\bfp}\text{KL}(t))\rangle)\\
           &=-\langle \nabla F(t), \nabla F(t)+\lambda(t)\nabla\text{KL}(t)\rangle_{\bfmu(t),P} \\
           &=-\|\nabla F(t)+\lambda(t)\nabla\text{KL}(t)\|^2_{\bfmu(t),P} +\lambda(t)\langle\nabla\text{KL}(t), \nabla F(t)+\lambda(t)\nabla\text{KL}(t)\rangle_{\bfmu(t),P}.
        \end{split}
    \end{equation*}
    By inserting the formula of $\lambda(t)$ given in equation \cref{eq:WCCGF2}, we observe that 
        \[
    \lambda(t)\langle\nabla\text{KL}(t), \nabla F(t)+\lambda(t)\nabla\text{KL}(t)\rangle_{\bfmu(t),P}=\alpha\lambda(t)\text{KL}(t).
        \]
Hence, we have 
    \[
    \frac{d}{dt}F(t) = -\|\nabla F(t)+\lambda(t)\nabla\text{KL}(t)\|^2_{\bfmu(t),P}
           +\alpha\lambda(t)\text{KL}(t).
    \]
By the fundamental theorem of calculus, 
    \begin{equation*}
            \int_0^T\|\nabla F(t)+\lambda(t)\nabla\text{KL}(t)\|^2_{\bfmu(t),P}dt
            = F(0)-F(T)+\int_0^T\alpha\lambda(t)\text{KL}(t)dt.
    \end{equation*}
    
To bound $\int_0^T\|\nabla F(t)+\lambda(t)\nabla\text{KL}(t)\|^2_{\bfmu(t),P}dt$, we show the following bound for $\lambda(t)\text{KL}(t)$:
    \begin{equation*}
        \begin{split}
            \lambda(t)\text{KL}(t)&= - \frac{\langle\nabla F(t),\nabla\text{KL}(t)\rangle_{\bfmu(t),P}}{\|\nabla\text{KL}(t)\|^2_{\bfmu(t),P}}\text{KL}(t) + \frac{\alpha\text{KL}(t)^2}{\|\nabla\text{KL}(t)\|^2_{\bfmu(t),P}}\\
            &\leq  \frac{|\langle\nabla F(t),\nabla\text{KL}(t)\rangle_{\bfmu(t),P}|}{\|\nabla\text{KL}(t)\|^2_{\bfmu(t),P}}\text{KL}(t) + \frac{\alpha\text{KL}(t)^2}{\|\nabla\text{KL}(t)\|^2_{\bfmu(t),P}}\\
            &\leq \frac{\|\nabla F(t)\|_{\bfmu(t),P}}{\|\nabla\text{KL}(t)\|_{\bfmu(t),P}}\text{KL}(t) + \frac{\alpha\text{KL}(t)^2}{\|\nabla\text{KL}(t)\|^2_{\bfmu(t),P}}\\
            &\leq  \frac{\|\nabla F(t)\|_{\bfmu(t),P}}{\|\nabla_{\bfmu}\text{KL}(t)\|_{\bfmu(t)}}\text{KL}(t) + \frac{\alpha\text{KL}(t)^2}{\|\nabla_{\bfmu}\text{KL}(t)\|^2_{\bfmu(t)}}
            \\
            &\leq \frac{1}{\underline{p}\sqrt{\kappa}}\sqrt{\text{KL}(t)}\|\nabla F(t)\|_{\bfmu(t),P} + \frac{\alpha}{\underline{p}^2\kappa}\text{KL}(t),
        \end{split}
    \end{equation*}
where the third line is due to Cauchy-Schwatz inequality, the fourth line is because $\|\nabla\text{KL}(t)\|_{\bfmu(t),P}\geq \|\nabla_{\bfmu}\text{KL}(t)\|_{\bfmu(t)}$, and the last line is due to statement (v) in \cref{lem:technical} and statement (i) in \cref{thm:convergence-2} (i.e., $p_k(t)\geq \underline{p}$ $\forall k\in[K], t\geq 0$).

By part (ii), $\text{KL}(\mubar(t)\|\pi)= e^{-\alpha t}\text{KL}(\mubar(0)\|\pi)$. Hence,
    \begin{equation*}
\lambda(t)\text{KL}(t)\leq \frac{1}{\underline{p}\sqrt{\kappa}}\sqrt{\text{KL}(0)}\|\nabla F(t)\|_{\bfmu(t),P}e^{-\frac{1}{2}\alpha t} + \frac{\alpha}{\underline{p}^2\kappa}\text{KL}(0)e^{-\alpha t}.
    \end{equation*}
By \cref{lem:technical2}, we have 
 \begin{equation*}
\lambda(t)\text{KL}(t)\leq \frac{1}{\underline{p}\sqrt{\kappa}}\sqrt{\text{KL}(0)} 2\sqrt{L_{\text{max}}^2K^2+K(\ell_{\text{max}}+\tfrac{\theta\beta}{\underline{p}^{\beta+1}})^2}e^{-\frac{1}{2}\alpha t} + \frac{\alpha}{\underline{p}^2\kappa}\text{KL}(0)e^{-\alpha t}.
    \end{equation*}
This implies
    \[
    \int^T_0\lambda(t)\text{KL}(t)dt\leq \frac{4}{\alpha\underline{p}\sqrt{\kappa}}\sqrt{\text{KL}(0)}\sqrt{L_{\text{max}}^2K^2+K(\ell_{\text{max}}+\tfrac{\theta\beta}{\underline{p}^{\beta+1}})^2} + \frac{1}{\underline{p}\kappa}\text{KL}(0).
    \]
Therefore, 
    \[
    \int_0^T\|\nabla F(t)+\lambda(t)\nabla\text{KL}(t)\|^2_{\bfmu(t),P}dt\leq F(0)+\ell_{\text{max}}+\frac{4}{\alpha\underline{p}\sqrt{\kappa}}\sqrt{\text{KL}(0)}\sqrt{L_{\text{max}}^2K^2+K(\ell_{\text{max}}+\tfrac{\theta\beta}{\underline{p}^{\beta+1}})^2} + \frac{1}{\underline{p}\kappa}\text{KL}(0).
    \]
    Note that the right-hand side is a constant that does not depend on $T$. Denote the upper bound by $C$. This implies
        \[
      \min_{t\leq T}  \|\nabla F(t)+\lambda(t)\nabla\text{KL}(t)\|^2_{\bfmu(t),P}\leq\frac{C}{T}.
        \]

Finally, we will show our statement (iii) as follows. Recall $\|\nabla F(\bfmu^*,\bfp^*)\|_{\mathcal{T}\calP_{\pi}(\bfmu^*,\bfp^*)}$ is defined as 
    \[
    \|\nabla F(\bfmu(t),\bfp(t))\|_{\mathcal{T}\calP_{\pi}(\bfmu(t),\bfp(t))} = \sup_{(\bfphi,\bfv)\in\calT\calP_\pi(\bfmu(t),\bfp(t))}\frac{\langle\nabla_{\bfmu}F(t),\bfphi\rangle_{\bfmu}+\langle\nabla_{\bfp}F(t),\bfv\rangle}{\|\bfphi\|^2_{\bfmu}+\|\bfv\|^2},
    \]
 and the vector $\nabla F(t)+\lambda(t)\nabla\text{KL}(t)$ is given by 
    \begin{equation*}
        \begin{split}
                   (&p_1(\nabla L(\mu_1(t))+\lambda(s_{\mubar(t)}-s_\pi),\ldots,p_K(\nabla L(\mu_K(t))+\lambda(s_{\mubar(t)}-s_\pi),\\ &\nabla_{\bfp}F(t)_1+\lambda(t)\langle\log\frac{\mubar(t)}{\pi}+1,\mu_1(t)\rangle,\ldots, \nabla_{\bfp}F(t)_K+\lambda(t)\langle\log\frac{\mubar(t)}{\pi}+1,\mu_K(t)\rangle).
        \end{split}
    \end{equation*}

For any pair $(\bfphi,\bfv)\in\calT\calP_\pi(\bfmu,\bfp)$, by \cref{lem:TP2}, we have 
    \[
    \sum_{k\in [K]} p_k\nabla \cdot(\mu_k \phi_k)=\sum_{k\in[K]}v_k\mu_k,\quad \sum_{k\in[K]}v_k=0.
    \]
Thus, 
    \begin{equation*}
        \begin{split}
\langle\nabla_{\bfmu}F(t),\bfphi\rangle_{\bfmu}+\langle\nabla_{\bfp}F(t),\bfv\rangle &= \langle\nabla_{\bfmu}F(t)+\lambda(t)\nabla_{\bfmu}\text{KL}(t),\bfphi\rangle_{\bfmu(t)}+\langle\nabla_{\bfp}F(t)+\lambda(t)\nabla_{\bfp}\text{KL}(t),\bfv\rangle\\
&- (\langle\lambda(t)\nabla_{\bfmu}\text{KL}(t),\bfphi \rangle_{\bfmu(t)}+ \langle\lambda(t)\nabla_{\bfp}\text{KL}(t),\bfv\rangle).
        \end{split}
    \end{equation*}
Note that 
    \[
    \langle\nabla_{\bfmu}F(t)+\lambda(t)\nabla_{\bfmu}\text{KL}(t),\bfphi\rangle_{\bfmu(t)}+\langle\nabla_{\bfp}F(t)+\lambda(t)\nabla_{\bfp}\text{KL}(t),\bfv\rangle = \langle\nabla F(t)+\lambda(t)\nabla\text{KL}(t),(\bfphi,\bfv)\rangle_{\bfmu(t),P}.
    \]
To see this, 
    \begin{equation*}
        \begin{split}
            \langle P\nabla_{\bfp}F(t)+\lambda(t)\nabla_{\bfp}\text{KL}(t),P\bfv\rangle &= \langle\nabla_{\bfp}F(t)+\lambda(t)\nabla_{\bfp}\text{KL}(t),\bfv\rangle - \frac{1}{K}\langle \nabla_{\bfp}F(t)+\lambda(t)\nabla_{\bfp}\text{KL}(t),\mathbf{1}\rangle\langle\bfv,\mathbf{1}\rangle\\
            &= \langle\nabla_{\bfp}F(t)+\lambda(t)\nabla_{\bfp}\text{KL}(t),\bfv\rangle,
        \end{split}
    \end{equation*}
because $\langle\bfv,\mathbf{1}\rangle=0$. Hence,
    \begin{equation*}
        \begin{split}
           &\langle\nabla_{\bfmu}F(t)+\lambda(t)\nabla_{\bfmu}\text{KL}(t),\bfphi\rangle_{\bfmu(t)}+\langle\nabla_{\bfp}F(t)+\lambda(t)\nabla_{\bfp}\text{KL}(t),\bfv\rangle \\&= \langle\nabla_{\bfmu}F(t)+\lambda(t)\nabla_{\bfmu}\text{KL}(t),\bfphi\rangle_{\bfmu(t)} +  \langle P\nabla_{\bfp}F(t)+\lambda(t)\nabla_{\bfp}\text{KL}(t),P\bfv\rangle\\
           &= \langle\nabla F(t)+\lambda(t)\nabla\text{KL}(t),(\bfphi,\bfv)\rangle_{\bfmu(t),P}.
        \end{split}
    \end{equation*}
    
Moreover, 
    \[
    \langle\lambda(t)\nabla_{\bfmu}\text{KL}(t),\bfphi \rangle_{\bfmu(t)}+ \langle\lambda(t)\nabla_{\bfp}\text{KL}(t),\bfv\rangle=0.
    \]
To see this,
    \begin{equation*}
        \begin{split}
            \langle\nabla_{\bfmu}\text{KL}(t),\bfphi \rangle_{\bfmu(t)} &= \sum_{k\in[K]} \langle\nabla_{\mu_k}\text{KL}(t),\phi_k \rangle_{\mu_k(t)} 
            =\sum_{k\in[K]} p_k\langle s_{\mubar(t)}-s_\pi,\phi_k \rangle_{\mu_k(t)}\\
            &=-\sum_{k\in[K]}v_k\langle\log\frac{\mubar(t)}{\pi}+1,\mu_k(t)\rangle
            = - \langle\nabla_{\bfp}\text{KL}(t),\bfv\rangle,
        \end{split}
    \end{equation*}
    where the first equality is the definition for the inner product, the second equality is due to the formula of $ \nabla_{\bfmu}\text{KL}(t)$, the third equality is because $\sum_{k\in [K]} p_k(t)\nabla \cdot(\mu_k(t) \phi_k)=\sum_{k\in[K]}v_k\mu_k(t)$, and the last equality is due to the definition of inner product and the formula of $\nabla_{\bfp}\text{KL}(t)$.

    To conclude, we get
        \begin{equation*}
            \begin{split}
                \langle\nabla_{\bfmu}F(t),\bfphi\rangle_{\bfmu}+\langle\nabla_{\bfp}F(t),\bfv\rangle &=\langle\nabla F(t)+\lambda(t)\nabla\text{KL}(t),(\bfphi,\bfv)\rangle_{\bfmu(t),P}\\
                &\leq \|\nabla F(t)+\lambda(t)\nabla\text{KL}(t)\|_{\bfmu(t),P}\|(\bfphi,\bfv)\|_{\bfmu(t),P},
            \end{split}
        \end{equation*}
        where the second inequality is by Cauchy-Schwartz inequality. Also, note that 
        \begin{equation*}
                \|(\bfphi,\bfv)\|_{\bfmu(t),P} = \|\bfphi\|^2_{\bfmu}+\|P\bfv\|^2=\|\bfphi\|^2_{\bfmu}+\|\bfv\|^2,
        \end{equation*}
        since $P\bfv = \bfv$ by the fact $\sum_{k\in[K]}v_k=0$.

Hence, for any pair $(\bfphi,\bfv)\in\calT\calP_\pi(\bfmu^*,\bfp^*)$,
    \begin{equation*}
        \begin{split}
\frac{\langle\nabla_{\bfmu}F(t),\bfphi\rangle_{\bfmu}+\langle\nabla_{\bfp}F(t),\bfv\rangle}{\|\bfphi\|^2_{\bfmu}+\|\bfv\|^2}\leq \|\nabla F(t)+\lambda(t)\nabla\text{KL}(t)\|_{\bfmu(t),P},
        \end{split}
    \end{equation*}
which implies 
    \[
     \|\nabla F(\bfmu(t),\bfp(t))\|_{\mathcal{T}\calP_{\pi}(\bfmu(t),\bfp(t))}\leq \|\nabla F(t)+\lambda(t)\nabla\text{KL}(t)\|_{\bfmu(t),P},
    \]
and thus,
 \[
     \min_{t\leq T}\|\nabla F(\bfmu(t),\bfp(t))\|_{\mathcal{T}\calP_{\pi}(\bfmu(t),\bfp(t))}\leq \min_{t\leq T}\|\nabla F(t)+\lambda(t)\nabla\text{KL}(t)\|_{\bfmu(t),P}\leq \frac{C}{\sqrt{T}}. \qedhere
    \]
\end{proof}






\end{appendices}

\end{document}